\documentclass[a4paper,reqno]{amsart}
\usepackage{amsmath,amssymb}
\usepackage{hyperref}
\usepackage{enumerate}
\usepackage{xspace}
\usepackage{boxedminipage}
\usepackage{algpseudocode}
\usepackage{graphicx,xcolor}
\usepackage{epsfig}
\usepackage{tikz}
\usepackage{rotating}
\usepackage{ifthen}
\usepackage{amsmath,amssymb}    
\usepackage{stmaryrd} 
\usepackage{eufrak}     
\usepackage{mathrsfs} 
\usepackage{bbold}    
\provideboolean{usemathrsfs}
\setboolean{usemathrsfs}{true}
\usepackage{enumerate}
\usepackage{xspace}
\usepackage{boxedminipage}
\usepackage[latin1]{inputenc}
\usepackage{verbatim}
\provideboolean{usebeamer}
\provideboolean{isthesis}
\provideboolean{isamsltex}
\provideboolean{issiamltex}

\newcommand{\Ignore}[1]{}
\newcommand{\freeze}[1]{}
\newcommand{\crossout}[1]{{\textcolor{red!20}{#1}}}
\newcommand{\highlight}[1]{{\textcolor{blue}{#1}}}
\newcommand{\standout}[1]{{\textcolor{magenta}{#1}}}
\provideboolean{shownotes}
\setboolean{shownotes}{false}
\newcounter{margnote}[page]
\newcommand{\margnotemark}{\standout{\upshape\texttt{>\arabic{margnote}<}}}
\newcommand{\margnote}[2][]{\ifthenelse{\boolean{shownotes}}%
  {\stepcounter{margnote}%
    \margnotemark%
    \marginpar{\texttt{\raggedright\tiny\margnotemark{#1}: #2}}}
  {}}
\provideboolean{showtodo}

\newcommand{\todo}[1]{
  \ifthenelse{\boolean{showtodo}}{\margnote[To do.]{#1}}{}}

\provideboolean{showcomments}
\newcommand{\margincomment}[1]{
  \ifthenelse{\boolean{showcomments}}{\marginpar{\tiny #1}}{}
}
\provideboolean{showchanges}
\setboolean{showchanges}{false}
\newcommand{\changes}[1]{
  \ifthenelse{\boolean{showchanges}} {{\highlight{#1}}} {#1}
}
\newcommand{\changefromto}[3][replace]{
  \ifthenelse{\boolean{showchanges}} {{\crossout{#2}\margnote{#1}}{\highlight{#3}}} {#3\xspace}
}

\newcommand{\ie}{\ensuremath{\text{ i.e., }}\xspace}

\newcommand{\aposteriori}{{aposteriori}\xspace}

\usepackage{mathrsfs}
\provideboolean{usemathrsfs}
\setboolean{usemathrsfs}{true}
\newcommand{\mathscript}{\mathscr}
	   
\newcommand{\cA}{\ensuremath{\mathscript A}\xspace}

\newcommand{\cC}{\ensuremath{\mathscript C}\xspace}

\newcommand{\cE}{\ensuremath{\mathscript E}\xspace}

\newcommand{\cR}{\ensuremath{\mathscript R}\xspace}

\newcommand{\cT}{\ensuremath{\mathscript T}\xspace}

\newcommand{\rN}{\ensuremath{\mathbb N}\xspace}

\newcommand{\rP}{\ensuremath{\mathbb P}\xspace}

\newcommand{\rR}{\ensuremath{\mathbb R}\xspace}

\newcommand{\rV}{\ensuremath{\mathbb V}\xspace}

\newcommand{\naturals}{\rN\xspace}

\newcommand{\reals}{\rR}
\newcommand{\R}[1]{\reals^{#1}}

\newcommand{\Tolto}{\smallsetminus}
\newcommand{\take}{\Tolto}

\newcommand{\closure}[1]{\overline{#1}}

\newcommand{\Th}{\ensuremath{\varTheta}\xspace}
\newcommand{\La}{\ensuremath{\varLambda}\xspace}

\newcommand{\Pg}{\ensuremath{\varPi}\xspace}

\newcommand{\Fi}{\ensuremath{\varPhi}\xspace}
\newcommand{\Ps}{\ensuremath{\varPsi}\xspace}
\newcommand{\W}{\ensuremath{\varOmega}\xspace}

\renewcommand{\a}{\ensuremath{\alpha}\xspace}
\renewcommand{\b}{\ensuremath{\beta}\xspace}
\newcommand{\g}{\ensuremath{\gamma}\xspace}

\newcommand{\ep}{\ensuremath{\varepsilon}\xspace}
\newcommand{\epsi}{\ensuremath{\epsilon}\xspace}

\newcommand{\la}{\ensuremath{\lambda}\xspace}

\newcommand{\w}{\ensuremath{\omega}\xspace}

\newcommand{\qp}[1]{\ensuremath{\!\left({#1}\right)}}

\newcommand{\qpbig}[1]{\ensuremath{\big(#1\big)}}
\newcommand{\qpBig}[1]{\ensuremath{\Big(#1\Big)}}
\newcommand{\qpbigg}[1]{\ensuremath{\bigg(#1\bigg)}}

\newcommand{\qb}[1]{\ensuremath{\!\left[{#1}\right]}}

\newcommand{\opinter}[2]{\ensuremath{\left(#1,#2\right)}\xspace}
\newcommand{\clinter}[2]{\ensuremath{\left[#1,#2\right]}\xspace}
\newcommand{\opclinter}[2]{\ensuremath{\left(#1,#2\right]}\xspace}

\newcommand{\powqp}[2]{\ensuremath{\qp{#2}^{\kern -.2em\lower .3ex\hbox{\scriptsize $#1$}}\kern-.3em}}
\newcommand{\psqrt}[1]{\powqp{1/2}{#1}}
\newcommand{\norm}[1]{\ensuremath{\left|#1\right|}}

\newcommand{\abs}[1]{\ensuremath{\left|#1\right|}}
\newcommand{\Norm}[1]{\ensuremath{\left\|#1\right\|}}

\newcommand{\Normbig}[1]{\ensuremath{\big\|#1\big\|}}

\newcommand{\ltwop}[2]{\ensuremath{\qa{#1,#2}}}

\newcommand{\ltwopbig}[2]{\ensuremath{\qabig{#1,#2}}}

\newcommand{\duality}[2]{\ensuremath{\left\langle #1\,\vert\,#2\right\rangle}}

\newcommand{\ensemble}[2]{\ensuremath{\left\{ #1:\;#2 \right\}}}

\newcommand{\setof}[1]{\qc{#1}}

\newcommand{\vecof}[1]{\begin{pmatrix}{#1}\end{pmatrix}}
\newcommand{\matof}[1]{\vecof{#1}}

\newcommand{\seqof}[1]{\qp{#1}}
\newcommand{\seq}[1]{\seqof{#1}}
\newcommand{\seqi}[2]{\seq{#1_{#2}}}
\newcommand{\sequ}[2]{\seq{#1^{#2}}}

\newcommand{\seqsu}[3]{{\sequ{#1}{#2}}_{#2\in#3}}

\newcommand{\seqsifromto}[4]{\seqi{#1}{#2}_{\rangefromto{#2}{#3}{#4}}}

\newcommand{\sumsu}[3]{\ensuremath{\sum_{#2\in#3}{#1}^{#2}}}
\newcommand{\sumifromto}[3]{\ensuremath{\sum_{#1=#2}^{#3}}}

\newcommand{\fromto}[2]{\ensuremath{\left[#1:#2\right]}}

\newcommand{\integerbetween}[2]{\ensuremath{\in\fromto{#1}{#2}}}
\newcommand{\rangefromto}[3]{\ensuremath{#1\integerbetween{#2}{#3}}}

\renewcommand{\d}{\ensuremath{\,\mathrm{d}}}

\newcommand{\AND}{\ensuremath{\text{ and }}}

\newcommand{\OR}{\ensuremath{\text{ or }}}

\newcommand{\constant}[1]{C_{#1}}
\newcommand{\constext}[1]{\constant{\text{#1}}}            

\renewcommand{\div}{\operatorname{div}}

\newcommand{\fraclpf}[2]{\qp{#1}/{#2}}

\newcommand{\fracl}[2]{{#1}/{#2}}

\newcommand{\old}{\mathrm{old}}
\newcommand{\new}{\mathrm{new}}
\newcommand{\meas}[1]{\operatorname{meas}#1}          

\newcommand{\Oh} {\operatorname{O}}                   
\newcommand{\pd}[2]{\ensuremath{\partial_{#1}{#2}}\xspace} 
\newcommand{\nthpd}[3]{\ensuremath{\partial_{#1}^{#2}{#3}}}
\newcommand{\secpd}[2]{\nthpd{#1}2{#2}} 

\newcommand{\pdt}[1]{\pd t{#1}}                       

\newcommand{\ddt}[1]{\dd t{#1}}    

\newcommand{\transpose}{{\boldsymbol\intercal}}   
\newcommand{\Transpose}[1]{\ensuremath{{#1}^{\transpose}}}
\newcommand{\Transposevec}[1]{\Transpose{\vec{#1}}}
\newcommand{\transposevec}[1]{\Transposevec{#1}}

\newcommand{\colvec}[1]{\ensuremath{\Transpose{\left(#1\right)}}}

\newcommand{\intersected}{\ensuremath{\cap}}
\newcommand{\meet}{\intersected}

\newcommand{\united}{\ensuremath{\cup}}
\newcommand{\join}{\united}
\newcommand{\union}[1]{\ensuremath{\bigcup}_{#1}}

\renewcommand{\vec}[1]{\ensuremath{\boldsymbol{#1}}}
\newcommand{\geovec}[1]{\ensuremath{\overrightarrow{#1}}}

\newcommand{\mat}[1]{\vec{#1}}
\newcommand{\num}[1]{\ensuremath{\mathsf{#1}}}
\newcommand{\numvec}[1]{{\vec{\num{\lowercase{#1}}}}}
\newcommand{\nummat}[1]{\numvec{\uppercase{#1}}}

\newcommand{\transnumvec}[1]{\Transpose{\numvec{#1}}}

\newcommand{\var}{\operatorname{var}}

\newcommand{\boundary}{\partial}

\newcommand{\CC}{\ensuremath{\operatorname C}\xspace}
\newcommand{\HH}{\ensuremath{\operatorname H}\xspace}
\newcommand{\LL}{\ensuremath{\operatorname L}\xspace}

\newcommand{\cont}[1]{\ensuremath{\CC^{#1}}}

\newcommand{\leb}[1]{\ensuremath{\LL_{#1}}}

\newcommand{\sobh}[1]{\ensuremath{\HH^{#1}}}
\newcommand{\sobhz}[1]{\sobh{#1}_0}

\newcommand{\poly}[1]{\ensuremath{\rP}^{#1}}

\newcommand{\fespace}{\rV}

\newcommand{\fes}[1]{\ensuremath{\fespace^{#1}}}

\newcommand{\fesh}{\ensuremath{\fespace_h}}

\newcommand{\EOC}{\ensuremath{\operatorname{EOC}}\xspace}

\newcommand{\Foreach}{\quad\forall\:}

\newcommand{\funk}[3]{\ensuremath{#1:#2\to#3}}

\newcommand{\dfunkmapsto}[6][]{\ensuremath{
    \begin{array}{rccl}
      {#2}: & {#4} &  \to   & {#6}
      \\
      & {#3} &\mapsto & {#5#1}
    \end{array}\quad}}

\newcommand{\downto}{\ensuremath{\searrow}}

\newcommand{\assignvalue}{\ensuremath{:=}}

\renewcommand{\restriction}[2]{\left.#1\right|_{#2}}

\newcommand{\Program}[1]{\textsf{#1}\xspace}

\newcommand{\Algoname}[1]{\ensuremath{\text{\textsf{#1}\xspace}}}

\newcommand{\codename}[1]{\texttt{#1}\xspace}

\newcommand{\alberta}{\Program{ALBERTA}}

\newcommand{\secsymbol}{\S}

\newcommand{\secref}[1]{\secsymbol\ref{#1}}

\newcommand{\eqncomment}[1]{\ensuremath{\qquad\qp{\text{{#1}}}}}

\providecommand{\ListParameters}{}
\renewcommand{\ListParameters}{
  \setlength{\topsep}{0em}
  \setlength{\leftmargin}{0em}
  \setlength{\itemsep}{0ex}
  \setlength{\parsep}{.5ex}
  \setlength{\itemindent}{\labelsep}
  \addtolength{\itemindent}{\labelwidth}
}

\newcounter{LetterListItem}
\renewcommand{\theLetterListItem}{(\alph{LetterListItem})}

\newcounter{NumberListItem}
\renewcommand{\theNumberListItem}{\arabic{NumberListItem}}

\newcounter{QuestionListItem}
\renewcommand{\theQuestionListItem}{\textbf{Question \arabic{QuestionListItem}}}

\newcounter{RomanListItem}
\renewcommand{\theRomanListItem}{(\roman{RomanListItem})}

\newcounter{StepsItem}
\newenvironment{Steps}{%
  \begin{list}%
    {Step \theStepsItem.\ }%
    {\usecounter{StepsItem}%
      \ListParameters
    }
}{%
\end{list}}

\providecommand{\grad}{\nabla}
\renewcommand{\grad}{\nabla}

\usepackage{amsthm} 
\ifthenelse{\boolean{isthesis}}{%
  \setcounter{secnumdepth}{1}%
}{
  \setcounter{secnumdepth}{2} 
}
\providecommand{\ListParameters}{}
\renewcommand{\ListParameters} {
  \setlength{\topsep}{0em}
  \setlength{\leftmargin}{0em}
  \setlength{\itemsep}{0ex}
  \setlength{\parsep}{.5ex}
  \setlength{\itemindent}{\labelsep}
  \addtolength{\itemindent}{\labelwidth}
}
\newtheoremstyle{plain}{}{}{\mdseries\slshape}{\parindent}{\bfseries}{.}{.5em}{}
		
\newtheoremstyle{note}{}{}{}{\parindent}{\bfseries}{.}{.5em}{}
				
\newtheoremstyle{claim}{}{}{\mdseries\slshape}{}{\bfseries}{}{.5em}{}
			
\newtheoremstyle{exercise}
		{}
		{}
		{}
 		{}
		{\bfseries}
		{.}
		{1em}
		{}
		
\newtheoremstyle{break}
		{}
		{}
		{}
		{}
		{\bfseries}
		{.}
		{\newline}
		{}
																								    
\newcommand{\Proofname}{Proof}

\ifthenelse{\boolean{isthesis}}{%
  
}{%
  
}

\newcommand{\pdfformat}[1]{
  \provideboolean{pdfoutput}
  \setboolean{pdfoutput}{#1}
  \ifthenelse{\boolean{pdfoutput}}%
	     {\typeout{using pdf}
	       \usepackage[pdftex]{graphicx,xcolor}
	       \newcommand{\graphext}{pdf}
	       \newcommand{\graphextex}{pdf_t}
	       \usepackage{epsfig}
	       \usepackage{tikz}
	       \usepackage{rotating}
	     }
	     {
	       \typeout{using eps}
	       \usepackage[dvips]{graphicx,xcolor}
	       \newcommand{\graphext}{eps}
	       \newcommand{\graphextex}{eps_t}
	       \usepackage{epsfig}
	       \usepackage{tikz}
	       \usepackage{rotating}
	     }
}
\renewcommand{\aposteriori}{a~posteriori\xspace}
\newcommand{\Apriori}{A~priori\xspace}
\newcommand{\apriori}{a~priori\xspace}
\newcommand{\estimator}{\cE}
\newcommand{\ef}[3]{\ensuremath{\estimator[#1,#2,#3]}}
\newcommand{\honezw}{{\sobhz1(\W)}} 

\newcommand{\tm}{{t_m}}

\newcommand{\tn}{{t_n}}
\newcommand{\tno}{{t_{n-1}}}
\newcommand{\tnp}[1]{\ensuremath{t_{n+#1}}}

\newcommand{\taun}{\ensuremath{\tau_n}}
\newcommand{\taunm}{\ensuremath{\tau_n^{-1}}}

\newcommand{\taunp}[1]{\ensuremath{\tau_{n+{#1}}}}

\newcommand{\hn}{\ensuremath{h_n}}
\newcommand{\hno}{\ensuremath{h_{n-1}}}

\newcommand{\elln}{\ensuremath{l_n}}
\newcommand{\ellno}{\ensuremath{l_{n-1}}}

\newcommand{\fn}{\ensuremath{f^n}}

\newcommand{\U}{\ensuremath{U}}

\newcommand{\un}{\U^n}
\newcommand{\uno}{\U^{n-1}}

\newcommand{\wn}{\ensuremath{\w^n}\xspace}

\newcommand{\an}{A^n}
\newcommand{\rec}{\cR}

\newcommand{\lproj}[1]{P_0^{#1}}

\newcommand{\discopera}{\ensuremath{A}}
\newcommand{\disca}[1]{\discopera^{#1}}
\newcommand{\elop}{\cA}
\renewcommand{\geovec}[1]{\boldsymbol{#1}}
\renewcommand{\vec}[1]{\ensuremath{\boldsymbol
    {\mathsf{#1}}}}
\renewcommand{\mat}[1]{\ensuremath{\boldsymbol
    {\mathsf{#1}}}}

\newenvironment{Proof}[1][{}]%
	       {\par\noindent{\bf \Proofname\ #1}}
	       {{\raggedright{{ }\hfill\qed}}} 
\newtheoremstyle%
    {plain}
    {}
    {}
    {\mdseries\slshape}
    {}
    {\bfseries}
    {.}
    {1.0ex}
    {}
    
\newtheoremstyle
    {note}
    {}
    {}
    {}
    {}
    {\bfseries}
    {.}
    {1.0ex}
    {}
    
\swapnumbers
\theoremstyle{plain}
\newtheorem{The}[subsection]{Theorem}
\newtheorem{Pro}[subsection]{Proposition}
\newtheorem{Pro*}{Proposition}
\newtheorem{Lem}[subsection]{Lemma}

\theoremstyle{note}
\newtheorem{Example}[subsection]{Example}
\newtheorem{Rem}[subsection]{Remark}
\newtheorem{Obs}[subsection]{Remark}
\newtheorem{Defn}[subsection]{Definition}
\newtheorem{Hyp}[subsection]{Assumption}

\renewcommand{\ddt}{\d_t}

\newcommand{\semidiscrete}{semidiscrete\xspace}
\newcommand{\Semidiscrete}{Semidiscrete\xspace}
\newcommand{\fullydiscrete}{fully discrete\xspace}
\newcommand{\Fullydiscrete}{Fully Discrete\xspace}

\newcommand{\Le}[1]{\ensuremath{\leb{#1}}} 
\newcommand{\hoz}{\sobhz1(\W)}

\newcommand{\V}[1]{\fes^{#1}}

\newcommand{\X}[1]{\mathcal{X}^{#1}}
\newcommand{\bi}[2]{\abil{#1}{#2}}

\newcommand{\threshhold}{\xi}
\newcommand{\coarsethresh}{\kappa}

\newcommand{\EI}{\ensuremath{\operatorname{EI}}\xspace}
\newcommand{\tol}{\ensuremath{\operatorname{tol}}\xspace}

\newcommand{\enorm}[1]{\ensuremath{\Norm{#1}}_a}

\newcommand{\dU}[1]{\ensuremath{\partial U^{#1}}}

\newcommand{\T}[1]{\cT^{#1}}
\newcommand{\coarseT}[1]{\cT_0^{#1}}
\newcommand{\patch}[1]{\hat{#1}}

\renewcommand{\ie}{i.e.,\xspace}

\newcommand{\ellop}{\ensuremath{\cA}}
\newcommand{\laginterpol}[1]{\La^{#1}}
\newcommand{\coarseinterpol}[1]{\La_0^{#1}}
\newcommand{\clement}{\Pg}

\newcommand{\hatcleinterpol}[1]{{\hat\clement}^{#1}}

\renewcommand{\meas}[1]{\lvert#1\rvert}
\renewcommand{\fesh}{\ensuremath{\fespace}}
\renewcommand{\integerbetween}[2]{=\ensuremath{#1,\ldots,#2}}
\renewcommand{\rangefromto}[3]{\ensuremath{#1\integerbetween{#2}{#3}}}
\providecommand{\figwidth}{\textwidth}
\renewcommand{\figwidth}{0.9\textwidth}
\newcommand{\figscale}{1.2}
\newcommand{\eint}[3]{\qpbigg{\int_{#1}^{#2}#3^2}^{\kern -.2em\lower .3ex\hbox{\scriptsize $1/2$}}\kern-.3em}
\newcommand{\eenorm}[3]{\eint{#1}{#2}{\enorm{#3}}}
\numberwithin{equation}{section}
\setboolean{showtodo}{false}
\setboolean{shownotes}{false}
\setboolean{showchanges}{false}
\setboolean{usemathrsfs}{false}
\renewcommand{\bi}[2]{\ensuremath{a\qp{#1,#2}}}
\renewcommand{\fes}{\ensuremath{\mathbb V}}
\newcommand{\Estfunk}{\ensuremath{\cE}}
\newcommand{\Est}[1]{\ensuremath{\Estfunk\qb{#1}}}
\newcommand{\projinter}[1]{\ensuremath{\Pg^{#1}}}
\author{Omar Lakkis}
\address{Department of Mathematics, University of Sussex, Falmer near Brighton, England, UK-BN1 9RF}
\email{o.lakkis@sussex.ac.uk}
\urladdr{http://www.maths.sussex.ac.uk/Staff/OL}
\author{Tristan Pryer}
\email{t.m.pryer@sussex.ac.uk}
\urladdr{http://www.sussex.ac.uk/maths/profile131964.html}
\thanks{T.P. was supported by his EPSRC D.Phil. scholarship grant.  
O.L. was partially supported by a Nuffield Young Researcher's Grant.}
\title[Gradient recovery in parabolic problems]{
  Gradient recovery in adaptive finite element methods for parabolic problems} 
\date{\today}
\usepackage{subfigure}
\allowdisplaybreaks
\begin{document}
  \begin{abstract}
    We derive energy-norm \aposteriori error bounds, using gradient
    recovery (ZZ) estimators to control the spatial error, for fully
    discrete schemes for the linear heat equation. This appears to be
    the first completely rigorous derivation of ZZ estimators for
    \emph{fully discrete} schemes for evolution problems, without any
    restrictive assumption on the timestep size.  An essential tool
    for the analysis is the elliptic reconstruction technique.
    
    Our theoretical results are backed with extensive numerical
    experimentation aimed at (a) testing the practical sharpness and
    asymptotic behaviour of the error estimator against the error, and
    (b) deriving an adaptive method based on our estimators.  

    An extra novelty provided is an implementation of a coarsening
    error ``preindicator'', with a complete implementation guide in
    \alberta in the appendix.
  \end{abstract}
  \maketitle
\section{Introduction}
\label{sec:intro}
\emph{Gradient recovery \aposteriori error estimators} have been
widely used since their dissemination in the engineering and
scientific computation community by \cite{Zienkiewicz:1987}, for which
we will often refer to them shortly as \emph{ZZ estimators}.  Since
their introduction they have constituted the most serious rival to
\emph{residual estimators} introduced earlier on in
\cite{Babuska:1978}.  The key to ZZ estimator's success is their
implementation's simplicity, mild dependence upon the problem's data,
and striking superconvergence and asymptotic exactness properties.
\changes{On the other hand, residual estimators, which are the main
  competitor to ZZ estimators, are a bit more involved in
  implementation and cost more to compute, but they are easier to
  handle from the mathematical analysis view-point in deriving
  rigorous upper and lower bounds.}  This situation has led to most of
the theoretical results for evolution equations being obtained in the
last two decades via residual estimators; we refer to
\cite{Lakkis:2006} for a review.  Meanwhile rigorous mathematical work
on recovery estimators has progressed, especially in the last decade,
but mostly for stationary elliptic equations, e.g.,
\cite{Ainsworth:2000,Bank:2003:1,Bank:2003:2,Carstensen:2002:1,Fierro:2006,
  Li:1999,Picasso:2003,Xu:2004}.  \changes{ In contrast, very little
  progress was made on evolution problems, where an exception is
  \cite{Leykekhman:2006}, where the main analytic difficulty comes
  from the singularly perturbed nature of the elliptic problems
  arising from time-stepping procedures.  }

The aim of our work is to bridge the gap between the practical use of
ZZ estimators in adaptivity for evolution equations, studied by
\cite{ziukas-wiberg:98,Picasso:2003}, and the rather mature error
control theory via recovery for stationary equations.  We focus on the
model problem provided by the linear heat equation.
\cite{Leykekhman:2006}, who are to our knowledge the only researchers
to have explored this issue in depth, while obtaining satisfactory
error bounds for spatially discrete schemes, must assume
unrealistically small time-steps for the fully discrete case.  In this
paper we push one step forward by thoroughly analysing the fully
discrete backward Euler schemes.  Namely, we provide reliable error
bounds.  The efficiency and asymptotic exactness of the bounds is
dealt with computationally.

Our main analytical tool to tackle the fully discrete scheme's
difficulties is the \emph{elliptic reconstruction} in the fully
discrete context, studied in \cite{Lakkis:2006}, which provides a way
to take advantage of elliptic \aposteriori error estimates based on
gradient recovery following the exposition of
\cite{Ainsworth:2000}.
  
The elliptic reconstruction technique, introduced under this name by
\cite{Makridakis:2003}, involves the
\emph{error's splitting} into two parts, a \emph{parabolic error}
and an \emph{elliptic error}, through the use of the \emph{elliptic
  reconstruction} of the discrete solution, defined in
(\ref{eqn:def:elliptic-reconstruction}).  This allows to utilise
existing elliptic \aposteriori estimators for the elliptic part and
standard parabolic energy estimates to control the second part.  
\changes{
  Despite this technique being initially introduced to derive sharp
  bounds for lower order spatial error norms, such as
  $\leb2(\W)$ \cite{Makridakis:2003,Lakkis:2006,lakkis-makridakis:07}
  and $\leb\infty(\W)$ \cite{demlow-lakkis-makridakis:09}, we apply
  it here as an analysis tool in an energy-norm framework, where a
  direct approach may lead to a highly complicated analysis for the
  fully discrete scheme.
} In fact, the single most interesting feature of the elliptic
reconstruction, is that \emph{the parabolic error's energy norm term
  is of higher order (with respect to the spatial mesh-size parameter)
  than the elliptic error}, as seen in \cite{Lakkis:2006}.  In this
paper we show, rigorously, that the \emph{full energy error can be
  accounted for only by the elliptic error}, as long as data and
time-step are resolved well enough (cf. Lemma \ref{Lem:SD}).  This
crucial observation is also used to obtain residual \aposteriori
estimates for nonconforming methods in \cite{georgoulis-lakkis:08}.
Note that, it is part of the adaptive methods practitioner's folklore
to employ heuristic versions of this argument.  By way of example, we
quote ``the [full parabolic] discretisation in
energy norm can be bounded by the [elliptic error]
estimator'' from \cite{ziukas-wiberg:98}.
  
Although we treat the case of the Laplace operator, for simplicity, in
this paper, our results can be extended to cover more general elliptic
operators, even time-dependent ones, by using appropriate elliptic
gradient recovery techniques, described by \cite{Fierro:2006}, and a
more careful time-step analysis, as in \cite{georgoulis-lakkis:08}.
  
The paper is organised as follows.  In \secref{sec:setup} we
introduce the model problem, and its discretisations via conforming
finite elements in space and backward Euler in time and we review
the known results, about recovery estimators for elliptic problems,
that will be used in the sequel.  In \secref{sec:semidiscrete} we
describe the elliptic reconstruction technique and illustrate its
use for the spatially \semidiscrete problem.  This paves the way to
tackle the \fullydiscrete problem in \secref{sec:fully-discrete},
where our main results are stated.  In \secref{Numerics}, using
numerical tests, we study the practical behaviour of the estimators
and in \secref{sec:adaptivity} we explore the adaptive schemes based
on our estimators.  

As we have used the finite element toolbox \alberta, written and
documented by \cite{alberta}, for the tests, we have taken the
opportunity to implement a \emph{coarsening preindicator}, previously
unavailable and (for space's sake) fully described in the
Appendix~\ref{sec:building-coarsening-estimator}.  This estimator
predicts the ``information loss'' error that will occur under
coarsening of the mesh at each timestep of the adaptive method and is
crucial in an adaptive code to control information loss during
coarsening.
\section{Set up}
\label{sec:setup}
\subsection{The model problem}
Let $\W\subset\R d$ be a bounded polyhedral domain and consider the
(generalised or weak) Laplace operator denoted by
\begin{equation}
  \dfunkmapsto[.]
	      {\ellop}
	      u
	      {\hoz}
	      {\ellop u:=-\Delta u
		:=-\div\grad u
                =-\sumifromto i1d\secpd i u}
	      {\sobh{-1}(\W)}
\end{equation}
We denote by $\leb2(\W)$ the space of square summable functions on
$\W$, with inner product and norm respectively defined by
\begin{equation}
  \label{eqn:def:ltwo-notation}
  \ltwop fg:=\int_\W f(\vec x)g(\vec x)\d x
  \AND
  \Norm f:=\ltwop ff^{1/2}.
\end{equation}
We will use the standard~\cite[]{Ciarlet:1978,Evans:1998} Sobolev
spaces
\begin{gather}
  \sobh1(\W):=\ensemble{\phi\in\leb2(\W)}{\grad\phi\in\leb2(\W)},
  \\
  \hoz:=\ensemble{\phi\in\sobh{1}(\W)}{\restriction\phi{\boundary\W}=0}
  \\
  \AND 
  \label{eqn:def:space:H-1:dualof:H10}
  \sobh{-1}(\W):=\operatorname{dual}\qp{\hoz}.
\end{gather}
  
Let $T>0$, the model parabolic problem consists in finding a function
$u\in\Le{2}(0,T;\hoz)$ and $\pdt{u}\in\Le{2}(0,T;\sobh{-1}(\W))$ such
that
\begin{equation}
  \begin{split}
    \label{eqn:heat}
    \pdt u(t)+\ellop u(t) 
    &=f(\cdot,t), 
    \text{ for all }t\in\opclinter0T,
    \\
    u(\geovec x,0) 
    &=u_0(\geovec x),\text{ for }\geovec x\in\overline{\W},
    \\
    u(\geovec x,t) 
    &=0,\text{ for }(\geovec x,t)\in\partial\W\times\opclinter0T.
  \end{split}
\end{equation}
We consider the case where $u_0 \in \Le{2}(\W)$ and
$f\in\Le2(0,T;\Le{2}(\W))$ for which the problem \eqref{eqn:heat}
admits a unique solution~\cite[]{Evans:1998}.

Problem \eqref{eqn:heat} is understood in the following weak form
\begin{equation}
  \begin{split}
    \label{eqn:weakheat}
    \ltwop{\pdt u(t)}{\phi}+\bi{u(t)}{\phi}
    &= \ltwop{f(t)}{\phi}
    \Foreach\phi\in\hoz,\,t\in\opclinter0T
    \\
    u(\cdot, 0)
    &=u_0(\cdot),
  \end{split}
\end{equation}
where $\ltwop\cdot\cdot$ is defined in (\ref{eqn:def:ltwo-notation})
and $\bi\phi\psi:=\ltwop{\nabla\phi}{\nabla \psi}$.  The form
$\bi\cdot\cdot$ is clearly bounded and coercive, \ie
\begin{equation}
  \label{coer}
  \bi{\phi}{\phi} \geq \alpha\Norm{\phi}_1^2 
  \Foreach \phi\in\changes{\sobhz1}(\W),
\end{equation}
where $\alpha=(1+\constext P^2)^{-1}$ and $\constext P$ is the 
Poincaré constant. The bilinear form defines an inner product on
$\hoz$ and hence we can denote the energy norm $\enorm{\cdot}^2 :=
\bi{\cdot}{\cdot}$.

\changes{
  \textsl{
    These observations justify our use of $\enorm{\cdot}$ (instead
    of $\Norm{\cdot}_{\sobh1(\W)}$) as the norm of $\honezw$ to be
    with the implied dual norm on $\sobh{-1}(\W)$ in (\ref{eqn:def:space:H-1:dualof:H10}).
  }
}
\subsection{Spatial discretisation}
Let $\T{}$ be a conforming, not necessarily quasiuniform,
triangulation of $\W$, i.e., (1) $K\in\T{}$ means $K$ is an open
simplex (triangle for $d=2$ or tetrahedron for $d=3$), (2) for any
$K,J\in\T{}$ we have that $\closure K\meet\closure J$ is a full
subsimplex (i.e., it is either $\emptyset$, a vertex, an edge, a
face, or the whole of $\closure K$ and $\closure J$) of both $\closure
K$ and $\closure J$ and (3) $\union{K\in\T{}}\closure K=\closure\W$.
The shape regularity of $\T{}$ is defined as the number
\begin{equation}
  \label{eqn:def:shape-regularity}
  \mu(\T{}) := \inf_{K\in\T{}} \frac{\rho_K}{h_K},
\end{equation}
where $\rho_K$ is the radius of the largest ball contained inside
$K$ and $h_K$ is the longest side of $K$. An indexed family of
triangulations $\setof{\T n}_n$ is called \emph{shape regular} if 
\begin{equation}
  \label{eqn:def:family-shape-regularity}
  \mu:=\inf_n\mu(\T n)>0.
\end{equation}
We will use henceforth the usual convention where $\funk h\W\reals$ denotes the
\emph{mesh-size function} of $\T{}$, i.e.,
\begin{equation}
  h(x):=h_{\T{}}(x):=\max_{\closure K\ni x}h_K,
  \AND
  h_n:=h_{\T n}.
\end{equation}

With a triangulation $\T{}$ as described above, and an integer
$p\geq1$ considered fixed in the sequel, we may consider the
\emph{finite element space}
\begin{align}
  \label{eqn:def:finite-element-space}
  \fes
  &:=\ensemble{\Fi \in \hoz}{\Fi\vert_{K} \in \poly p\Foreach K\in\T{}};
\end{align}
and $\poly k$ denotes the linear space of polynomials in $d$
\changes{variables} of degree no higher than a positive integer $k$.  The
\emph{spatially discrete finite element solution} in $\fes$, is the
function $U:[0,T]\rightarrow\fes$ such that
\begin{equation}
  \begin{split}
    \label{eqn:semidiscrete-scheme}
    \ltwop{\pdt U}{\Fi} + \bi{U}{\Fi} 
    &=\ltwop f\Fi 
    \Foreach\Fi \in \fes,
    \\
    U(\geovec x,0) &=U^0:=\projinter\fes u_0(\geovec x)
    \Foreach\geovec x\in\W,
  \end{split}
\end{equation}
where $\funk{\projinter\fes}{\leb2(\W)}\fes$ is a suitable projector
(or an interpolator if the data $u_0$ is in a higher regularity
subspace of $\leb2(\W)$, e.g., $\cT$-wise continuous) and
$\ltwop\cdot\cdot$ is the same as in (\ref{eqn:weakheat}) and
(\ref{eqn:def:ltwo-notation}).

We will often write the scheme \eqref{eqn:semidiscrete-scheme} in its \emph{pointwise form}
\begin{equation}
  \label{eqn:semidiscrete-scheme:pointwise}
  \pdt U+\disca{} U=\lproj{} f\AND U(0)=U^0,
\end{equation}
where the finite dimensional space operator
$\disca{}:\fes\rightarrow\fes$ is the \emph{discrete
  Laplacian} defined, through the Riesz representation in $\fes$, by
\begin{equation}
  \label{defn:discretelap}
  \ltwop{\disca{} V}\Fi = \bi{V}{\Fi}\Foreach\Fi\in\fes,
\end{equation}
and $\lproj{}:\Le2(\W)\rightarrow\fes$ is the
$\leb2(\W)$-projection operator such that, for each $v\in\Le2(\W)$,
\begin{equation}
  \label{defn:l2proj}
  \ltwop{\lproj{} v}{\Fi} = \ltwop{v}{\Fi}\Foreach\Fi\in\fes.
\end{equation}
The pointwise form is convenient as it allows for a more compact
notation.
\subsection{Fully discrete scheme}
\label{sec:fully-discrete-scheme}
Subdivide the time interval $[0,T]$ into a partition of $N$
consecutive adjacent subintervals whose endpoints are denoted
$t_0=0<t_1<\ldots<t_{N}=T$.  The $n$-th timestep is defined as $\tau_n
:= t_n - t_{n-1}$.  We will consistently use the shorthand
$F^n(\cdot):=F(\cdot,t_n)$ for a generic time function $F$. 

The backward Euler method consists in finding a sequence of functions,
$U^n\in\V{n}$, such that for each $\rangefromto n1N$ we have:
\begin{equation}
  \begin{split}
    \label{eqn:fully-discrete}
    \frac1{\tau_n}\ltwop{U^n - \laginterpol n U^{n-1}}{\Fi} + \bi{U^{n}}{\Fi} 
    &=
    \ltwop{f^n}{\Fi}\Foreach \Fi\in \V{n},
    \\ 
    U^0&=\projinter0 u_0,
  \end{split}
\end{equation}
where $\funk{\laginterpol\fes}{\cont0(\W)}{\fes}$ denotes the
Lagrange interpolation operator, $\laginterpol n:=\laginterpol{\V
  n}$, and $\projinter0$ is defined as $\projinter\fes$.

\changes{
  Note our nonrestrictive use of the Lagrange interpolator as a
  ``data-transfer'' operator from a finite element space to the
  next.  We do this to reflect exactly what we do in practical
  computations (where interpolation is faster than averaging).  All
  our analysis applies, however to a different data-transfer
  operator, including the $\leb2(\W)$ projector, if necessary.
}

As with the \semidiscrete scheme the \fullydiscrete scheme can
be written in a pointwise form as follows:
\begin{equation}
  \label{eqn:fully-discretepw}
  \frac{\un-\laginterpol n \uno}{\taun}+\disca{\V n}\un
  =\changes{\lproj n}\fn \AND U^0
  =\changes{\projinter0} u_0,
\end{equation}
\changes{where $\disca n=\disca{\V n}$ and $\lproj n=\lproj{\V n}$ (cf. (\ref{defn:discretelap})).}  
\subsection{Recovery \aposteriori estimators}
\label{Recovery}
The stationary elliptic problem corresponding to a steady state of the
evolution equation \eqref{eqn:heat} is,
\begin{equation}
  \label{eqn:elliptic}
  \text{given $g\in\leb2(\W)$, find $w\in\hoz$ such that } \ellop w = g,
\end{equation}
where the operator is understood in a generalised sense and the
solution is a weak one.  The finite element discretisation of the
elliptic problem \eqref{eqn:elliptic} consists in
\begin{equation}
  \label{eqn:elliptic-discrete}
  \text{finding $W\in\fes$ such that } 
  \bi W\Fi=\ltwop g\Fi\Foreach\Fi\in\fesh.
\end{equation}
We shall henceforth denote by $w$ and $W$ the solutions of
\eqref{eqn:elliptic} and \eqref{eqn:elliptic-discrete}.

From the literature on elliptic \aposteriori estimation
\cite[]{Ainsworth:2000,Babuska:1978,Ciarlet:1978,verfuerth:book,Bank:2003:1,Zienkiewicz:1987}
there is a variety of ways to compute upper and lower bounds 
\changes{for the error in some functional space $\X{}$}
(e.g., $\hoz$, $\leb2(\W)$ and $\leb\infty(\W)$).  For instance, a
generic upper \aposteriori $\X{}$-norm error bound takes the form
\begin{equation}
  \label{elliptic}
  \Norm{w - W}_\X{} \leq \Est{W,g,\X{},\fes},
\end{equation}
where $\Estfunk$ is an appropriate \emph{estimator functional}.

One way of providing an estimator functional consists, for example, in
starting by applying a \emph{gradient postprocessing operator
  (postprocessor)}, say $G$, to the approximate solution $W$. And then
proving that $\Norm{G W-\grad W}$ is equivalent to the error
$\Norm{\grad w-\grad W}$.  \emph{Gradient recovery operators} form a
subclass of gradient postprocessors.

Recovery operators can be built in a variety of ways such as local
weighted averaging (where the gradient is sampled from neighbouring
elements)~\cite[]{Picasso:2003}, discrete $\Le2(\W)$-projection (using
least squares fitting)~\cite[]{Zienkiewicz:1987} or global
$\Le2(\W)$-projection (where a full discrete problem is
solved)~\cite[]{Bank:2003:1}. In our numerical results we use local
weighted averaging, defined explicitly in \eqref{eqn:recovery-average}, but our
theoretical results can be applied with any choice of recovery
operator that provide upper and lower bounds for the elliptic problem.
The fundamental idea behind these approaches is to build an
approximation $G{W}$ of $\grad w$ which is more regular than the
piecewise discontinuous gradient $\grad W$; the extra regularity is
aimed at obtaining a higher approximation order.
\begin{Defn}[gradient recovery operator, from \cite{Ainsworth:2000}]
  \label{def:gradient-recovery}
  A \emph{gradient recovery (ZZ) operator} on $\fes$ is a linear
  operator $\funk G\fes{\fes^d}$ which enjoys the following properties:
  \begin{description}
  \item[Consistency] we have, with
    $\funk{\laginterpol\fes}{\cont 0(\W)}{\fes}$ denoting the Lagrange
    interpolator,
    \begin{equation}
      \restriction{G(\laginterpol{\fes} v)}{K} = \restriction{\nabla v}{K}
      \Foreach v\in\poly{p+1},\,K\in\T{}.
    \end{equation}
  \item[Local bound] there exists a $\constext{ZZ}>0$ such that
    \begin{equation}
      \Norm{G{V}}_{\Le\infty(K)} 
      \leq\constext{ZZ}\Norm{\grad V}_{\Le{\infty}(\patch{K})}
      \Foreach V\in\fes,\,K\in\T{},
    \end{equation}
    where $\patch K$ is the \emph{patch} generated by $K$ (the union
    of all $L\in\T{}$ such that $\closure L\meet\closure
    K\neq\emptyset$).
  \end{description}
\end{Defn}
For simplicity, we assume that the operator is in a mesh-local
relation with $\grad V$ noting, nonetheless, that global methods such
as the global $\Le2(\W)$-projection proposed by
\cite{Bank:2003:1,Bank:2003:2} exist and can be included in our
discussion.

Under certain regularity assumptions recovery estimators are shown to be
asymptotically exact.  For instance, \cite{Zlamal:1977} shows
that if $w\in\sobh{s+1}(\W)$, with reference to (\ref{eqn:elliptic}) and
(\ref{eqn:elliptic-discrete}), its approximation $W\in\fes$ satisfies the
following \emph{superconvergence property}:
\begin{equation}
  \label{super}
  \Norm{\nabla (W-\laginterpol{\fes}w)} 
  =\Oh(h^{\changes p+\zeta})\text{ for some }\zeta\in\opclinter01.
\end{equation}
A review of superconvergence results is given by
\cite{Krizek:1987}. If (\ref{super}) is satisfied then the recovered
gradient also satisfies the following superconvergence
property~\cite[]{Ainsworth:2000}:
\begin{equation}
  \label{recovery:super}
  \Norm{\nabla w - G W} 
  =\Oh(h^{\changes p+\zeta})\text{ for some }\zeta\in\opclinter01.
\end{equation}
The reach of Zlámal's result is appreciated by stating the following
consequence.
\begin{Lem}[gradient recovery \aposteriori estimate from \cite{Ainsworth:2000}]
  \label{lem:gradient-recovery-estimate}
  Let $\fes$ be the finite element space defined in
  (\ref{eqn:def:finite-element-space}) and $\funk G\fes{\fes^d}$ a
  gradient recovery operator according to
  \secref{def:gradient-recovery}. If $w,W$ are the solutions of
  (\ref{eqn:elliptic}) and (\ref{eqn:elliptic-discrete}), respectively, and
  \eqref{recovery:super} holds then the recovery operator is
  asymptotically exact, in the sense that
  \begin{equation}
    \label{eqn:gradient-recovery-asymptotically-exact}
    \lim_{h_{\T{}}\to0}
    \frac{\Norm{\nabla W - G W}}{\Norm{\nabla(W - w)}}=1.
  \end{equation}
  Thus, there exist $\delta_0\geq0$, such that $\delta_0(h)\to0$ as $h\to0$ and
  \begin{equation}
    (1-\delta_0)\Norm{\grad W-G{W}}\leq\Norm{\grad[W-w]}\leq (1+\delta_0)\Norm{\grad W-G W}
  \end{equation}
  for all partitions $\T{}$ of $\W$ satisfying $h_{\T{}}<h_0$.
\end{Lem}
\begin{Obs}[recovery in absence of regularity]
  Lacking Zlámal's regularity assumption, recovery-based estimators
  are empirically observed to be efficient, reliable estimators, even
  on meshes with low
  shape-regularity~\cite[]{Carstensen:2004:1}.
  
  For more details about recovery-based estimators we refer to the
  available literature~\cite[]{Bank:2003:1,Bank:2003:2,Xu:2004,Li:1999,Ainsworth:2000,Fierro:2006}.
\end{Obs}
\begin{Defn}[gradient recovery \aposteriori estimator functional]
  \label{def:estimator-functional}
  Lemma \ref{lem:gradient-recovery-estimate} then justifies the use of
  the \emph{recovery estimator} in the $\hoz$-norm (and by equivalence
  the energy norm) by defining, for the rest of this paper, the
  \emph{gradient recovery \aposteriori estimator functional}
  \begin{equation}
    \label{eqn:def:estimator-functional}
    \Est V:=\Est{V,\hoz,\fes}:=\Norm{G{V}-\grad V},\text{ for }V\in\fes,
  \end{equation}
  where $G$ is a gradient recovery operator as defined
  in~\S\ref{def:gradient-recovery}.
\end{Defn}
\begin{Hyp}[elliptic \aposteriori error estimates]
  \label{hyp:aposteriori-estimates}
  We will consider henceforth the blanket assumption that \changes{
    for a fixed $h_0$, there are some $c_0<C_0$, such that for any $\fes$
    with mesh-size $h<h_0$, for $w$ and $W$ solutions of
    (\ref{eqn:elliptic}) and (\ref{eqn:elliptic-discrete}),
    respectively and $\Estfunk$ defined in
    \ref{def:estimator-functional} the following bounds are true
    \begin{equation}
      \label{eqn:hyp:elliptic-aposteriori-estimates}
      c_0\Est W\leq\Norm{\grad\qb{W-w}}\leq C_0\Est W.
    \end{equation}
  }
  
  \changes{
    Optionally, we will assume \emph{asymptotic exactness}, in which
    case
    \begin{equation}
      \label{eqn:hyp:asymptotic exactness}
      C_0\leq1+B(h_0)\AND c_0\geq1+\beta(h_0),
    \end{equation}
    for some continuous functions $B$ and $\beta$ that vanish at
    $0$.
  }
  
  \changes{ Assumptions~(\ref{eqn:hyp:elliptic-aposteriori-estimates})
    and~(\ref{eqn:hyp:asymptotic exactness}) are true, modulo
    hierarchic \emph{oscillation terms} of the data function $g$
    in~(\ref{eqn:elliptic-discrete}). These assumptions are thus
    justified, for example, when $g$ is in a finite dimensional space,
    for example $g\in\fespace$ as we shall assume in the sequel, by
    isolating the bulk of the oscillations in data-approximation
    terms.  For a thorough discussion of the oscillation in the
    context of recovery, we refer to \cite{Fierro:2006}.}
  
  The lower bound is not needed for the theory to be developed
  herein, as we will prove only upper bounds.  Nonetheless, this
  property is required for the efficiency of the parabolic
  estimators in practical situations.
\end{Hyp}
\section{\Semidiscrete scheme}
\label{sec:semidiscrete}
To make the link between the parabolic problem and the elliptic
recovered gradient estimates we utilise the elliptic reconstruction
technique \cite[]{Makridakis:2003,Lakkis:2006}.  To make the
discussion more accessible, we first do this for the spatially
(semi)discrete scheme.  We divide the error into two parts---one
called elliptic error the other parabolic error---via the
\emph{elliptic reconstruction of the discrete solution}.  Because the
elliptic error can be directly bounded under the blanket Assumption
\ref{hyp:aposteriori-estimates}, it is enough to show that the full
error can be bounded in terms of the elliptic error only.  This result
is in accordance with the fact that the parabolic error on uniform
meshes is of higher $h$-order in the energy norm with respect to the
elliptic (and thus the full) error, as observed by \cite{Lakkis:2006}.
The main result of this section is summarised in
Theorem~\ref{the:aposteriori-semidiscrete}.
\begin{Defn}[elliptic reconstruction]
  \label{def:ert}
  The \emph{elliptic reconstruction operator} is defined as
  $\funk\rec\fes\hoz$ such that
  \begin{equation}
    \label{eqn:def:elliptic-reconstruction}
    \ellop[\rec V]=\disca{} V,
  \end{equation}
  where $\disca{}$ is the discrete elliptic operator defined in
  (\ref{defn:discretelap}).  In weak form, equation
  (\ref{eqn:def:elliptic-reconstruction}) reads
  \begin{equation}
    \bi{\rec  V}\Fi = \ltwop{\disca{} V}\Fi
    \Foreach\Fi\in\hoz,
  \end{equation}
  and it is well defined in virtue of the elliptic problem's well
  Sydney's.  We will refer to the function $\rec V$ as the
  \emph{elliptic reconstruction} of $V$, while the elliptic
  reconstruction operator $\rec$ will be called the
  \emph{reconstruction operator} (or just the \emph{reconstructor})
  from $\fes$.
  
  If $U(t)$ denotes the solution of
  (\ref{eqn:semidiscrete-scheme:pointwise}) at time $t$, we shall
  indicate by $\w(t)$ its reconstruction $\rec U(t)$.
  
  \label{def:ert:remark}
  Thus, posing $g(t):=\disca{} U(t)$, we then see $U(t)$ is the finite element
  solution corresponding to the elliptic problem of finding
  $\w(t)\in\hoz$ such that $\ellop\w(t)=g$.
\end{Defn}
\subsection{The error and its splitting}
\label{sec:error-its-splitting}
For the whole of this section we shall consider $u$ to be the solution
of (\ref{eqn:heat}), understood in the weak sense, and $U$ its
semidiscrete approximation given by (\ref{eqn:semidiscrete-scheme:pointwise}).  The
corresponding \emph{semidiscrete error} is defined by
\begin{equation}
  e(t):= U(t)-u(t),
\end{equation}
and can be split, using the elliptic reconstruction $\w=\rec U$
defined in \S\ref{def:ert}, as follows:
\begin{equation}
  \label{eqn:error-splitting}
  \begin{split}
    e(t)=
    \left(\w(t) - u(t)\right) - \left(\w(t) - U(t)\right)
    =:\rho(t) - \epsilon(t).
  \end{split}
\end{equation}
We shall refer to $\epsilon$ and $\rho$ here defined as the
\emph{elliptic (reconstruction) error} and the \emph{parabolic error}
respectively.
Using this notation we have the estimate
\begin{equation}
  \begin{split}
    \Norm{\nabla\qb{U-u}(t)} 
    &\leq \Norm{\nabla\rho(t)} + \Norm{\nabla\epsilon(t)},
  \end{split}
\end{equation}
where, following the remarks made in Definition \ref{def:ert:remark} and
Assumption~\ref{hyp:aposteriori-estimates}, the elliptic error can be
bounded by the computable elliptic \aposteriori estimator functional
$\Estfunk$:
\begin{equation}
  \enorm{\epsilon(t)}=\Norm{\nabla\epsilon(t)} \leq C_0\Est{U(t)}.
\end{equation}
It is therefore sufficient to bound the error's energy norm using the
elliptic error's energy norm.
\begin{Lem}[elliptic energy bound for parabolic \semidiscrete error]
  \label{Lem:SD}
  If $e,\,\epsilon$ are defined as in
  \S\ref{sec:error-its-splitting} then, for each $t\in\clinter0T$,
  we have
  \begin{equation}
    \Norm{e(t)}^2+
    \int_0^t \enorm{e(s)}^2 \d s \leq \Norm{e(0)}^2 
    +\int_0^t \enorm{\epsilon(s)}^2
    +2\ltwop{\lproj{} f(s)-f(s)}{e(s)}\d s.
  \end{equation}
\end{Lem}
\begin{Proof}
  From the the exact problem (\ref{eqn:heat}), the \semidiscrete
  scheme \eqref{eqn:semidiscrete-scheme:pointwise}, and the splitting \eqref{eqn:error-splitting}
  \begin{equation}
    \begin{split}
      \label{eqn:proof:lem:semidiscrete:identity}
      \pdt e+\ellop\rho
      =
      \pdt{\qb{U-u}}+\ellop\qb{\w-u}
      =
      \pdt U+\disca{} U
      -
      \pdt u-\ellop u
      =
      \lproj{} f-f.
    \end{split}
  \end{equation}
  Testing with $e$ we obtain
  \begin{equation}
    \ltwop{\pdt e}{e}+\bi\rho{e}
    =\ltwop{\lproj{} f-f}{e}
  \end{equation}
  and thus
  \begin{equation}
    \frac{1}{2}\ddt{\Norm{e}^2}
    +\enorm{e}^2
    =
    \ltwop{\lproj{} f-f}{e}-\bi{\epsi}e.
  \end{equation}
  Integration from $0$ to $t$ yields 
  \begin{equation}
    \Norm{e(t)}^2+2\int_0^t \enorm{e}^2
    =
    \Norm{e(0)}^2
    +2\int_0^t\ltwop{\lproj{} f-f}e
    -2\int_0^t\bi{\epsilon}e
    \Foreach t\in\clinter0T.
  \end{equation}
  Hence, by Young's inequality on $\bi\epsi{e}$, we have
  \begin{equation}
    \Norm{e(t)}^2+2\int_0^t\enorm{e}^2
    \leq
    \Norm{e(0)}^2
    +2\int_0^t\ltwop{\lproj{} f-f}{e}
    +\int_0^t\enorm{e}^2
    +\int_0^t\enorm\epsilon^2,
  \end{equation}
  whereby the claim is verified.
\end{Proof}
\begin{Obs}[proliferation of $\sqrt2$ syndrome]
  Let $a,b,c\geq0$ such that $a^2\leq c^2+ab$, then, by Young's
  inequality, it follows that $a^2\leq 2c^2+b^2$.  Note however that the
  factor ``$2$'' in $2c^2$ is not needed, in that we also have that
  $a\leq c+b$.  If $1>a\sim c\gg b>0$, then the first bound provides
  $a/c\approx\sqrt2$ whereas the second bound gives $a/c\approx1$, which
  is tighter.
  
  The following result, which generalises $a\leq c+b$, is extremely
  simple yet useful in avoiding this ``proliferation of $\sqrt2$
  syndrome'' from repeated usage of Young's inequality.
\end{Obs}
\begin{Pro}[$\leb2$ simplification rule]
  \label{pro:L2-simplify}
  If $\numvec a,\numvec b\in\R N$, $N\in\naturals$, $c\in\reals$ and
  $f,g\in\leb2(D)$, for some measurable domain $D$, are such that
  \begin{equation}
    \label{eqn:L2-simplify:hypothesis}
    \norm{\numvec a}^2+\Norm{f}^2\leq c^2
    +\Transpose{\numvec a}\numvec b
    +\int_D f g,
  \end{equation}
  then
  \begin{equation}
    \label{eqn:L2-simplify:conclusion}
    \qp{\norm{\numvec a}^2+\Norm f^2}^{1/2}
    \leq
    \abs c+\qp{\norm{\numvec b}^2+\Norm g^2}^{1/2},
  \end{equation}
  where all the vector norms are Euclidean, and the function norms
  $\leb2(D)$.
\end{Pro}
\begin{Proof}
  \newcommand{\va}{\vec\a}
  \newcommand{\vb}{\vec\b}
  \newcommand{\nva}{\norm{\va}}
  \newcommand{\nvb}{\norm{\vb}}
  
  Denote by $\va:=(\norm{\numvec a},\Norm{f})$
  and $\vb:=(\norm{\numvec b},\Norm{g})$.
  
  If $\nva\leq\nvb$ then (\ref{eqn:L2-simplify:conclusion}) is
  trivially satisfied.  Otherwise we have $\nva>\nvb$ whereby
  (\ref{eqn:L2-simplify:hypothesis}) and the
  Cauchy--Bunyakovskiy--Schwarz inequality imply that
  \begin{equation}
    \begin{split}
      \nva^2
      &\leq c^2
      +\norm{\numvec a}\norm{\numvec b}
      +\Norm f\Norm g
      +\nvb\qp{\nva-\nvb}
      \\
      &\leq
      c^2+2\nva\nvb-\nvb^2.
    \end{split}
  \end{equation}
  Hence $\qp{\nva-\nvb}^2\leq c^2$, and thereby
  \begin{equation}
    \nva\leq \abs c+\nvb,
  \end{equation}
  as claimed.
\end{Proof}
\begin{The}[\aposteriori \semidiscrete error estimate]
  \label{the:aposteriori-semidiscrete}
  With $u$ and $U$ as defined by (\ref{eqn:heat}) and (\ref{eqn:semidiscrete-scheme}),
  respectively, and an estimator functional $\Estfunk$ as
  defined in (\ref{eqn:def:estimator-functional}), we have
  \begin{multline}
    \label{eqn:semidiscrete-aposteriori-estimate}
    \psqrt{\Norm{U(t)-u(t)}^2+\int_0^t\enorm{U-u}^2}
    \\
    \leq
    \Norm{U(0)-u(0)}
    +
    C_0\Norm{\Est{U}}_{\leb2\clinter0T}
    +
    2\Norm{\lproj{} f-f}_{\leb2(0,T;\sobh{-1}(\W))}
    .
  \end{multline}
\end{The}
\begin{Proof} Using Lemma \ref{Lem:SD} we have
  \begin{equation}
    \Norm{e(t)}^2+
    \int_0^t \enorm{e}^2
    \leq 
    \Norm{e(0)}^2 
    +\int_0^t \enorm{\epsilon}^2
    +2\int_0^t\ltwop{\lproj\fespace f-f}{e}.
  \end{equation}
  Using Proposition \ref{pro:L2-simplify}, we
  obtain
  \begin{equation}
    \psqrt{\Norm{e(t)}^2+\int_0^t\enorm{e}^2}
    \leq
    \psqrt{\Norm{e(0)}^2+\int_0^t \enorm{\epsilon}^2}
    +2\psqrt{\int_0^t\Norm{\lproj{} f-f}_{\sobh{-1}(\W)}^2}.
  \end{equation}
  Assumption (\ref{eqn:hyp:elliptic-aposteriori-estimates}) and the
  discussion in \S\ref{sec:error-its-splitting} ensure then that
  \begin{multline}
    \psqrt{\Norm{e(t)}^2+\int_0^t\enorm{e}^2}
    \\
    \leq
    \psqrt{\Norm{e(0)}^2+C_0^2\int_0^t \Est U^2}
    +2\psqrt{\int_0^t\Norm{\lproj{} f-f}_{\sobh{-1}(\W)}^2},
  \end{multline}
  which implies the claim.
\end{Proof}
\begin{Obs}[short versus long integration times]
  The bound for the pointwise in time $\leb2(\W)$ error,
  $\Norm{e(t)}$, appearing on the \changes{left-hand side}
  of~(\ref{eqn:semidiscrete-aposteriori-estimate}), is tight only for
  very short times.  For example, it is well-known that on a uniform
  mesh of size $h\to0$ on a convex domain $\W$  the energy term
  $\psqrt{\int_0^t\enorm e^2}$ is $\Oh(h^p)$, while $\Norm{e(t)}$ is
  $\Oh(h^{p+1})$.
\end{Obs}
\begin{Obs}[dealing with the $\sobh{-1}(\W)$ norm]
  \label{obs:negative-sobolev-norms}
  \changes{
    In practise the $\sobh{-1}(\W)$ norm can be well approximated as
    shown by Lemma \ref{lem:computing-negative-sobolev-norms}, so, in
    the lack of a priori information, the last term in
    (\ref{eqn:semidiscrete-aposteriori-estimate}) may be replaced
    using the Poincaré inequality
    \begin{equation}
      2\Norm{\lproj{} f-f}_{\leb2(0,T;\sobh{-1}(\W))}
      \leq
      2\constext{P}(\W)\Norm{\lproj{} f-f}_{\leb2(\W\times\opinter0T)}.
    \end{equation}
  }
  
  It is also possible to obtain bounds by using the
  Cauchy--Bunyakovskiy--Schwarz inequality for $\leb2(\W)$ on the
  term $\ltwop{\lproj{} f-f}{e}$---rather than the $(\sobh{-1},\sobhz1)$
  duality---and ``absorb'' the resulting $\sup_{\clinter0t}\Norm e$
  into the first term on the right hand side of
  (\ref{eqn:semidiscrete-aposteriori-estimate}).  However, whenever
  possible, we shy away from this procedure as it incurs in
  artificially higher constants and a $\leb1\clinter0T$ accumulation
  on the right-hand side while the energy term on the left-hand side
  accumulates like $\leb2\clinter0T$.  This time-accumulation
  disparity between the error and the estimator is likely to result
  in an error--estimator ratio bound that has the order of $\sqrt
  T$, that is, although having the right order of convergence, the
  estimator will overestimate the error over long integration times.
\end{Obs}
\changes{
  We show now how to practically approximate the $\sobh{-1}(\W)$ norm of
  an arbitrary given function $v\in\leb2(\W)$.
}
\changes{
  \begin{Lem}[computing the $\sobh{-1}(\W)$ norm]
    \label{lem:computing-negative-sobolev-norms}
    Let $v\in\leb2(\W)$, consider the functions $\psi\in\sobhz1(\W)$ and  $\Ps\in\fespace$ such
    that
    \begin{equation}
      \label{eqn:duality-continuous-discrete}
      \elop\psi=v
      \AND
      \disca{}\Ps=\lproj{} v,
    \end{equation}
    where $\disca{}$ and $\lproj{}$ are the discrete Laplacian and the
    $\leb2(\W)$ projection on $\fespace$, respectively.  Then,
    recalling our convention whereby $\Norm v_{\honezw}=\Norm{\grad
      v}$ we have
    \begin{equation}
      \begin{split}
        \Norm v_{\sobh{-1}(\W)}^2
        =
        \Norm\psi_{\honezw}^2
        =
        \Norm{\psi-\Psi}_{\honezw}^2
        +\Norm \Psi_{\honezw}^2.
      \end{split}
    \end{equation}
  \end{Lem}
}
\changes{
  \begin{Proof}
    With $\psi$ and $\Ps$ as given in
    \eqref{eqn:duality-continuous-discrete} we have $\Fi\in\fespace$
    \begin{equation}
      \duality{\elop\psi-\disca{}\Ps}\Fi=
      \ltwop{v-\lproj{} v}\Fi=0,
    \end{equation}
    i.e., that $\psi-\Psi$ is Galerkin-orthogonal to $\fespace$.      
    Also, we have
    \begin{equation}
      \Norm v_{\sobh{-1}(\W)}=\Norm{\psi}_{\sobhz1(\W)}.
    \end{equation}
    Indeed, on the one hand
    \begin{equation}
      \begin{split}
        \Norm v_{\sobh{-1}(\W)}
        :=
        \sup_{\phi\in\honezw}
        \frac{
          \ltwop v\phi}{
          \Norm\phi_\honezw}
        =
        \sup_{\phi\in\honezw}
        \frac{
          \ltwop{\grad v}{\grad \phi}}{
          \Norm\phi_\honezw}
        \\
        \leq
        \sup_{\phi\in\honezw}
        \frac{
          \Norm\psi_\honezw\Norm\phi_\honezw}{
          \Norm\phi_\honezw}
        =
        \Norm\psi_\honezw,
      \end{split}
    \end{equation}
    and, on the other hand
    \begin{equation}
      \Norm v_{\sobh{-1}(\W)}
      :=
      \sup_{\phi\in\honezw}
      \frac{
        \ltwop{\grad\psi}{\grad\phi}}{
        \Norm\phi_\honezw}
      \geq
      \frac{
        \ltwop{\grad\psi}{\grad\psi}}{
        \Norm\psi_\honezw}
      =\Norm\psi_\honezw.
    \end{equation}
    By the above, Galerkin-orthogonality and Pythagoras's Theorem, we have
    \begin{equation}
      \begin{split}
        \Norm v_{\sobh{-1}(\W)}^2
        =
        \Norm\psi_{\honezw}^2
        =
        \Norm{\psi-\Psi}_{\honezw}^2
        +\Norm \Psi_{\honezw}^2.
      \end{split}
    \end{equation}
  \end{Proof}
}
\begin{Obs}[the $\sobh{-1}(\W)$ norm approximation]     
\changes{
    The next-to-last term $\Norm{\psi-\Psi}_{\honezw}$ is the error of a
    function and its Ritz projection.  This can be easily estimated with a
    fully computable \aposteriori estimator functional $\estimator$ such
    that
    \begin{equation}
      \Norm{\psi-\Psi}_{\honezw}
      \leq
      \ef \Psi v\fespace
      =\Oh(h_\fespace^r),
    \end{equation}
    where $h_\fespace$ is the ``mesh-size'' of the space $\fespace$.
}

\changes{
  Hence the $\sobh{-1}(\W)$ can be computed using the relation:
  \begin{equation}
    \begin{split}
      \Norm v_{\sobh{-1}(\W)}^2
      =
      \Norm \Ps_{\honezw}^2
      +
      \zeta[\Ps,v]^2,
    \end{split}
    \text{ where }
    \zeta[\Ps,v]\leq\Est{\Ps}.
  \end{equation}
}

\changes{    
  The term $\Norm{\Ps}_{\honezw}$ is clearly computable, by
  computing $\Ps$, which involves one $\leb2(\W)$-projection, one
  stiffness matrix inversion and one (discrete) energy norm
  computation.  Furthermore
  \begin{equation}
    \Norm v_{\sobh{-1}(\W)}^2
    =\Norm \Ps_{\honezw}^2+\Oh(h_\fespace^{2r}).
  \end{equation}
  Hence, if $\Ps$ is finite with respect to the mesh-size
  $h_\fespace$, i.e., $\Norm{\Ps}=\Oh(h_\fespace^0)$, then we can
  approximate the $\sobh{-1}(\W)$ of a function with as much precision
  as the finite element method allows it for the energy norm.  On the
  other hand if $\Ps$ is small, precisely,
  $\Ps=\Oh(h_\fespace^s)$ with $s>0$ (implying that $v$ is
  small as well), then this result has to be handled with more care
  for the error to be of some order of $h_\fespace$ higher than the
  computed quantity.  }
\end{Obs}
\begin{Obs}[sharper versions of Theorem~\ref{the:aposteriori-semidiscrete}]
  The error estimate (\ref{eqn:semidiscrete-aposteriori-estimate}) can
  be tightened further to
  \begin{multline}
    \psqrt{\frac12\Norm{e(t)}^2+\int_0^t\enorm e^2}
    \\
    \leq
    \frac1{\surd2}\Norm{e(0)}
    +
    \psqrt{\int_0^t\Norm{\lproj{} f-f}_{\sobh{-1}(\W)}+C_0^2\Est U^2}.
  \end{multline}
  But this estimate becomes noticeably better only when one of the
  terms $\Norm{e(0)}$ or $\Norm{\lproj{} f-f}_{\sobh{-1}(\W)}$ dominates
  the $\Est U$ term, which should not be allowed to happen.  So
  there is no need to lengthen the discussion by insisting on such
  tight bounds, as long as it is possible to obtain the elliptic
  \aposteriori estimate constant $C_0$ in the leading term on the
  right-hand side.
\end{Obs}
\section{\Fullydiscrete scheme}
\label{sec:fully-discrete}
The main result of this section---and the
paper---is the \aposteriori error bound, stated in Theorem
\ref{the:fully-discrete-aposteriori-bound}, on the error
between the approximate solution $U$ of the fully discrete problem
(\ref{eqn:fully-discretepw}) and that of the exact problem (\ref{eqn:heat}).

The analysis in this section follows narrowly the one we performed in
\secref{sec:semidiscrete}, albeit with the complications that the fully
discrete scheme imports. We will first extend the discrete solution
sequence to a continuous-time function.  Then we derive an error
identity on which we mimic the energy techniques of \secref{sec:semidiscrete}
to bound the error's energy norm in terms of some residual terms and
the elliptic error's energy norm, which is finally controlled via gradient
recovery estimators.
\subsection{Time extension of the discrete solution}
Recalling the fully discrete scheme (\ref{eqn:fully-discretepw}),
the fully discrete solution is the sequence of finite element
functions $\un\in\V n$ defined at each discrete time $\tn$,
$\rangefromto n0N$.  Define the piecewise linear (affine) extension
\begin{equation}
  \label{eqn:time-extension}
  U(t):=\sumifromto n0N \un l_n(t),
\end{equation}
where we use the one-dimensional piecewise linear continuous Lagrange
basis functions\changes{, defined for $t\geq0$}, as
\begin{equation}
  l_n(t):=
  \begin{cases}
    \fraclpf{t-\tno}\taun,
    &\text{ for }t\in\opclinter\tno\tn\text{ \changes{(and $n>0$)}},
    \\
    \fraclpf{\tnp1-t}{\taunp1},
    &\text{ for }t\in\opclinter\tn{\tnp1}
    \\
    0,
    &\text{ otherwise.}
  \end{cases}
\end{equation}
\textsl{We warn the reader that we use the same symbol, $U$, to indicate the
  fully discrete solution's extension to $\clinter0T$, as the one we
  used for its semidiscrete counterpart in \S~\ref{sec:semidiscrete}.}
\subsection{Elliptic reconstruction and error splitting}
Next we define the elliptic reconstruction, needed for the following analysis,
similarly to that of the semidiscrete scheme
(cf. \eqref{eqn:def:elliptic-reconstruction}).  For each
$n\in\fromto0N$, with the discrete elliptic operator $\an$ as in
\ref{sec:fully-discrete-scheme}, we define the corresponding
elliptic reconstruction operator $\rec^n:\V{n}\rightarrow\hoz$, for
each $V\in\V n$, by solving for $\rec^n V$ the elliptic problem
\begin{equation}
  \label{eqn:def:fully-discrete:elliptic-reconstruction-n}
  \ellop\rec^n V=\an V,
\end{equation}
which can be read in the weak form as
\begin{equation}
  \bi{\rec^nV}\Fi=\ltwop{\an V}\Fi\Foreach\Fi\in\hoz.
\end{equation}
We denote 
\begin{equation}
  \wn:=\rec^n\un,\text{ for each }\rangefromto n0N, 
\end{equation}
and this sequence's piecewise linear extension by
$\funk\w{\clinter0T}\hoz$, i.e.,
\begin{equation}
  \label{eqn:def:fully-discrete:elliptic-reconstrution-ext}
  \w(t):=\sum_{n=0}^{N}\wn l_n(t).
\end{equation}

As in the semidiscrete analysis we introduce symbols for the
\emph{full error} $e:=U-u$, the \emph{elliptic error} $\epsi:=\w-U$ and the
\emph{parabolic error} $\rho:=\w-u$, whereby
\begin{equation}
  e=\rho-\epsi,
\end{equation}
and, based on the Assumption \ref{hyp:aposteriori-estimates},
\begin{equation}
  \label{eqn:basic-elliptic-error-bound-by-estimator-functional}
  \begin{split}
    \enorm{\epsi(t)}
    &\leq 
    C_0\Est{\un\elln(t)+\uno\ellno(t)}
    \\
    &\changes{\leq 
      C_0\qp{\Est\un\elln(t)+\Est\uno\ellno(t)}}
    \text{ for }t\in\clinter\tno\tn.
  \end{split}
\end{equation}
The last step is guaranteed by the linearity of the operators $G$ and
$\grad$, hence the homogeneity $\Est{\la V}=\abs\la\Norm{G{V}-\grad V}$,
\changes{and by the triangle inequality}.
\begin{Lem}[parabolic error identity]
  \label{lem:paraberrid}
  For each $\rangefromto n1N$ and each $t\in\opinter\tno\tn$ we have
  \begin{equation}
    \label{eqn:parabolic-error-identity}
    \pdt e(t)+\ellop\rho(t)= 
    \fraclpf{ 
      \laginterpol n U^{n-1}
      -
      U^{n-1}
    }\taun
    +\ellop\qb{\w(t) - \w^n}+\lproj n f^n-f(t).
  \end{equation}
\end{Lem}
\begin{Proof}
  By the definition of $U$, (\ref{eqn:time-extension}),
  for each $\rangefromto n1N$ and $t\in\opinter\tno\tn$
  we have
  \begin{equation}
    \pdt U(t)=\un\elln'(t)+\uno\ellno'(t)=\fraclpf{\un-\uno}\taun
  \end{equation}
  and using the fully discrete scheme~(\ref{eqn:fully-discretepw}), we have
  \begin{equation}
    \begin{split}
      \pdt U(t)+\ellop\wn
      &=
      \changes{
	\fraclpf{
          \laginterpol n\uno
          -
          \uno
        }\taun
	+\fraclpf{
          \un
          -
          \laginterpol n\uno
        }\taun
      }
      +{\an\un}
      \\
      &=
      \fraclpf{
        \laginterpol n\uno
        -
        \uno
      }\taun
      +\lproj n f^n.
    \end{split}
  \end{equation}
  Hence
  \begin{equation}
    \pdt U(t)+\ellop\w(t)
    =
    \changes{\fraclpf{
        \laginterpol n\uno
        -
        \uno
      }\taun}
    +\ellop\qb{\w(t)-\wn}
    +\lproj n f^n
  \end{equation}
  and, using the exact PDE (\ref{eqn:heat}), we get
  \begin{equation}
    \begin{split}
      \pdt e(t)+\ellop\rho(t)
      &=
      \pdt U(t)+\ellop\w(t)-\pdt u(t)-\ellop u(t)
      \\
      &=
      \changes{
	\fraclpf{
          \laginterpol n\uno
          -
          \uno}
        \taun}
      +\ellop\qb{\w(t)-\wn}
      +\lproj n f^n-f(t),
    \end{split}
  \end{equation}
  as stated.
\end{Proof}
\begin{Defn}[\aposteriori error indicators]
  \label{def:fully-discrete:estimators}
  The notation we introduce here will be valid for the rest of the
  article.
  \begin{description}
  \item
    [{elliptic error indicator via recovery}]
    \begin{equation}
      \label{recoveryestimator}
      \varepsilon_n
      :=\Est{U^n,\hoz,\V n}
      =\changes{C_0}
      \Norm{\grad U^n - G^n[U^n]},
    \end{equation}
    with the functional $\Estfunk$ as defined in
    \S\ref{def:estimator-functional}, and\footnote{In the numerical
      experiments we use $(\ep_n^2+\ep_{n-1}^2)/2$ instead of
      $\tilde\ep_n$.}
    \begin{equation}
      \label{eqn:def:recovery-estimator:mixed}
      \tilde\ep_n^2
      :=\frac13\qp{\ep_n^2+\ep_{n-1}^2+\ep_n\ep_{n-1}}
      \leq\frac12\qp{\ep_n^2+\ep_{n-1}^2}.
    \end{equation}
  \item
    [{time-discretisation error indicators}]
    \begin{equation}
      \label{eqn:def:time-indicator}
      \theta_n : =
      \frac1{\sqrt3}
      \begin{cases}
	\Norm{
	  {\lproj n f^n-\laginterpol n\dU n} 
	  -\qp{\lproj{n-1}f^{n-1}-\laginterpol{n-1}\dU{n-1}}
	}_{\sobh{-1}(\W)}
	&\text{ for }n \geq 2,\\
	\Norm{{\lproj 1f^1-\laginterpol 1\dU1}-A^0U^0}_{\sobh{-1}(\W)}
	&\text{ for }n = 1,
      \end{cases}
    \end{equation}
    \changes{where $\dU n:=\fraclpf{\un-\uno}\taun$,
      (cf. Lemma~\ref{lem:computing-negative-sobolev-norms}), also
      possible to use in its alternative (faster to compute but not as
      sharp)} version
    \begin{equation}
      \label{eqn:def:time-indicator:alternative}
      \tilde{\theta}_n := \constant\mu\enorm{U^{n-1} - U^n},
    \end{equation}
    where $\constant\mu$ is dependent on the shape regularity $\mu$ of the
    family of triangulations defined in (\ref{eqn:def:family-shape-regularity}).
  \item
    [{mesh-change (coarsening) error indicators}]
    a main mesh-change indicator
    \begin{equation}
      \label{eqn:def:coarsening-indicator}
      \gamma_n:=  
      \tau_n^{-1}
      \Norm{\laginterpol n U^{n-1}-U^{n-1}}_{\sobh{-1}(\W)},
    \end{equation}
    and a \emph{higher order} mesh-change indicator
    \begin{equation}
      \label{eqn:def:coarsening-indicator:higher-order}
      \tilde{\gamma}_n
      :=
	    {\constant\mu}'
	    \begin{cases}
	      \Norm{\smash{\hat h_n}\qp{{\lproj n f^n-\laginterpol n\dU n} 
	          -\lproj{n-1}f^{n-1}+\laginterpol{n-1}\dU{n-1}}}
	      ,
	      &\: n\geq2,
	      \\
	      \Norm{\smash{\hat h_1}\qp{{\lproj 1f^1-\laginterpol 1\dU1}-A^0U^0}}
	      ,
	      &\: n=1,
	    \end{cases}
    \end{equation}
    where
    $\hat{h}_n(\geovec{x})
    =
    \max\setof{h_{n-1}(\geovec{x}),h_n(\geovec{x})}$
    for $\geovec x\in\W$ and a constant $\constant\mu'$.
  \item
    [{data approximation error indicator}]
    \begin{equation}
      \label{eqn:def:data-approx-indicator}
      \beta_n := 
      \tau_n^{-1} 
      \int^{t_n}_{t_{n-1}}\Norm{\lproj n f^n - f(t)}_{\sobh{-1}(\W)} \d t.
    \end{equation}
  \end{description}
\end{Defn}
\changes{
  \begin{Obs}[computing $\sobh{-1}(\W)$ norms]
    Clearly the $\sobh{-1}(\W)$ norms appearing in Definition
    \ref{def:fully-discrete:estimators} cannot be computed in practise.
    The corresponding indicators can be replaced by upper bounds using
    the (dual) Poincaré inequality
    \begin{equation}
      \label{eqn:dual-Poincare}
      \Norm\phi_{\sobh{-1}(\W)}\leq\constext P\Norm\phi.
    \end{equation}
    Other alternatives, as described in Remark
    \ref{obs:negative-sobolev-norms} are possible but will not be
    described here.
  \end{Obs}
}
\begin{The}[\aposteriori estimate for fully discrete scheme]
  \label{the:fully-discrete-aposteriori-bound}
  Let the sequence $\seqsu U n{\fromto0N}$, $U^n\in\V{n}$, be the solution of the
  \fullydiscrete problem \eqref{eqn:fully-discrete} and $U$ its piecewise
  linear time-extension as in \eqref{eqn:time-extension}. Let $u$ be the
  exact solution of the exact problem \eqref{eqn:heat} then
  \begin{equation}
    \begin{split}
      \psqrt{\frac{\Norm{U^N-u(T)}^2}2
	+\int_0^T\enorm{U(t)-u(t)}^2 \d t}
      \leq&
      \frac{\Norm{U(0)-u(0)}}{\surd2}
      +\eta_N
    \end{split}
  \end{equation}
  where the \emph{(global) error estimator} is given by the following
  discrete $\leb2(0,T)$ summation of the error indicators defined in
  \S\ref{def:fully-discrete:estimators}:
  \begin{gather}
    \label{estimatorstheta}
    \eta_N^2=\sumifromto n1N\qp{\tilde\ep_n+\gamma_n+\beta_n+\theta_n}^2\taun.
  \end{gather}
\end{The}
\begin{Proof}
  The proof shadows that of Lemma \ref{Lem:SD} and Theorem
  \ref{the:aposteriori-semidiscrete}, but we must take into account
  the complications arising from the time discretisation.  For the
  reader's convenience we divide it into steps.
  \begin{Steps}
  \item
    Using the notation from Lemma \ref{lem:paraberrid} and identity
    (\ref{eqn:parabolic-error-identity}) therein we have that
    \begin{equation}
      \begin{split}
	\pdt e(t)+\ellop e(t)=
	\ellop\epsi(t)
	+\fraclpf{\laginterpol n U^{n-1} - U^{n-1}}\taun
	\\
	+\ellop\qb{\w(t) - \w^n}+\lproj n f^n-f(t).
      \end{split}
    \end{equation}
    Testing this with $e$ we obtain
    \begin{equation}
      \begin{split}
	\frac12\ddt\Norm{e(t)}^2+\enorm{e(t)}^2
	&=
	\bi{\epsi(t)}{e(t)}
	+\ltwop{\fraclpf{\laginterpol n U^{n-1} - U^{n-1}}\taun}{e(t)}
	\\
	&\phantom=
	+\ltwop{\ellop\qb{\w(t) - \w^n}}{e(t)}
	+\ltwop{\lproj n f^n-f(t)}{e(t)},
      \end{split}
    \end{equation}
    for all $t\in\opinter\tno\tn$ and each $n\integerbetween 1N$.
    Integrating over $[0,T]$ gives
    \begin{equation}
      \label{eqn:def:fully-discrete:lemma:proof:basic}
      \begin{split}
	\Norm{e^N}^2/2
	&
	+\int_0^T \enorm{e(t)}^2 \d t 
	=
	\Norm{e^0}^2/2+\int_0^T \bi{\epsilon(t)}{e(t)}\d t
	\\&
	+
	\sumifromto n1N
	\int_\tno^\tn
	\ltwop{\fraclpf{\laginterpol n U^{n-1}-U^{n-1}}\taun}{e(t)} 
	\\&
	\phantom{\sumifromto n1N \int_\tno^\tn}
	+\bi{\w(t)-\w^n}{e(t)} 
	+\ltwop{\lproj n f^n - f(t)}{e(t)} 
	\d t
	\\
	&=:
	  {\mathcal{B}_1+\mathcal{B}_2
	    +\mathcal{B}_3+\mathcal{B}_4}+\Norm{e^0}^2/2.
      \end{split}
    \end{equation}
    We proceed by bounding each of the terms $\mathcal{B}_j$,
    $\rangefromto j14$, appearing in the right-hand side of
    (\ref{eqn:def:fully-discrete:lemma:proof:basic}).
  \item
    The first term to be bounded in
    (\ref{eqn:def:fully-discrete:lemma:proof:basic}) yields the spatial
    discretisation error indicator as follows:
    \begin{equation}
      \label{eqn:fully-discrete:lemma:proof:space}
      \begin{split}
	\mathcal{B}_1
	&=
	\int_0^T \bi{\epsilon(t)}{e(t)} \d t 
	= \sum_{n=1}^N \int_{t_{n-1}}^{t_n} \bi{\epsilon(t)}{e(t)} \d t
	\\&
	\leq
	\sumifromto n1N
	\eenorm\tno\tn\epsi\eenorm\tno\tn e
	\leq
	\sumifromto n1N 
	\tilde\ep_n\taun^{1/2}
	\eenorm\tno\tn e
      \end{split}
    \end{equation}
    where we have used~(\ref{eqn:def:recovery-estimator:mixed}) and
    in view
    of~(\ref{eqn:basic-elliptic-error-bound-by-estimator-functional})
    and~(\ref{eqn:def:recovery-estimator:mixed}), we may write
    \begin{equation}
      \int_\tno^\tn\enorm{\epsi}^2
      \leq
      \ep_{n-1}^2\int_\tno^\tn l_{n-1}^2+2\ep_{n-1}\ep_n\int_\tno^\tn l_{n-1}l_n+\ep_n^2\int_\tno^\tn l_n^2
      =
      \tilde\ep_n^2\taun
      .
    \end{equation}
    \par
    The second term in
    (\ref{eqn:def:fully-discrete:lemma:proof:basic}) contains mesh-change
    term which we bound as follows:
    \begin{equation}
      \label{eqn:fully-discrete:lemma:proof:coarsening}
      \begin{split}
	\mathcal{B}_2
	&=
	\sum_{n=1}^{N}\int_{t_{n-1}}^{t_n} 
	\ltwop{\fraclpf{\laginterpol n U^{n-1}-U^{n-1}}\taun}{e(t)}\d t
	\\
	&\leq 
	\sum_{n=1}^N 
	    {\Norm{\laginterpol n U^{n-1}-U^{n-1}}_{\sobh{-1}(\W)}}\taunm
	    \int_{t_{n-1}}^{t_n} \enorm{e(t)} \d t
	    \\
	    &\leq 
	    \sum_{n=1}^N\gamma_n
	    \taun^{1/2}
	    \eenorm\tno\tn e
      \end{split}
    \end{equation}
    where $\gamma_n$ is defined by (\ref{eqn:def:coarsening-indicator}).
    \par
    Similarly the data error term is bounded as follows
    \begin{equation}
      \label{eqn:fully-discrete:lemma:proof:data}
      \begin{split}
	\mathcal{B}_4 
	&= \int_0^T \ltwop{\lproj n f^n - f(t)}{e(t)} \d t
	\leq
	\sum_{n=1}^N \beta_n\taun^{1/2}\eenorm\tno\tn e,
      \end{split}
    \end{equation}
    where $\beta_n$ is defined in (\ref{eqn:def:data-approx-indicator}).
  \item 
    The third term in (\ref{eqn:def:fully-discrete:lemma:proof:basic})
    yields a time discretisation term and is a bit more involved to estimate.
    Using the definition of $\w^n$, $\w$ and $\rec^n$, given
    in~(\ref{eqn:def:fully-discrete:elliptic-reconstruction-n}) and
    (\ref{eqn:def:fully-discrete:elliptic-reconstrution-ext}),
    we observe that
    \begin{equation}
      \label{eqn:fully-discrete:lemma:proof:time}
      \begin{split}
	\mathcal{B}_3
	=& \sum_{n=1}^N\int_{t_{n-1}}^{t_n} \bi{\w - \w^n}{e(t)}
	\d t
	\\
	=& \sum_{n=1}^N \int_{t_{n-1}}^{t_n} \bi{l_{n-1}(t)
	  \rec^{n-1}U^{n-1} + l_n(t) \rec^n U^n -
	  \rec^n U^n}{e(t)} \d t\\
	=& \sum_{n=1}^N \int_{t_{n-1}}^{t_n} l_{n-1}(t)
	\bi{\rec^{n-1}U^{n-1} - \rec^n U^n}{e(t)} \d t
	\\ 
	=& \sum_{n=1}^N 
	\int_{t_{n-1}}^{t_n} l_{n-1}(t)
	\ltwop{A^{n-1}U^{n-1}-A^n U^n}{e(t)}
	\d t
	\\ 
	\leq&
	\sum_{n=1}^N
	\Norm{A^{n-1}U^{n-1} - A^n U^n}_{\sobh{-1}(\W)}
	     {\eint\tno\tn{l_{n-1}}}
	     \eenorm\tno\tn e
	     \\
	     \leq&
	     \sum_{n=1}^N 
	     \theta_n{\taun^{1/2}}\eenorm\tno\tn e
	     ,
      \end{split}
    \end{equation}
    where in the last passage we use the discrete scheme (\ref{eqn:fully-discretepw}) for the substitution
    \begin{equation}
      \label{eqn:fully-discrete:lemma:proof:time:substitution}
      A^n\un=\fraclpf{\laginterpol n \uno-\un}{\taun}+\lproj n\fn
      \text{ for }n\geq1.
    \end{equation}
  \item Grouping together
    (\ref{eqn:def:fully-discrete:lemma:proof:basic}),
    (\ref{eqn:fully-discrete:lemma:proof:space}),
    (\ref{eqn:fully-discrete:lemma:proof:coarsening}),
    (\ref{eqn:fully-discrete:lemma:proof:data}) and
    (\ref{eqn:fully-discrete:lemma:proof:time}), we obtain
    \begin{multline}
      \Norm{e^N}^2/2
      +\int_0^T \enorm{e(t)}^2 \d t
      \\
      \leq
      \Norm{e^0}^2/2
      +\sumifromto n1N
      \qp{\tilde\ep_n+\gamma_n+\beta_n+\theta_n}\taun^{1/2}
      \eenorm\tno\tn e.
    \end{multline}
    Using an $\leb2$ simplification (cf. \secref{pro:L2-simplify}), we
    conclude that
    \begin{equation}
      \psqrt{
	\frac{\Norm{e^N}^2}2
	+\int_0^T \enorm{e(t)}^2 \d t
      }
      \leq
      \frac{\Norm{e^0}}{\surd2}
      +
      \psqrt{
	\sumifromto n1N
	\qp{\tilde\ep_n+\gamma_n+\beta_n+\theta_n}^2\taun
      }.
    \end{equation}
  \end{Steps}
  Referring to the notation in (\ref{eqn:time-extension}) and Definition
  \ref{def:fully-discrete:estimators}, we obtain the result.
\end{Proof}
\begin{Rem}[the alternative time indicator]
  \label{alternativetime}
  \changes{ 
    Assuming there is no mesh change from time $\tno$ to
    time $\tn$, then the discrete Laplacians $\disca{n-1}$ and
    $\disca n$, defined in~(\ref{defn:discretelap}), coincide.  Thus
    the time discretisation error indicator $\theta_n$, which is
    part of the estimator $\eta_N$ in
    Theorem~\ref{the:fully-discrete-aposteriori-bound}, can be
    written as
    \begin{equation}
      \theta_n
      =\frac1{\surd2}\Norm{A^n\un-A^{n-1}\uno}_{\sobh{-1}(\W)}
      =
      \frac{\taun}{\surd2}\Norm{A^n\pdt U}_{\sobh{-1}(\W)}.
    \end{equation}
    In the form given in~(\ref{eqn:def:time-indicator}) and using the
    dual Poincaré inequality~(\ref{eqn:dual-Poincare}), this indicator
    is easily bounded.
  }
  
  \changes{
    A more precise, but slightly more expensive,
    calculation can be done using
    Lemma~\ref{lem:computing-negative-sobolev-norms}.  The same
    idea, will be used in the next result where we show that the
    indicator $\theta_n$ is equivalent, up to higher order
    terms, to the alternative time indicator $\tilde\theta_n$,
    defined in~(\ref{eqn:def:time-indicator:alternative}), which
    requires only an energy norm computation. 
  } 
  This alternative time indicator, which is more common in energy
  estimates \cite[e.g.]{picasso:98}, $\tilde\theta_n$ is also more
  ``natural'', as it measures the time derivative in the energy norm as
  opposed to the $\sobh{-1}$ norm of the time derivative of $AU$.
  Due to mesh-change effects, this simpler indicator comes at the
  (affordable) price of having to add the higher order mesh change term
  $\tilde\gamma_n$ to the otherwise simpler $\gamma_n$.
\end{Rem}
\begin{The}[alternative time estimator]
  \label{the:fully-discrete-aposteriori-bound:alternative}
  With the same assumptions and notation of
  Theorem~\ref{the:fully-discrete-aposteriori-bound}
  we have
  \begin{equation}
    \begin{split}
      \psqrt{\frac{\Norm{U^N-u(T)}^2}2
	+\int_0^T\enorm{U(t)-u(t)}^2 \d t}
      \leq&
      \frac{\Norm{U(0)-u(0)}}{\surd2}
      +\tilde\eta_N
    \end{split}
  \end{equation}
  where the \emph{(alternative global) error estimator} is given by the following
  discrete $\leb2(0,T)$ summation of the error indicators defined in
  \S\ref{def:fully-discrete:estimators}:
  \begin{gather}
    \tilde\eta_N^2:=\sumifromto n1N\qp{\tilde\ep_n+\gamma_n+\tilde\gamma_n+\beta_n+\tilde\theta_n}^2\taun.
  \end{gather}
\end{The}
\begin{Proof}
  We proceed similarly to the proof of
  Theorem~\ref{the:fully-discrete-aposteriori-bound}, in steps. The
  notation is the same and steps 1 and 2 are identical.
  \begin{Steps}
    \setcounter{StepsItem}2
  \item
    This step starts similarly to its homologue in the proof of
    Theorem~\ref{the:fully-discrete-aposteriori-bound} by observing
    that
    \begin{equation}
      \mathcal{B}_3 
      =
      \sum_{n=1}^N \int_\tno^\tn l_{n-1}(t) 
      \ltwop{A^{n-1}U^{n-1}-A^n U^n}{e(t)}\d t.
    \end{equation}
    The function $A^{n-1}U^{n-1}-A^n U^n$ belongs to $\V{n}+\V{n-1}$,
    but in general it is in neither of $\V n$ nor $\V{n-1}$). Thus, to
    proceed, we use the $\leb2(\W)$-projection and the
    Clément--Scott--Zhang interpolator denoted respectively by
    \begin{equation}
      \funk{\check{P}^n}{\Le2(\W)}{\V{n}+\V{n-1}}
      \AND
      \funk{\hatcleinterpol n}{\Le2(\W)}{\V{n}\meet\V{n-1}}.
    \end{equation}
    We recall that the operators $\hatcleinterpol n$ and $\check P^n$ are both
    known~\cite[resp.]{Scott:1990,Carstensen:2002:L2stab} to enjoy the
    following stability properties for all $n\integerbetween0N$:
    \begin{gather}
      \label{eqn:clement-stable}
      \enorm{\hatcleinterpol n\phi} 
      \leq 
      \constant{1,\mu}\enorm{\phi}
      \Foreach\phi\in \sobh1(\W),
      \\
      \label{eqn:l2-projection-stable}
      \enorm{ \check{P}^n \phi} \leq 
      \constant{2,\mu} \enorm{\phi} 
      \Foreach\phi\in \sobh1(\W),
    \end{gather}
    where $\mu$ is the shape-regularity of the triangulation family
    $\setof{\T n}_{n\integerbetween0N}$ defined in
    (\ref{eqn:def:family-shape-regularity}).  Furthermore, the following
    interpolation inequality is valid~\cite[\S{B.3}]{Lakkis:2006}
    \begin{equation}
      \label{eqn:clement:interpolation-inequality}
      \Norm{
        \fraclpf{
          \psi-\hatcleinterpol{n}\psi
        }{
          \hat h_n
        }
      }
      \leq
      \constant{3,\mu}\enorm{\psi}
      \Foreach\psi\in\hoz,\,n\integerbetween1N,
    \end{equation}
    where $\hat h_n:=\max\setof{\hn,\hno}$.
  \item
    Using these operators, we derive that
    \begin{equation}
      \begin{split}
	\mathcal{B}_3 
	=& 
	\sum_{n=1}^N \int_{t_{n-1}}^{t_n}
	\ltwopbig{A^{n-1}U^{n-1} - A^n U^n}{\check{P}^n e(t)} 
	l_{n-1}(t)\d t
	\\
	=&
	\sum_{n=1}^N \int_{t_{n-1}}^{t_n} 
	\qpBig{
	  \ltwopbig{A^{n-1}U^{n-1}-A^n U^n}{
	    \check{P}^n e(t)-\hatcleinterpol n\check{P}^n e(t)} 
	  \\
	  &
	  \phantom{\smash{\sum_{n=1}^N \int_{t_{n-1}}^{t_n}}}
	  + \ltwopbig{A^{n-1}U^{n-1}-A^n U^n}{
	    \hatcleinterpol n\check{P}^n e(t)} 
	}
	l_{n-1}(t)\d t
	\\
	\leq&
	\sum_{n=1}^N \int_{t_{n-1}}^{t_n}
	\qpBig{
	  \Norm{\smash{\hat{h}_n}\qpbig{A^{n-1}U^{n-1} - A^n U^n}}\Normbig{\smash{\hat{h}_n^{-1}}
	    \qpbig{\check{P}^n e(t) - \hatcleinterpol n \check{P}^n e(t)}} 
	  \\
	  & \qquad
	  + 
	  \bi{U^{n-1} - U^n}
	     {\smash{\hatcleinterpol n \check{P}^n e(t)}}
	}
	l_{n-1}(t)\d t.
      \end{split}
    \end{equation}
    
    Using inequalities (\ref{eqn:clement-stable}),
    (\ref{eqn:l2-projection-stable}) and
    (\ref{eqn:clement:interpolation-inequality}), we get the bound
    \begin{equation}
      \begin{split}
	\mathcal{B}_3 
	&\leq
	\sumifromto n1N\int_\tno^\tn 
	\qpbig{
	  \constant{3,\mu}\Norm{\smash{\hat h_n}(A^{n-1}U^{n-1}-A^n{U}^n)}
	  \enorm{\check{P}^n e(t)}
	  \\
	  &\phantom{\smash{\leq\sumifromto n1N\int_\tno}}
	  +\constant{1,\mu}
	  \enorm{U^{n-1}-U^n}
	  \enorm{\check{P}^n{e}(t)}
	}
	l_{n-1}(t)\d t
	\\
	&\leq
	\sum_{n=1}^N
	\Big(
	\constant{3,\mu}
	\Norm{\smash{\hat{h}_n}(A^{n-1}U^{n-1}-A^n{U}^n)}
	+
	\constant{1,\mu}
	\enorm{U^{n-1} - U^n}
	\Big)
	\\
	&\phantom{\smash{\leq\sumifromto n1N}}
	\times
	\constant{2,\mu}\int_{t_{n-1}}^{t_n}
	\enorm{e(t)}l_{n-1}(t)\d t
	\\
	&\leq
	\sumifromto n1N
	\qp{\tilde\gamma_n+\tilde\theta_n}
	\eenorm\tno\tn e
      \end{split}
    \end{equation}
    by taking $\constant\mu:=\constant{1,\mu}\constant{2,\mu}/3$,
    ${\constant\mu}':=\constant{3,\mu}\constant{2,\mu}/3$ in
    (\ref{eqn:def:time-indicator:alternative}) and
    (\ref{eqn:def:coarsening-indicator:higher-order}) for the last step.
  \end{Steps}
  We may now conclude exactly like the last step in the proof of
  Theorem~\ref{the:fully-discrete-aposteriori-bound}, albeit with
  $\theta_n$ replaced by $\tilde\gamma_n+\tilde\theta_n$.
\end{Proof}
\section{Computer experiments: convergence rates}
\label{Numerics}
In this section and in \S\ref{sec:adaptivity} we study the numerical
behaviour of the error indicators and estimators and compare this
behaviour with the true error on three model problems.  The
\Program{C} code that we used includes the
adaptive FEM library \alberta~\cite[]{alberta}.  The quadrature formal error is
made negligible with respect to other error by using overkill
quadrature formulas (exact on polynomials of degree $17$ and less).
\subsection{Benchmark problems}
Consider three benchmark problems, the solution of which is
known. Namely, take $d=2$, each problem's data $f,u_0$ is then chosen
such that the exact solution to \ref{eqn:heat} is given by:
\begin{align}
  \label{eqn:Problem:1}
  u(\geovec{x},t)&= \sin\qp{\pi t}\exp\qp{-10\norm{\geovec{x}}^2},
  \\
  \label{eqn:Problem:3}
  u(\geovec{x},t)&=  \sin\qp{20\pi t}\exp\qp{-10\norm{\geovec{x}}^2},
  \\
  \label{eqn:Problem:2}
  u(\geovec x,t)&= t\sin\frac{2\arctan(x_2/x_1)}3
  {\geovec{x}}^{\fracl23}\exp\qp{\frac{-1}{1-\norm{\geovec x}^2}},
\end{align}
The domain $\W$ for Problems \eqref{eqn:Problem:1} and
\eqref{eqn:Problem:3} is the square
$S:=\opinter{-1}1\times\opinter{-1}1$.  Problem
\eqref{eqn:Problem:2}, whose solution's gradient is singular at the
origin, is considered on the L shaped domain
$\W=S\take[0,1]\times[-1,0]$.
The benchmark problems \eqref{eqn:Problem:1} and
\eqref{eqn:Problem:3} have been chosen such that they can be
compared with previous numerical studies~\cite[]{Lakkis:2006}.

For all Problems \eqref{eqn:Problem:1}--\eqref{eqn:Problem:2}, we
take zero initial condition, $u_0 = 0$ to avoid dealing with the
initial adaptivity which is a side issue here.

The solution \eqref{eqn:Problem:3} has a time dominant
discretisation error, while \eqref{eqn:Problem:2} was constructed
to have a dominant spatial error. It is the product of a linear
function in time, a well known solution to Laplace's equation
producing the spatial singularity and a mollifier.

Problem \eqref{eqn:Problem:1} is used to test asymptotic behaviour
of the indicators under uniform space-time refinements further in
\S\ref{sec:experimental-order-of-convergence}.  Problems
\eqref{eqn:Problem:2} and \eqref{eqn:Problem:3} will be used to test
the adaptive strategies in \S\ref{sec:adaptivity}.
\subsection{Gradient recovery implementation}  
The recovery operator, $G^n$, is obtained by taking the
discontinuous gradients of the numerical solution at the
super convergent sampling points \cite[]{Ainsworth:2000} (and references
therein). The recovery operator used here is built in the following
way: fixing $V\in\V n$, for each degree of freedom $\geovec x$, we define
\begin{equation}
  \label{eqn:recovery-average}
  G^n[V](\geovec x):=
  \frac
      {\sum_{K\in\T n:\geovec x\in K}\meas K\restriction{\grad V}K(\geovec x)}
      {\sum_{K\in\T n:\geovec x\in K}\meas K},
\end{equation}
This defines a unique piecewise polynomial field $G^n[V]\in\fesh^d$.
(Note that formula~(\ref{eqn:recovery-average}) is non trivial for
only for those DOF that are are on the boundary of an element; for
the internal DOF, that arise in using $\poly p$ elements for
$p\geq3$, it is not necessary to calculate anything.)
\begin{Defn}[experimental order of convergence]
  \label{def:EOC}
  Given two sequences $a(i)$ and $h(i)\downto0$,
  $i\integerbetween l{}$, we define experimental order of convergence
  (\EOC) to be the local slope of the $\log a(i)$ vs. $\log h(i)$
  curve, i.e.,
  \begin{equation}
    \EOC(a,h;i):=\frac{ \log(a(i+1)/a(i)) }{ \log(h(i+1)/h(i)) }.
  \end{equation}
\end{Defn}
\begin{Defn}[effectivity index]
  \label{def:EI}
  The main tool deciding the quality of an estimator is the
  effectivity index (\EI) which is the ratio of the error and the
  estimator, \ie
  \begin{equation}
    \EI(\tn):= 
    \fracl{\eta_n}{\Norm{U-u}_{\Le{2}(0,t_n;\hoz)}}.
  \end{equation}
  If $\EI(\tn)\to1$ as $\sup_{x,n} h_n(x)\to0$ then we say the estimator is
  \emph{asymptotically exact}.
\end{Defn}
\subsection{Indicator's numerical asymptotic behaviour}
\label{sec:experimental-order-of-convergence}
In the following convergence rate tests we discuss the practical
realisation of Theorems~\ref{the:fully-discrete-aposteriori-bound}
and~\ref{the:fully-discrete-aposteriori-bound:alternative}, to which
we refer for notation.

We use a uniform timestep and uniform meshes that are fixed with
respect to time.  Hence for each test we have $\V n=\V 0=\fes$ and
$\taun=\tau(h)$ for all $n\integerbetween 1N$.  For each test we fix
the polynomial degree $p$ and two parameters $k,c$ and then compute
a sequence of solutions with $h=h(i)=2^{-i/2}$, and $\tau=ch^k$ for a
sequence of refinement levels $i=l,\dotsc,L$.

Due to the finite element space invariance in time, the coarsening
indicator $\gamma_n$ vanishes and is thus not computed (this indicator
will be discussed in~\secref{sec:adaptivity}).

The initial value being zero makes the initial error $U(0)-u(0)$
zero.  Thus we do not need to calculate this term in the estimator.

For all solutions the boundary values are not exactly zero, but of
a negligible value, hence little interpolation error is committed
here (nonetheless some care is taken when dealing with very small
errors).  Finally, the data approximation error term, $\beta_n$, though
important for highly oscillatory data, will not be studied here
given the regularity of our data.

Therefore, what we compute on a space-time uniform mesh are the
indicators $\ep_n$ and $\theta_n$ (or $\tilde\theta_n$,$\gamma_n$),
defined in~\S\ref{def:fully-discrete:estimators}, and the
corresponding \emph{cumulative indicators} $\seqsifromto En1N$ and
$\seqsifromto\Th n1N$ defined by:
\begin{equation}
  \label{eqn:def:cumulative-indicators}
  \begin{split}
    &E_m^2:=\sumifromto n1m\qp{\ep_n^2+\ep_{n-1}^2}\taun/2\eqncomment{for space},
    \\
    \AND
    &\Th_m^2:=\sumifromto n1m\theta_n^2\taun
    \OR
    \sumifromto n1m\qp{\tilde\theta_n^2+\tilde\gamma_n^2}\taun
    \eqncomment{for time}.
  \end{split}
\end{equation}
From the Theorems~\ref{the:fully-discrete-aposteriori-bound}
and~\ref{the:fully-discrete-aposteriori-bound:alternative}, we know
that
\begin{equation}
  \Norm{U^n-u(\tn)}^2\leq E_m^2+\Th_m^2+{\sumifromto n1m\beta_n^2\taun}.
\end{equation}

Our results and the comments are reported in the captions of
figures.

In Figures \ref{fig:P1}--\ref{fig:P4} we visualise the results and
comment them, for Problem~(\ref{eqn:Problem:1}) for conforming
finite elements of polynomial degree $p=1,\dotsc,4$, respectively.
Having fixed $p,k,c$ such that $\tau=ch^k$, for each level $i$, we
plot $\Th_m$ and $E_m$, $\Norm{U-u}_{\leb2(0,\tm;\hoz)}$, their
experimental order of convergence and the effectivity index
$\EI(\tm)$ versus (discrete) time $\tm=0,\dotsc,T$.  The conclusion
is that the estimator is sharp and reliable, but to achieve
asymptotic exactness (or close) the time indicator must be made
smaller than the space indicator by taking $\tau\ll h^p$.  In all
these tests we used the first form for $\Th_m$ appearing in
(\ref{eqn:def:cumulative-indicators}).

In Figure \ref{Fig:comparetheta} we summarise a comparison between
the two time indicators $\theta_n$ and $\tilde\theta_n$, showing
that the latter yields a much sharper bound, but with the added cost
of having to compute the higher order term $\tilde\gamma_n$.
\renewcommand{\figscale}{0.75}
\begin{figure}[ht]
  \caption{\label{fig:P1}
    Numerical Results for Problem \eqref{eqn:Problem:1} with
    $\poly1$ and $h=h(i)=2^{-i/2}$, $i=4,\dotsc,9$ (details
    in~\S\ref{sec:experimental-order-of-convergence}).
  }
  \begin{center}
    \subfigure[{\label{fig:P1h:tau}
	Mesh-size is $h$ and timestep $\tau=0.1\,h$.  On top we plot
	the \EOC's of the single cumulative indicators $E$ and
	$\Th$.  Below we plot their logs.  Both indicators have
	$\EOC\to1$, but the cumulative time error indicator $\Th_m$
	is dominant.  The estimator is reliable and sharp, but not
	asymptotically exact and results in $\EI\gg1$.
    }]{
      \includegraphics[scale=\figscale,width=\figwidth]{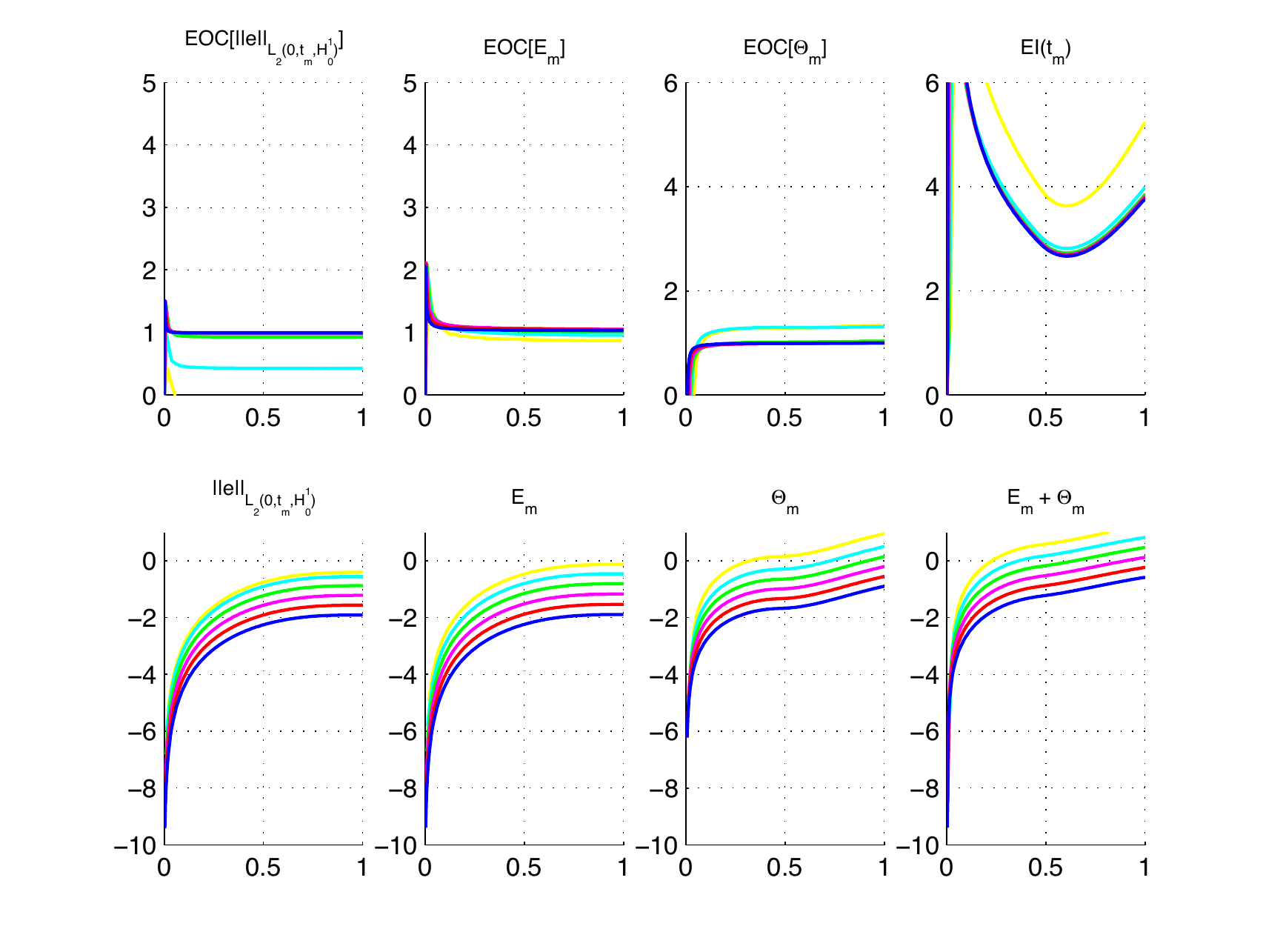}
    }
    \\
    \subfigure[{\label{fig:P1h2}
	Timestep is $\tau=0.1\,h^2$.  This choice leads to
	$\EOC[\Th_m]\to2$ and $\EOC[E_m]\approx1$, i.e., the time
	indicator $\Th_m$ is of higher order than the spatial
	indicator $E_m$ which leads the estimator's order.  Thus we
	obtain asymptotic exactness $\EI\to1$, as expected from ZZ
	estimators for $p=1$.
    }]{
      \includegraphics[scale=\figscale,width=\figwidth]{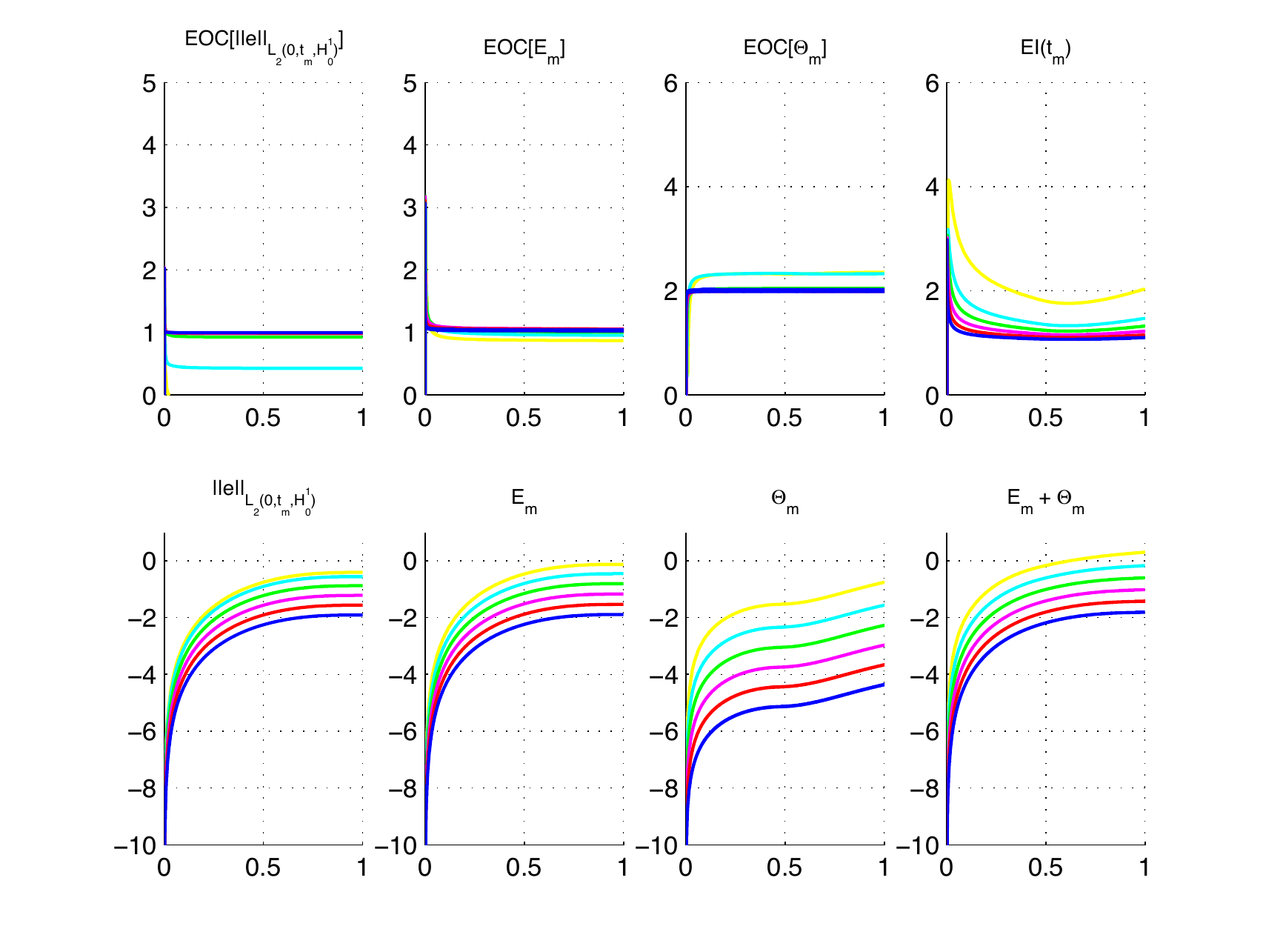}
    }
  \end{center}
\end{figure}
\begin{figure}[ht]
  \caption{\label{fig:P2}
    Numerical Results for \eqref{eqn:Problem:1} with $\poly2$
    elements and $h=h(i)=2^{-i/2}$ with $i=3,\dotsc,8$.  We compute
    the same quantities as in Figure \ref{fig:P1}.
  }
  \begin{center}
    \subfigure[{\label{fig:P2h2}
	Timestep $\tau=0.1\,h^2$.  The cumulative time error
	indicator $\Th_m$ is dominant with $\EOC[\Th_m]\to2$, but
	$\EI\gg1$.
    }]{
      \includegraphics[scale=\figscale,width=\figwidth]{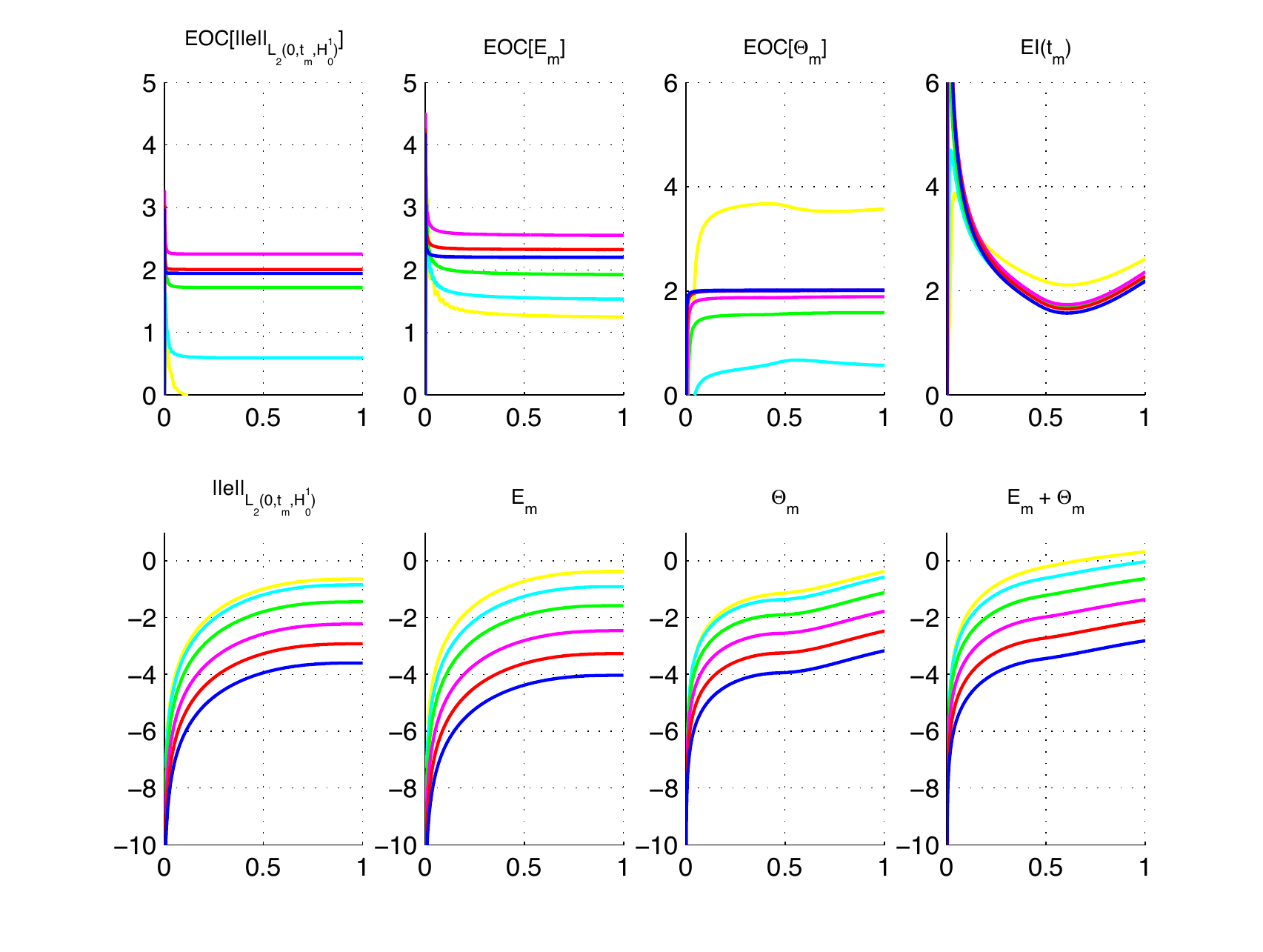}
    }
    \\
    \subfigure[{\label{fig:P2h3}
	Timestep is $\tau=0.1\,h^3$, with In the bottom set of
	results the spatial is dominant ($\EOC\approx2$) showing the
	estimator is sharp and reliable for higher order polynomials
	as well, and close to asymptotically exact ($\EI$ just
	smaller than $1$).
    }]{
      \includegraphics[scale=\figscale,width=\figwidth]{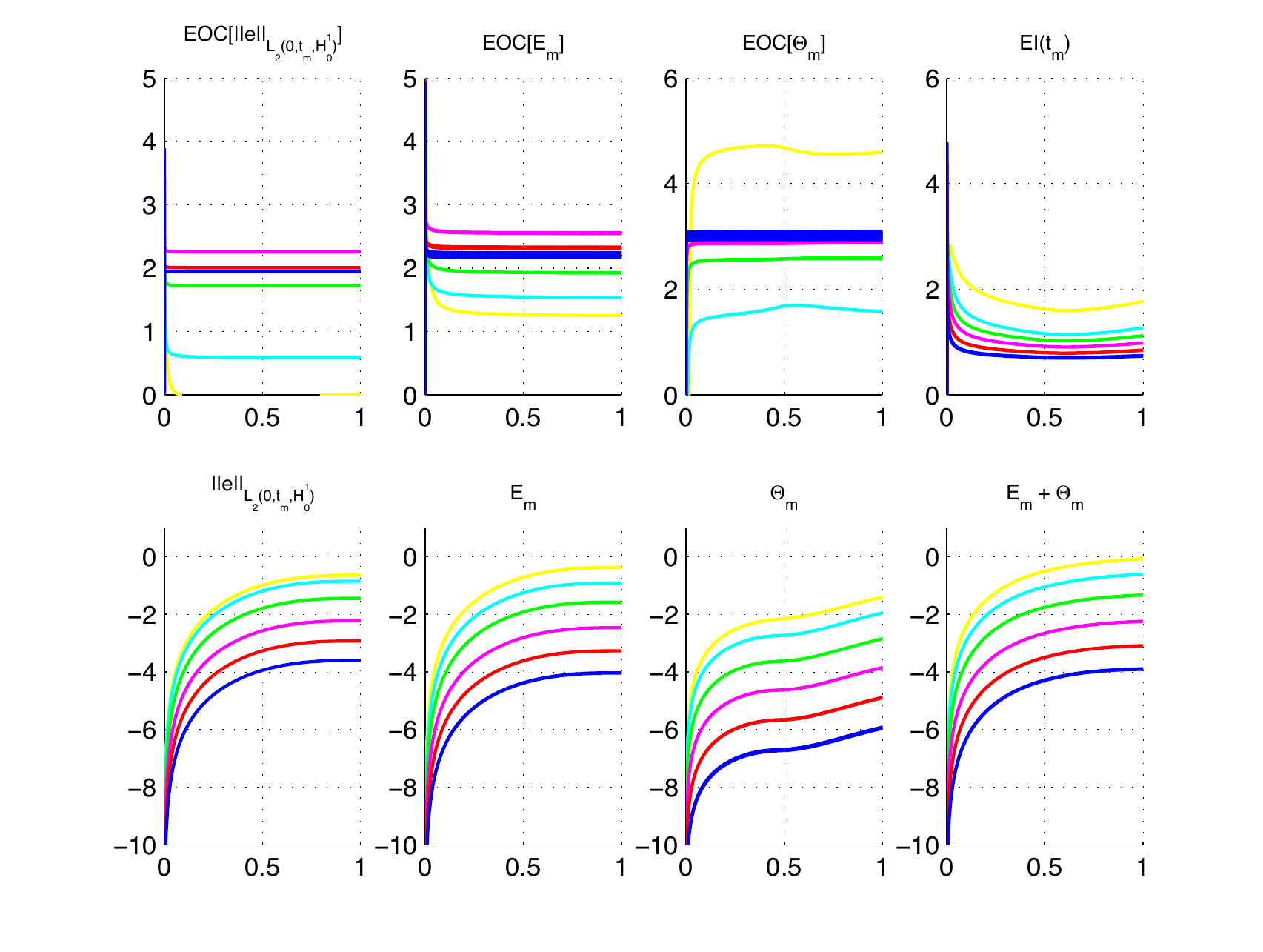}
    }
  \end{center}
\end{figure}
\begin{figure}[ht]
  \caption{\label{fig:P3}
    Numerical Results for \eqref{eqn:Problem:1} with $\poly3$
    elements for mesh-sizes $h(i) = 2^{-i/2}$, $i = 2,\dotsc,6$.. We
    compute the same quantities as in Figures \ref{fig:P1} and
    \ref{fig:P2}.
  }
  \begin{center}
    \subfigure[{\label{fig:P3h3}
	Timestep is $\tau=0.1\,h^3$. Again, the time indicator is
	dominant and $\EOC[\Th_m]\to3$, but $\EI\gg1$.
    }]{
      \includegraphics[scale=\figscale,width=\figwidth]{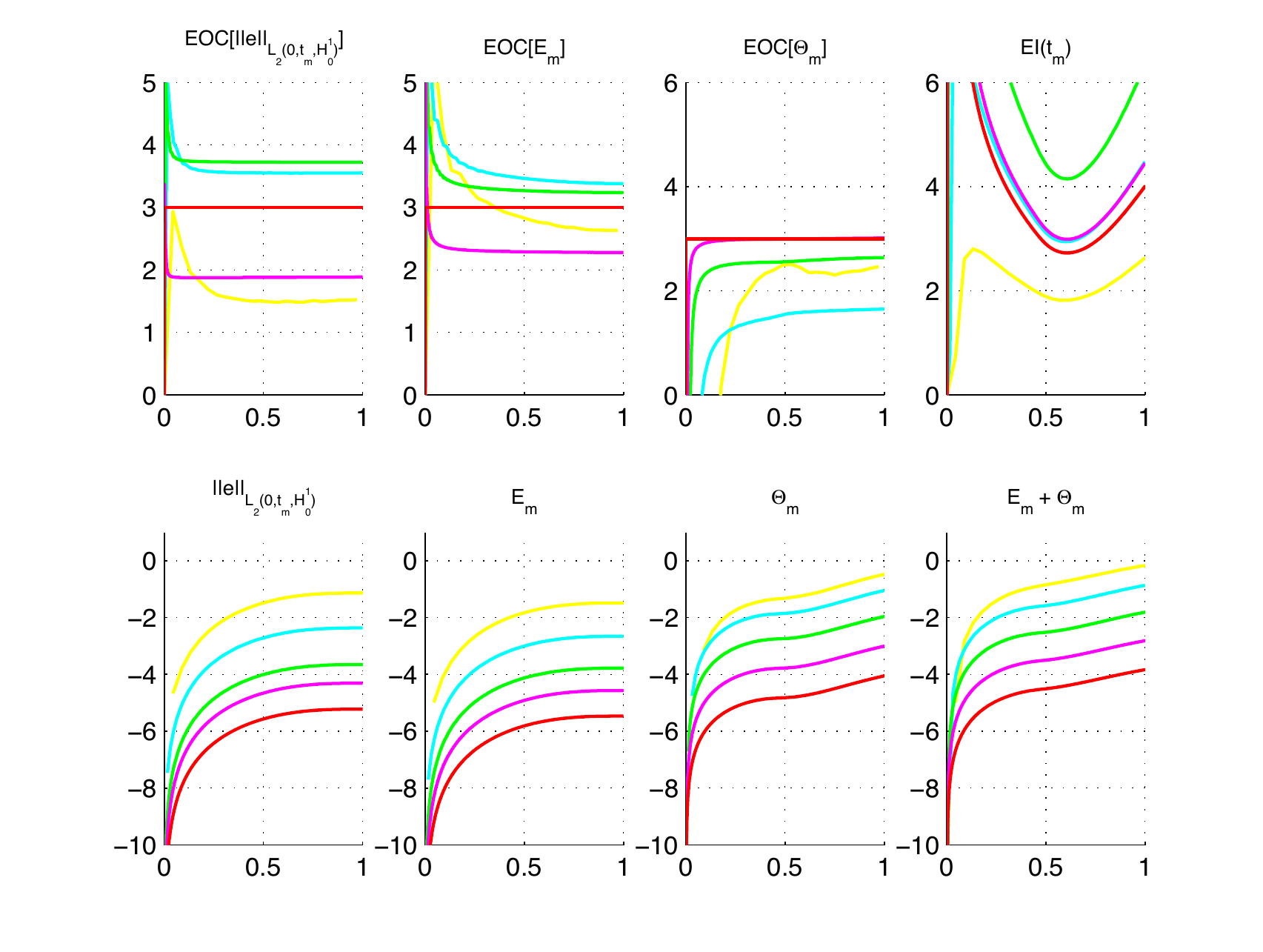}
    }
    \\
    \subfigure[{\label{fig:P3h4}
	Timestep is $\tau=0.1\,h^4$.  The elliptic error is dominant
	($\EOC[E_m]\to3$) and the estimator is sharp and reliable
	with very good $\EI$.
    }]{
      \includegraphics[scale=\figscale,width=\figwidth]{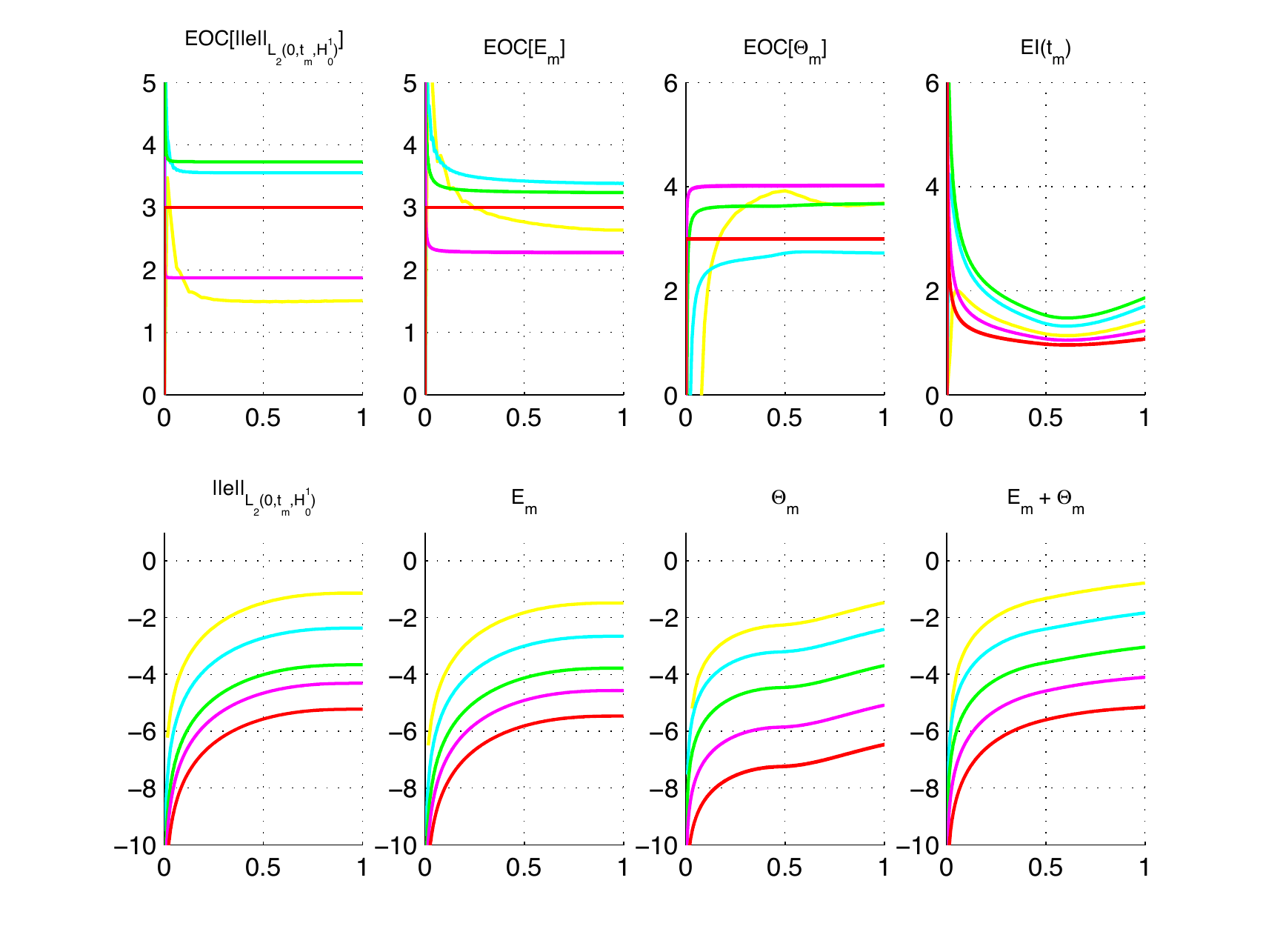}
    }
  \end{center}
\end{figure}
\begin{figure}[ht]
  \caption{\label{fig:P4}
    Results for \eqref{eqn:Problem:1} with $\poly4$
    elements and $h(i) = 2^{-i/2}$, $i = 2,\dotsc,6$.  We
    compute the same time accumulation quantities as in Figures
    \ref{fig:P1}--\ref{fig:P3}.
  }
  \begin{center}
    \subfigure[{\label{fig:P4h4}
	Mesh-size is $\tau=0.1\,h^4$.  Again, the time indicator is
	dominant with order $\EOC[\Th_m]\to4$) and a quite good
	$\EI$ in this case.
    }]{
      \includegraphics[scale=\figscale,width=\figwidth]{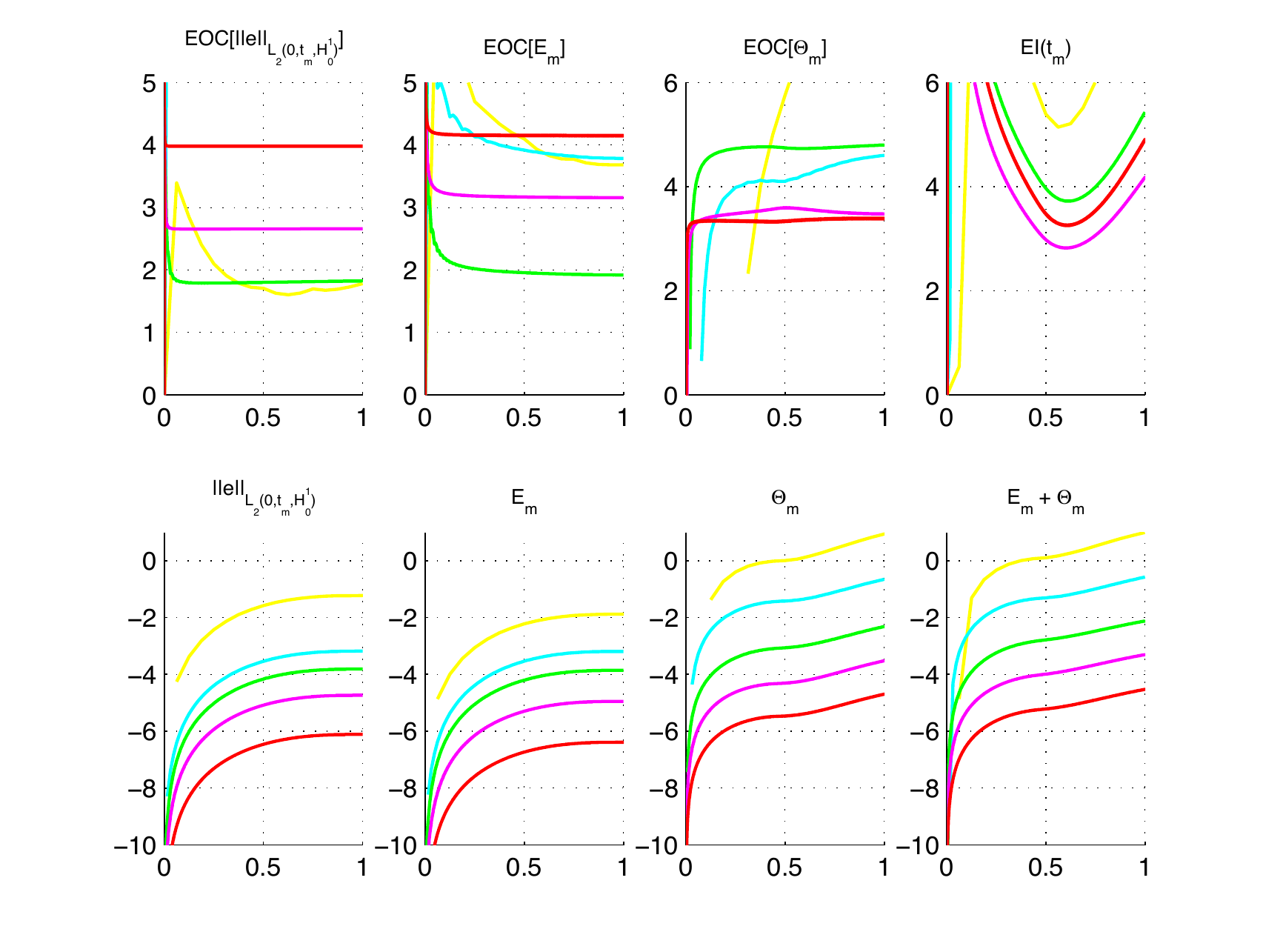}
    }
    \\
    \subfigure[{\label{fig:P4h5}
	Mesh-size is $\tau = 0.1 h^5$. The spatial error is dominant
	and $\EOC[E_m]\to4$.  Effectivity index improves slightly
	over previous case.
    }]{
      \includegraphics[scale=\figscale,width=\figwidth]{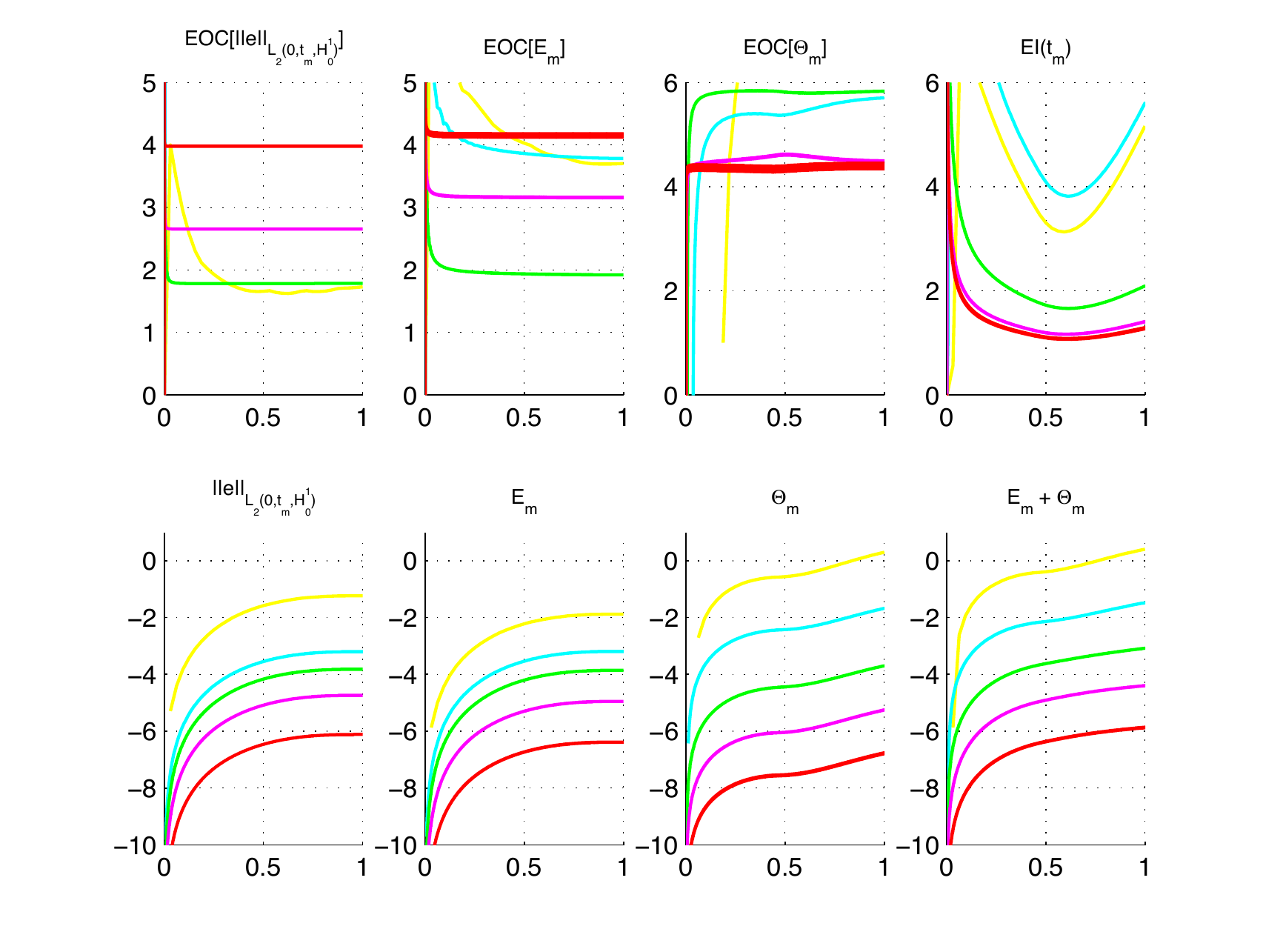}
    }
  \end{center}
\end{figure}
\begin{figure}[ht]
  \caption{\label{Fig:comparetheta}
    For each $m=1,\dotsc,N$ we plot values and \EOC's of two
    alternative time indicators
    $\psqrt{\sum_{n=1}^m\tau_n\tilde\theta_n^2}$ (above) and
    $\psqrt{\sum_{n=1}^m\tau_n\theta_n^2}$ (below) and the
    alternative mesh-change indicator
    $\sum_{n=1}^m\tau_n\tilde{\gamma}^2_n$ (above-right).  All
    quantities are plotted against time.  We took a uniform timestep
    $\tau =0.1\,h$ and mesh-size $h=2^{-i}$, $i = 4,\dotsc,9$. The
    numerical results show (1) that the two time indicators are equivalent in
    order, as expected, and (2) that the term
    $\sum_{n=1}^m\tau_n\tilde{\gamma}^2_n$ is indeed a higher order
    term and can be safely ignored in most practical schemes.  The
    indicators $\tilde\theta_n$ have a better effectivity index.
  }
  \begin{center}
    \includegraphics[scale=\figscale,width=\figwidth]{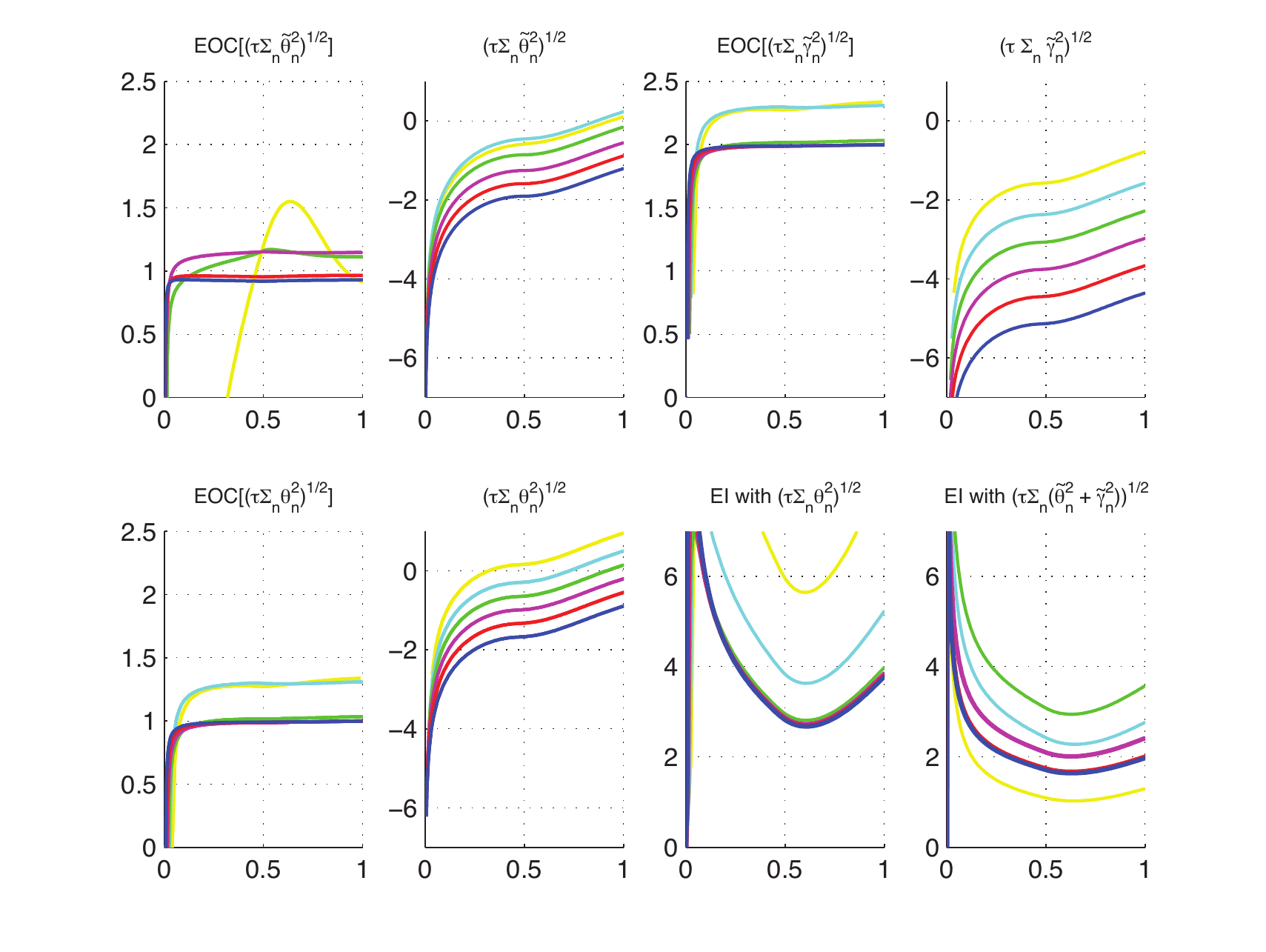}
  \end{center}
\end{figure}
\section{Computer experiments: adaptive schemes}
\label{sec:adaptivity}
We present now an adaptive algorithm based on the error indicators
defined in~\S\ref{def:fully-discrete:estimators}.  
\changes{
  As with many
  adaptive methods for time-dependent
  problems~\cite[]{picasso:98,alberta,chen-jia:03}, we perform space and
  time adaptivity separately.  Adaptivity is controlled via the
  indicators $\eta_n$ and $\eta_n$ (or $\tilde\eta_n$)---see
  Theorems~\ref{the:fully-discrete-aposteriori-bound}
  and~\ref{the:fully-discrete-aposteriori-bound:alternative}---which
  are kept under a given tolerance $\tol$.}

Namely, at each timestep $\tno\to\tn$, we use adaptive schemes for
elliptic problems as to minimise the indicators $\tilde\ep_n$ and
$\beta_n$.  There are different strategies to perform the timestep
adaptivity, all geared towards minimising $\theta_n$ (or
$\var\theta_n$).  Finally, the coarsening estimator $\gamma_n$ is
minimised by precomputing it and performing only one coarsening
operation at the beginning of each timestep.

Note that it is not in the scope of this paper to prove any rigorous
result about the adaptive algorithm and, based on heuristic
arguments only, we use it for illustration purposes.
\subsection{Space adaptivity via maximum strategy}
At each timestep an elliptic problem is solved.  For linear
elliptic problems, convergence of adaptive schemes is reasonably well
understood~\cite[]{morin-nochetto-siebert:2002:review,Binev:2004}
so we follow the criteria given therein, namely the Maximum Strategy.

The algorithm we used can be pseudocoded as follows. 
\subsection{{$\Algoname{Space Adapt}$}}
\label{alg:one-time-step}
\begin{algorithmic}
    \Require $(U^{\old},\V{\old},\tol_\ep, k_{\max},t,\tau,\threshhold,\tol_\g)$
    \Ensure $(U^{\new}, \V{\new})$ solution of \eqref{eqn:fully-discrete}
    
    \Procedure{Coarsening}{}
    
    \State $\numvec\gamma=(\gamma^K)_{K\in\cT}
    \assignvalue\Algoname{Coarsening Preindicator}(U\old,\V\old)$ 
    (cf.~\S\ref{alg:coarsening-error}).
    
    \State $\cT \assignvalue\Algoname{Mesh}(\V\old)$
        
    \State find $\cC \subset \cT$ such that $\sum_{K\in\cC}\smash{(\gamma^K)}^2 \leq \tol_\g^2$
    
    \State $\cT \assignvalue\Algoname{Coarsen}(\cT ,\cC)$ using \cite[\S1.1.2--1.1.3]{alberta}
    \EndProcedure

    \Procedure{Maximum Strategy Refinement~\cite{alberta}}{}
    \State $k\assignvalue0$
    \State compute $\ep_n$ using~(\ref{recoveryestimator})
    \State $\cR\assignvalue\emptyset$ \Comment{refinement set}

    \While{$\varepsilon_n > \tol_\ep$ and $k\leq k_{\max}$}
    
    \For{$K\in\T{n}$}
    
    \If{$\varepsilon_{K,{n}}^2 \geq \threshhold \max_{L\in\T{n}}{\varepsilon_{L,n}^2}$}
    
    \State $\cR\assignvalue \setof K\join\cR$ \Comment{mark $K$ for refinement}
    
    \EndIf
    
    \EndFor
    
    \State $\cT :=\Algoname{Refine}(\cT ,\cR)$ using \cite[\S1.1.1]{alberta} \Comment{hence update $(U^\old,\fes)$}
    
    \State set $\laginterpol n U^{n-1}:=U^\old$, $\taun=\tau$, $\tn=t$ and solve for $U^n$ in~(\ref{eqn:fully-discretepw}) 

    \State $U:=U^n$
    
    \State compute $\ep_n$ using~(\ref{recoveryestimator})
    
    \State $k\assignvalue k+1$
    \EndWhile
    \EndProcedure
    \State return $(U,\fes)$
  \end{algorithmic}
\subsection{Coarsening}
\label{coarsening}
In time-dependent problems mesh coarsening , which is not to be
confused with the coarsening needed in proving optimal complexity
for adaptive schemes \cite[]{Binev:2004}, is used to reduce DOF that
become redundant in time.

Mesh coarsening is a delicate procedure and should be used sparingly
as to avoid needless overhead computing time.  In
Algorithm~\ref{alg:one-time-step}, coarsening is performed only
once, at the beginning, for each time-step.  

The coarsening strategy we propose is based on \emph{predicting the
  effect of a possible removal of degrees of freedom}.  The reason for
this is that in \alberta (and many other finite element codes) upon
coarsening, all DOF-dependent vectors (encoding finite element
function coefficients) are ``coarsened'' via interpolation.  This
makes it possible to compute the effect of coarsening, and the
coarsening estimator $\gamma_n$ defined in
(\ref{eqn:def:coarsening-indicator}), \emph{before} mesh-change
occurs.  The details of this procedure are discussed
in \S~\ref{sec:building-coarsening-estimator}.
\subsection{Timestep control}
Timestep control can be achieved using two different strategies.

An \emph{implicit timestep control} strategy used is ready
implemented in \alberta \cite[]{alberta} using Algorithm
\ref{alg:one-time-step} upon each timestep.

Here we propose an \emph{explicit timestep control} strategy which
we have implemented in \alberta.  The reason for this is that the
implicit strategy, though better in terms of timestep determination,
is very time-consuming as it requires the repeated solution of the
timestep.  In contrast, the explicit strategy has a
rougher---nonetheless still satisfactory--- control over the
timestep, but it is much faster.  The conclusion is that the ideal
control should be a smart implicit/explicit-switching algorithm.

The explicit strategy can be described as follows. 
\subsection{$\Algoname{Explicit Timestep Adapt}$}
\label{alg:time}
\begin{algorithmic}
  
  \Require $(\tau_0,t_0,T,\T0,u^0,\tol_\ep,k_{\max},\threshhold,\tol_\g,\tol_{\theta,\min},\tol_{\theta})$
  
  \Ensure $(\taun,\V n,U^n)_{\rangefromto n1N}$ satisfying (\ref{eqn:fully-discrete}) 
  and possibly $\int_0^T\Norm{U-u}^2\leq\tol^2$
  
  \State $(U^0, \V0) = \Algoname{Initial Space Adapt}
  (\T0,u^0,k_{\max},\threshhold,\coarsethresh)$
  \Comment{data interpolation}
  
  \State $n:=1$
  
  \State $\tau_n := \tau_{n-1}$
  
  \State $t_n := t_{n-1} + \tau_n$
  
  \While{$t_n\leq T$}
    
  \State $(U^n,\V{n}) := 
  \Algoname{Space Adapt}
  (U^{n-1},\V{n-1}, \tol_\ep, k_{\max}, \tau_n, t_n, \threshhold, \tol_\g)$
  
  \State compute $\theta_n$ 
  
  \If{$\theta_n > \tol_{\theta}$}
  
  \State $\tau_{n+1} := \tau_n/\sqrt{2}$
    
  \ElsIf{$\theta_n \leq \tol_{\theta,\min}$}
  
  \State $\tau_{n+1} := \sqrt{2}\tau_n$
  
  \EndIf
  
  \State $t_{n+1} := t_n + \tau_{n+1}$
  
  \State $n := n+1$
  
  \EndWhile
  
  \State return $(U^n)_{\rangefromto n1N}$,
\end{algorithmic}
where the \emph{global tolerance} $\tol$ is given by the relation
\begin{equation}
  \tol^2=T\qp{\tol_\theta^2+\tol_\ep^2+\tol_\g^2}.
\end{equation}
Note that this algorithm does not guarantee reaching a tolerance,
unlike more sophisticated ones found in the
literature~\cite[e.g.]{chen-jia:03}, but it guarantees termination
in reasonable CPU times.
\subsection{Numerical results}  
In Tables \ref{Tbl:Prob1}--\ref{Tbl:Prob3} we compare the implicit
timestep control strategy described by algorithm \ref{alg:time} with
a uniform timestep scheme.  For the uniform strategy we take a
stationary mesh in time and set $\tau = 0.04 h^2$.  We calculate the
error for various numerical simulations using differing values of
$h$ using the uniform strategy and set those values as tolerances
for the adaptive scheme varying $\threshhold$ appropriately.

Each column displays results for either the uniform strategy or the
adaptive strategy using various thresholds. These columns are
further subdivided into two, the first containing
$\sum_{n=1}^N\dim{\V{n}}$ (\ie the total number of degrees of
freedom from all meshes over time) which we denote DOF and the
second containing CPU time (secs) for all model problems
\eqref{eqn:Problem:1}--\eqref{eqn:Problem:2}.
\begin{table}[ht]
  \begin{center}
    {\tiny
    \begin{tabular}{|c|c|c|c|c|c|c|c|c|}
      \hline
        & \multicolumn{2}{c|}{Uniform} & \multicolumn{6}{c|}{Adaptive}\\
      \hline
      & \multicolumn{2}{c|}{  } & \multicolumn{2}{c|}{$\threshhold = 0.65$} 
      & \multicolumn{2}{c|}{$\threshhold = 0.70$} 
      & \multicolumn{2}{c|}{$\threshhold = 0.75$} \\
      \hline
      $\tol$ & DOF's & CPU   &DOF's &CPU  &DOF's &CPU  &DOF's&CPU \\
      \hline
      0.573    & 232,290 &3  &24,080 & 4  &22,792 &5  &22,240& 4 \\
      \hline
      0.295   & 3,489,090 &49  & 42,042 &8 & 39,414 & 8 & 38,630&6 \\
      \hline
      0.149  &54,097,020& 598 &82,172&15 &77,932 &15 &76,452&16 \\
      \hline
      0.0625 & OOM   &OOM &206,709&39&195,810&37&191,650&37 \\
      \hline
    \end{tabular}
    }
  \end{center}
  \caption{
    \label{Tbl:Prob1}
    Explicit timestep control with various spatial maximum strategy
    thresholds for Problem \eqref{eqn:Problem:1}.  The adaptive
    method clearly saves DOF and CPU time over the uniform
    method.}
\end{table}
\begin{table}[ht]
  \begin{center}
    {\tiny
    \begin{tabular}[subsection]{|c|c|c|c|c|c|c|c|c|c|c|}
      \hline
      & \multicolumn{2}{c|}{Uniform} & \multicolumn{6}{c|}{Adaptive}\\
      \hline
      & \multicolumn{2}{c|}{  } & \multicolumn{2}{c|}{$\threshhold = 0.65$} & \multicolumn{2}{c|}{$\threshhold = 0.7$} & \multicolumn{2}{c|}{$\threshhold = 0.75$}\\
      \hline
      $\tol$      & DOF's  & CPU    & DOF's & CPU & DOF's & CPU & DOF's & CPU\\
      \hline
      0.296      & 3,489,090  & 47    & 12,092& 5 & 11,430 & 5 & 11,498 & 5\\
      \hline
      0.21     & 13,940,289 & 196   & 17,038 & 7 & 16,140 & 8 & 16,201 & 7\\
      \hline
      0.104    & 54,097,020 & 602  & 106,188 & 32 & 100,058 & 29 & 22,597 & 10\\
      \hline
      0.03125   & OOM& OOM & 513,694 & 120 & 460,637 & 118 & 449,568 & 115\\
      \hline
    \end{tabular}
    }
  \end{center}
  \caption{
    \label{Tbl:Prob2}
    Explicit timestep control with various spatial maximum strategy
    thresholds for spatial-error dominant Problem
    \eqref{eqn:Problem:2}.  Adaptivity saves DOF and CPU.}
\end{table}
\begin{table}[ht]
  \begin{center}
    {\footnotesize
      \begin{tabular}[subsection]{|c|c|c|c|c|c|c|c|c|}
        \hline
        & \multicolumn{2}{c|}{Uniform} & \multicolumn{4}{c|}{Adaptive}\\
        \hline
        & \multicolumn{2}{c|}{  } & \multicolumn{2}{c|}{$\threshhold = 0.7$} & \multicolumn{2}{c|}{$\threshhold = 0.75$}\\
        \hline
        $\tol$ &DOF's & CPU & DOF's & CPU & DOF's & CPU\\
        \hline
        1.000         &925,809 & 12 & 159,070& 43 & 127,610&58\\
        \hline
        0.569     &3,489,090 & 49& 237,960 & 142 & 204,376&180 \\
        \hline
        0.295     &54,097,020& 605 & 471,733&755& 471,542 &920\\
        \hline
        0.149     &OOM&OOM& 940,618 & 1410 & 940,138& 1850\\
        \hline
      \end{tabular}    
    }
  \end{center}
  \caption{
    \label{Tbl:Prob3}
    Implicit timestep control with various spatial maximum strategy
    thresholds for spatial-error dominant Problem
    \eqref{eqn:Problem:3}.  Adaptivity saves DOF (even better than
    explicit control) but the CPU time grows very quickly due to
    overhead.}
\end{table}
\begin{table}[ht]
  \begin{center}
    {\footnotesize
      \begin{tabular}{|c|c|c|c|c|c|c|c|c|}
        \hline
        & \multicolumn{2}{c|}{Uniform} & \multicolumn{4}{c|}{Adaptive}
        \\
        \hline
        & \multicolumn{2}{c|}{  } & \multicolumn{2}{c|}{$\threshhold = 0.7$} & \multicolumn{2}{c|}{$\threshhold = 0.75$}
        \\
        \hline
        $\tol$ &DOF's & CPU & DOF's & CPU & DOF's & CPU
        \\\hline
        1.000&925,809&12&135,788&5&127,004&4
        \\\hline
        0.569&3,489,090 & 49&198,628&7&194,311&8
        \\\hline
        0.295&54,097,026&605&397,716&15&395,876&16
        \\\hline
        0.149&OOM&OOM&2,177,666&79&2,079,081&76
        \\\hline
      \end{tabular}
    }
  \end{center}
  \caption{
    \label{Tbl:Prob3-explicit}
    Explicit timestep control with various spatial maximum strategy
    thresholds for time-error dominant Problem
    \eqref{eqn:Problem:3}
  }
\end{table}
\begin{Rem}[implicit timestep control on fast oscillating solutions]
  We take note of the CPU times from the results for Problem
  \eqref{eqn:Problem:3}. These show that implicit timestep control is
  undesirable for fast oscillating functions. This is because the
  timestep searching becomes computationally inefficient. Numerical
  simulations for an explicit timestep control strategy is given in Table
  \ref{Tbl:Prob3-explicit}. This algorithm is described in detail in
  the \alberta manual \cite[\S 1.5.4]{alberta} The results show
  although for a method with low tolerance we use more degrees of
  freedom we make a substantial gain on the CPU time.
\end{Rem}

We then fix a value of $\threshhold$ and compare an adaptive
strategy with uniform for a single value of $\tol$. This is to
illustrate how the number of degrees of freedom of the mesh change
over time, and how the implicit timestep control affects the timestep
size for all test problems in Figures~\ref{fig:adaptivity}.
\renewcommand{\figscale}{0.30}
\begin{figure}
  \caption{\label{fig:adaptivity}
    Adaptive (green) against uniform (red) degrees of freedom and
    timestep sizes.  In each pair of graphs we plot the (log of) the
    DOF against time on the left, and the timestep against time on
    the right.
  }
  \begin{center}
    \subfigure[{
	Implicit timestep control for Problem
	(\ref{eqn:Problem:1}).  The explicit timestep control
	yields the same results (but is much more CPU efficient),
	thus it is not shown.
    }]{
      \includegraphics[scale=\figscale]{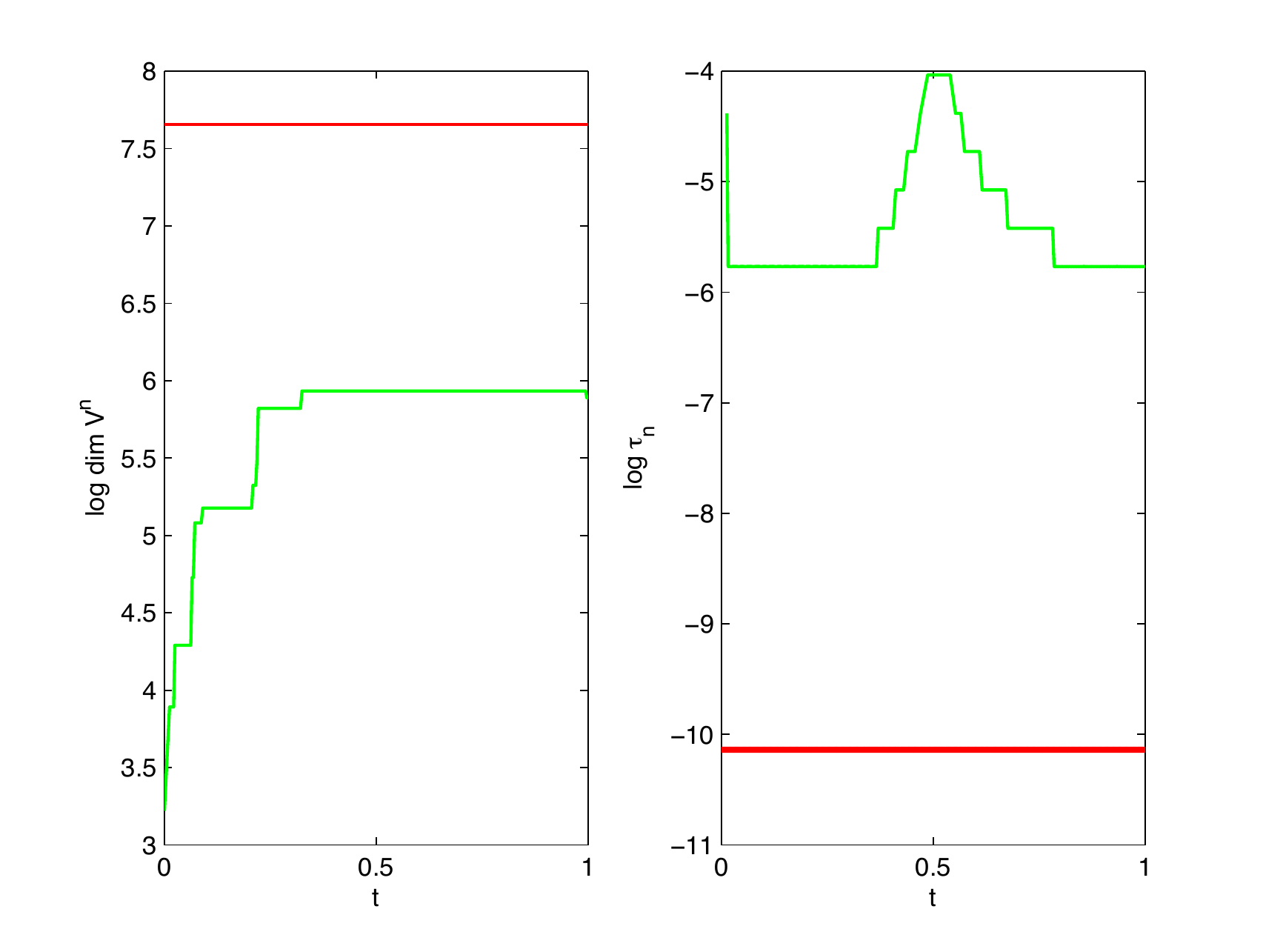}
    }
    \hfill
    \subfigure[{
	Implicit timestep control for Problem
	(\ref{eqn:Problem:3}), where the spatial error dominates.
	The explicit timestep control yields the same meshes and
	time-steps, thus not shown.
    }]{
      \includegraphics[scale=\figscale]{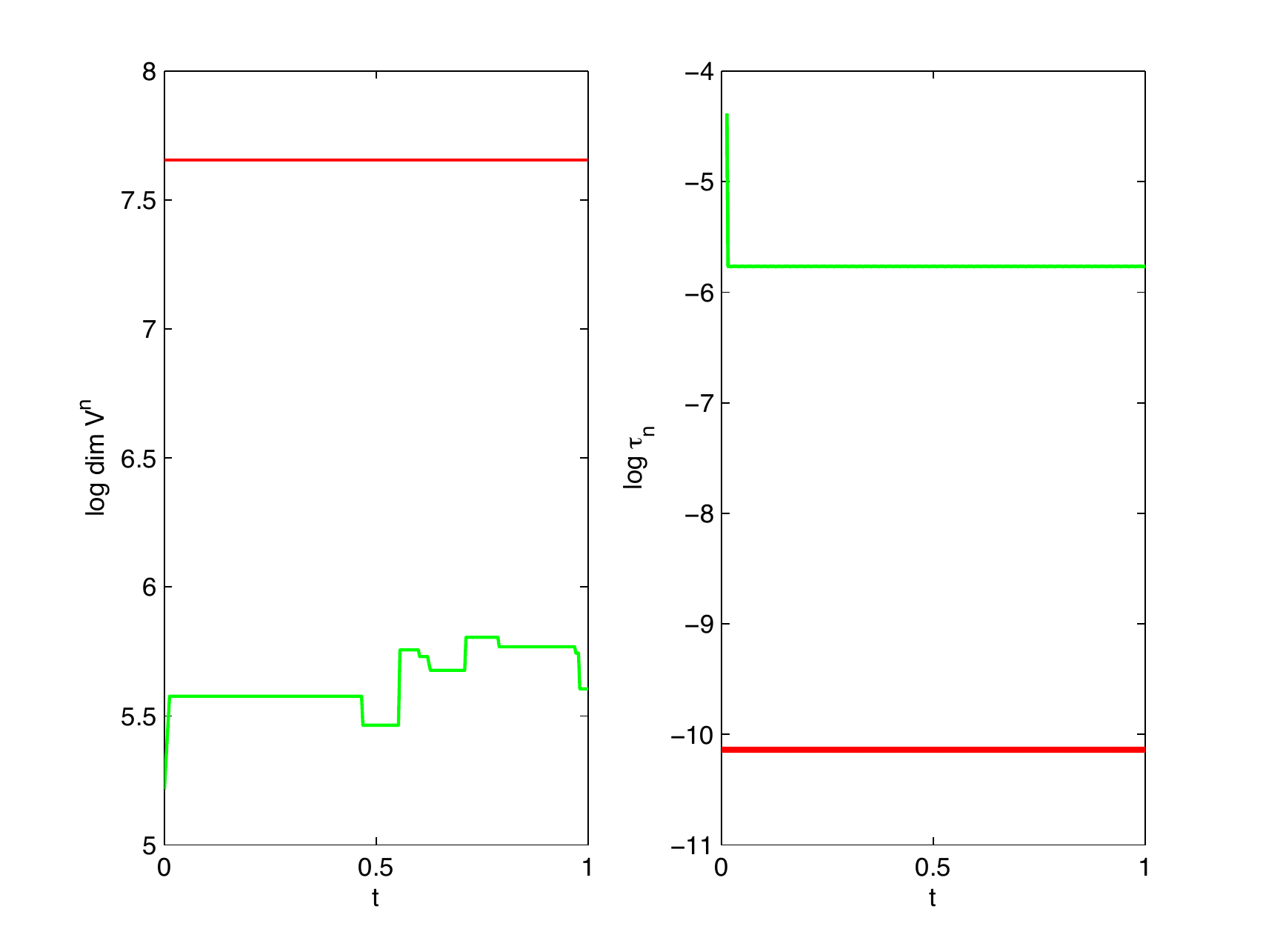}
    } 
    \\
    \subfigure[{
	\label{fig:adaptivity-time-explicit}
	Explicit timestep control for Problem
	(\ref{eqn:Problem:2}), where the time discretisation error
	dominates.  Interesting when compared with Figure \ref{fig:adaptivity-time-implicit}.
    }]{
      \includegraphics[scale=\figscale]{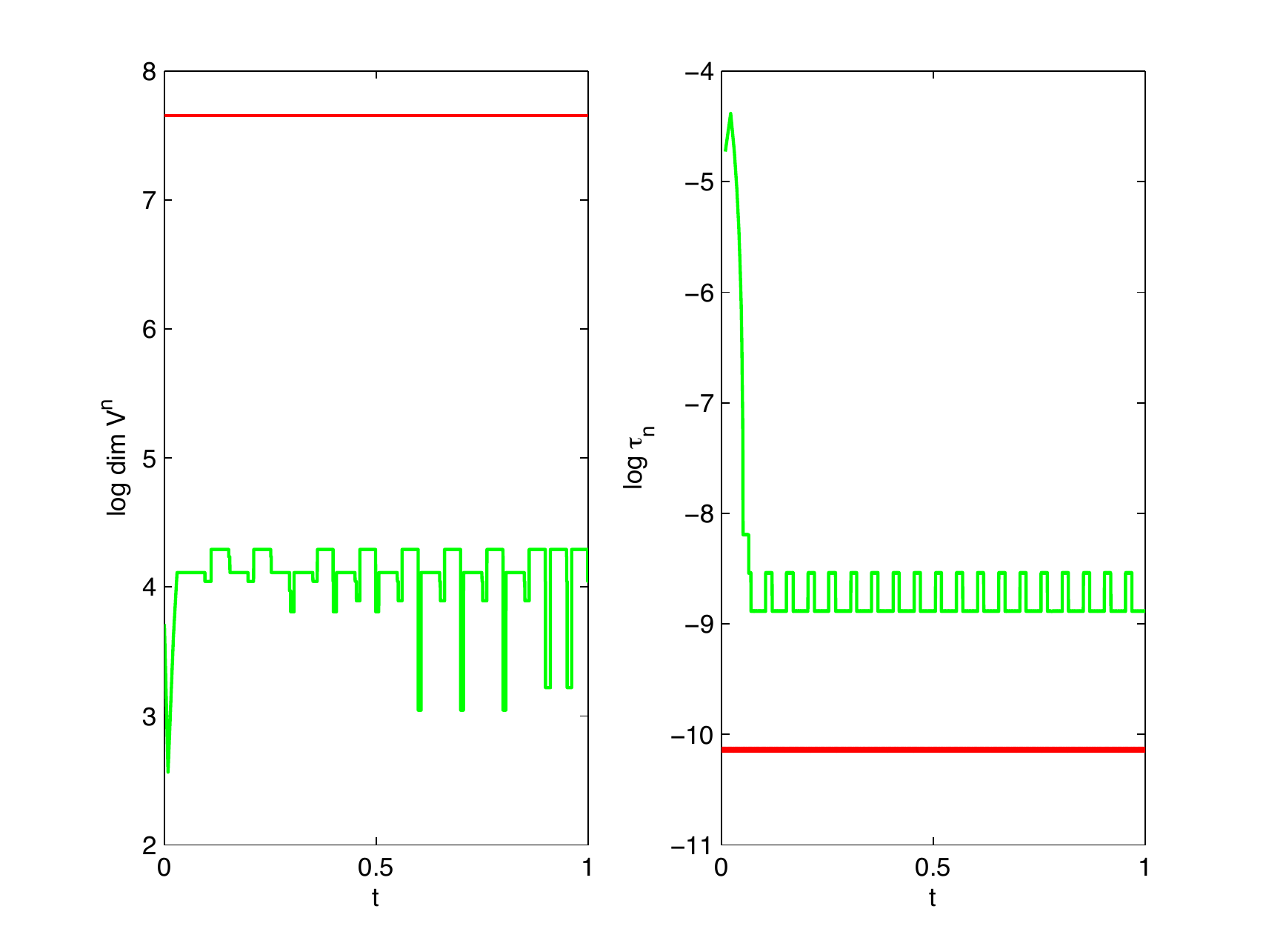}
    }
    \hfill
    \subfigure[{
	\label{fig:adaptivity-time-implicit}
	Implicit timestep control for Problem
	(\ref{eqn:Problem:2}).  Comparing with Figure
	\ref{fig:adaptivity-time-explicit} shows that the implicit
	timestep control yields more efficient timestep and
	meshes, but at a much higher CPU cost (cf. Tables
	\ref{Tbl:Prob3} and \ref{Tbl:Prob3-explicit}).
    }]{
      \includegraphics[scale=\figscale]{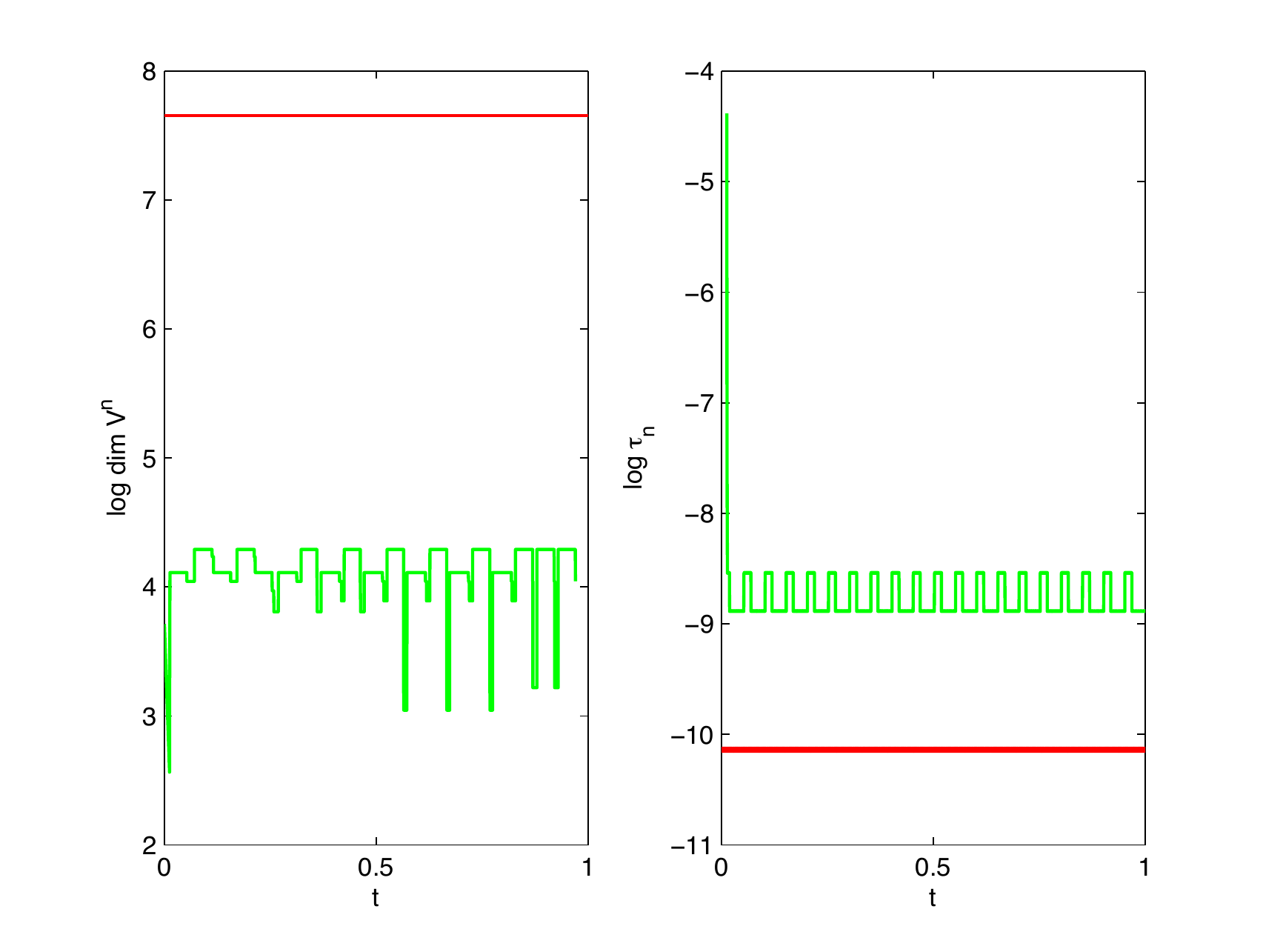}}
  \end{center}
\end{figure}

\subsection{Incompatible data singular solution}
\label{sec:incompatible-data}
\changes{ We close the paper by testing the adaptive algorithm on
  an example with incompatible initial and boundary conditions,
  which is the type of situation where adaptivity is really needed
  in practise. Consider problem (\ref{eqn:heat}) with $\W =
  (0,1)\times (0,1)$, $f=0$ and $u_0=1$. The initial conditions are
  thus incompatible with the homogeneous Dirichlet boundary valid for
  all positive times. The exact solution $u$, though singular at all
  points of $\boundary\W\times\setof0$, can be readily evaluated ``by
  hand'' and may be represented in terms of Fourier series of the
  Laplacian's eigenvalues. Namely, we have
  \begin{equation}
    \label{eq:incompat-ics}
    u(\geovec x, t) 
    =
    \sum_{m,n = 1}^\infty 
    C_{m,n} 
    \exp(-\left(m^2+n^2\right)\pi^2t)
    \sin(m\pi x_1)
    \sin(n\pi x_2)
    ,
    \text{ for }t>0,
  \end{equation}
  where the constant $C_{m,n}$ is given by
  \begin{equation}
    C_{m,n} = \frac{4}{n m\pi^2}
    \qp{
      1 - \cos\left(m\pi\right) - \cos\left(n\pi\right) 
      + \cos\left(n\pi\right)\cos\left(m\pi\right)
    }
    .
  \end{equation}
  Since the solution (\ref{eq:incompat-ics}) is an infinite
  Fourier series it cannot be computed exactly, but its rapid decay
  allows to truncate early with machine-epsilon precision. 
}

\changes{ 
  In order to generate a reference tolerance, which is common for
  the uniform and the adaptive scheme we couple $h=0.05\tau$ and run
  the uniform refinement code.  We use then the error computed as a
  tolerance for the adaptive scheme, results of this are shown in
  Figure \ref{Fig:incompat-adapt-uniform-dof-tau}.  In Figure
  \ref{Fig:incompatible} we visualise the
  adapted FE mesh for Problem (\ref{eq:incompat-ics}) at various
  times.
}

\renewcommand{\figscale}{0.50}      
\begin{figure}
  \caption{
    \label{Fig:incompat-adapt-uniform-dof-tau}
    Implicit timestep control for Problem
    (\ref{eq:incompat-ics}).}
  \begin{center}
    \includegraphics[scale=\figscale,trim=100 120 50 200, clip]{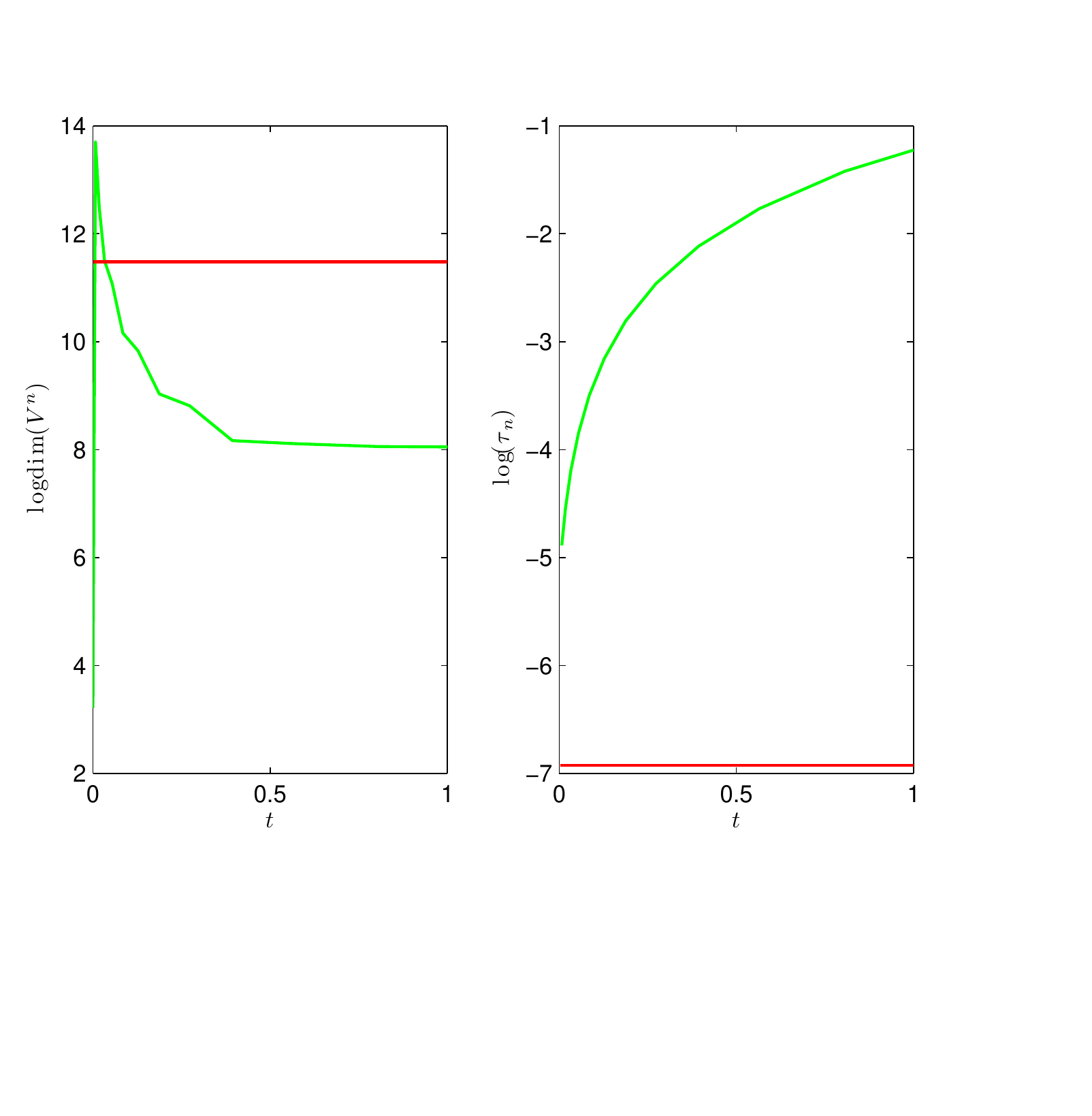}
  \end{center}
\end{figure}
\renewcommand{\figscale}{0.14}
\begin{figure}
  \caption{
    \label{Fig:incompatible}
    The adaptive scheme for (\ref{eq:incompat-ics}) using
    implicit timestep control.
  }
  \begin{center}
    \subfigure[{
        \label{Fig:incompat-01}
        Mesh at time $t_n=0.007544$ with
        $\dim(\V{n}) = 894,677$
    }]{
      \includegraphics[scale=\figscale]{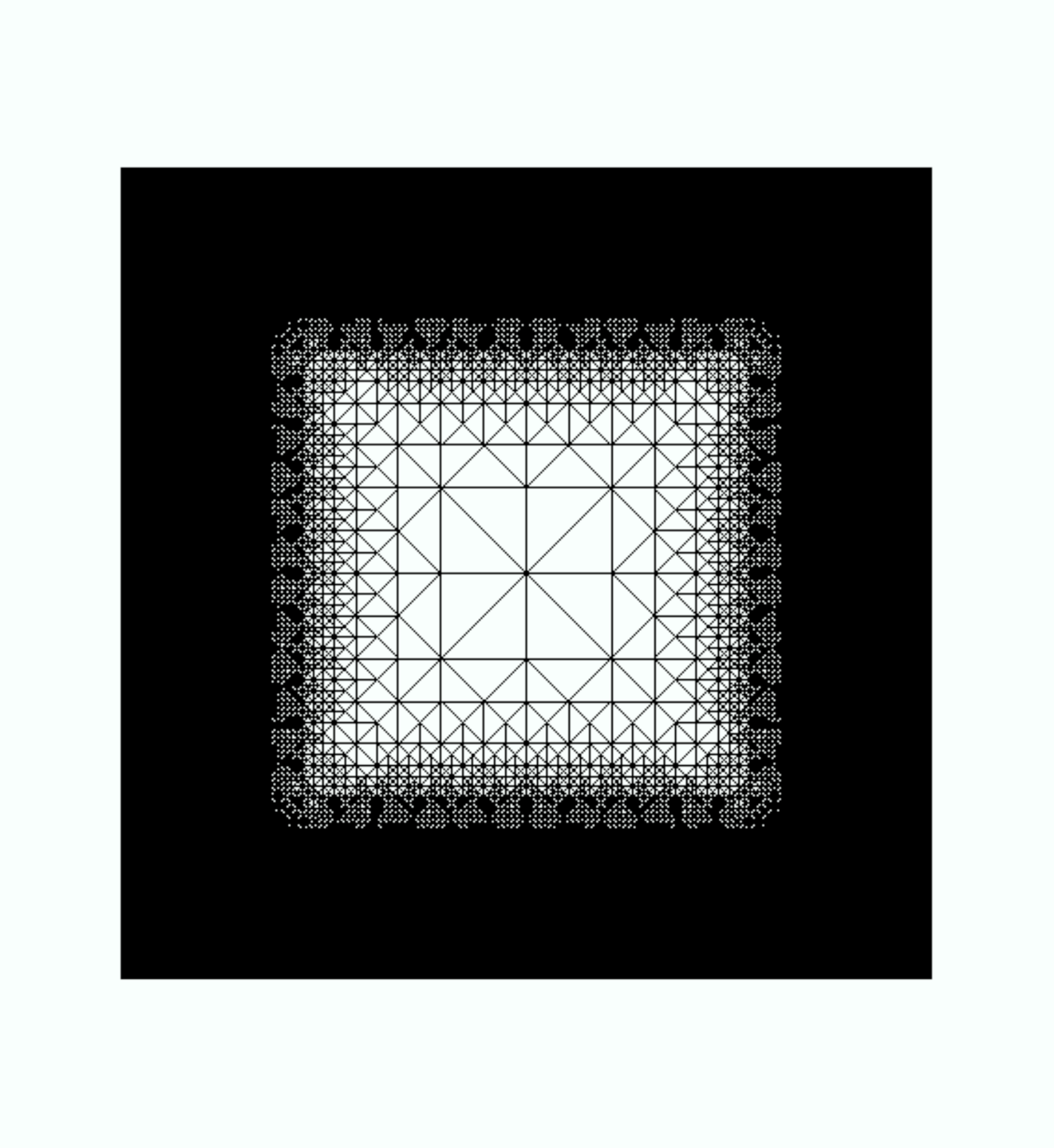}
    }
    \hfill 
    \subfigure[{
        \label{Fig:incompat-02}
        Mesh at time $t_n=0.033302$ with $\dim(\V{n}) = 98,773$ 
    }]{
      \includegraphics[scale=\figscale]{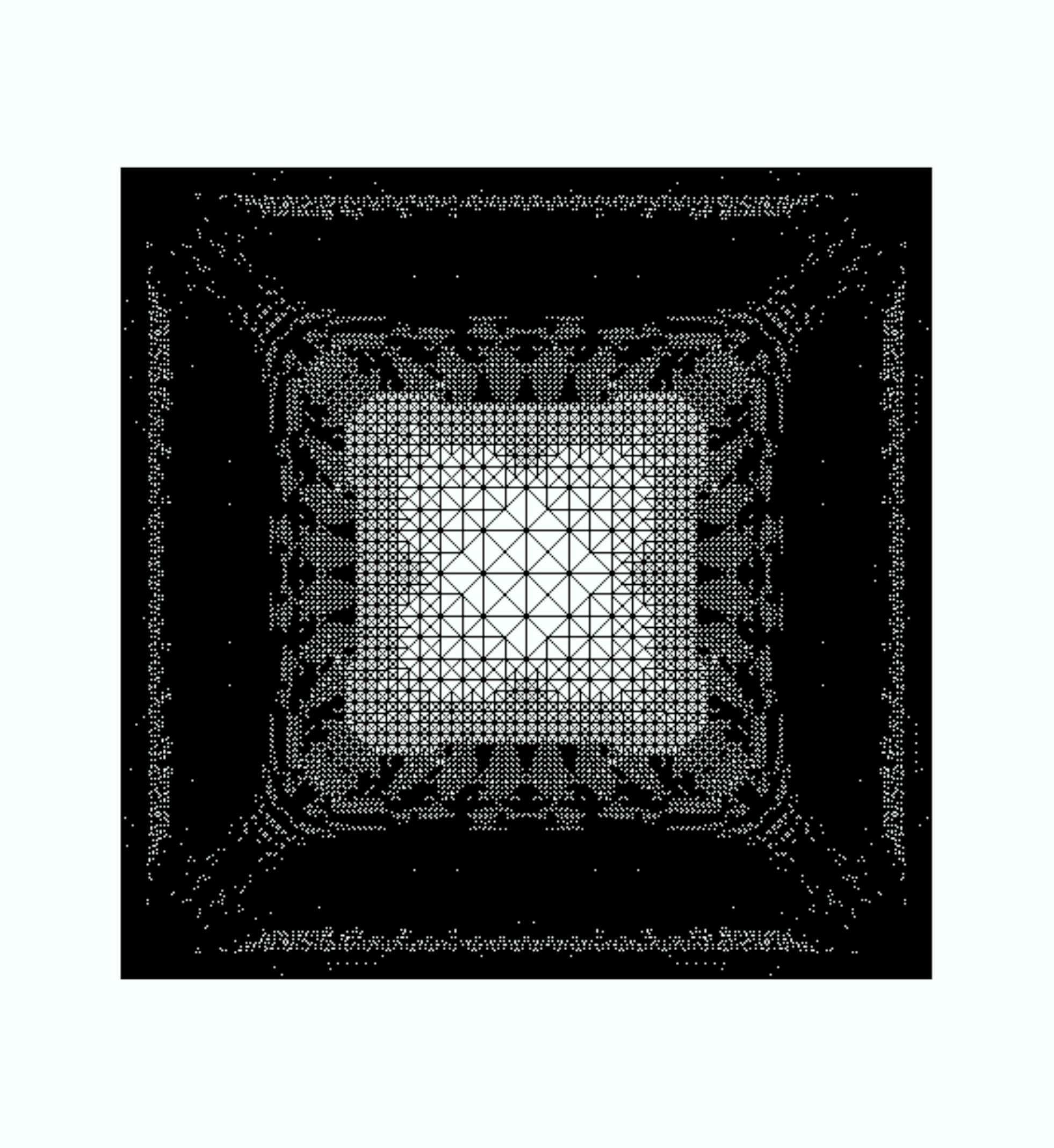}
    }
    \hfill 
    \subfigure[{
        \label{Fig:incompat-03}Mesh at time $t_n=0.127492$ with $\dim(\V{n}) = 18,613$ 
    }]{
      \includegraphics[scale=\figscale]{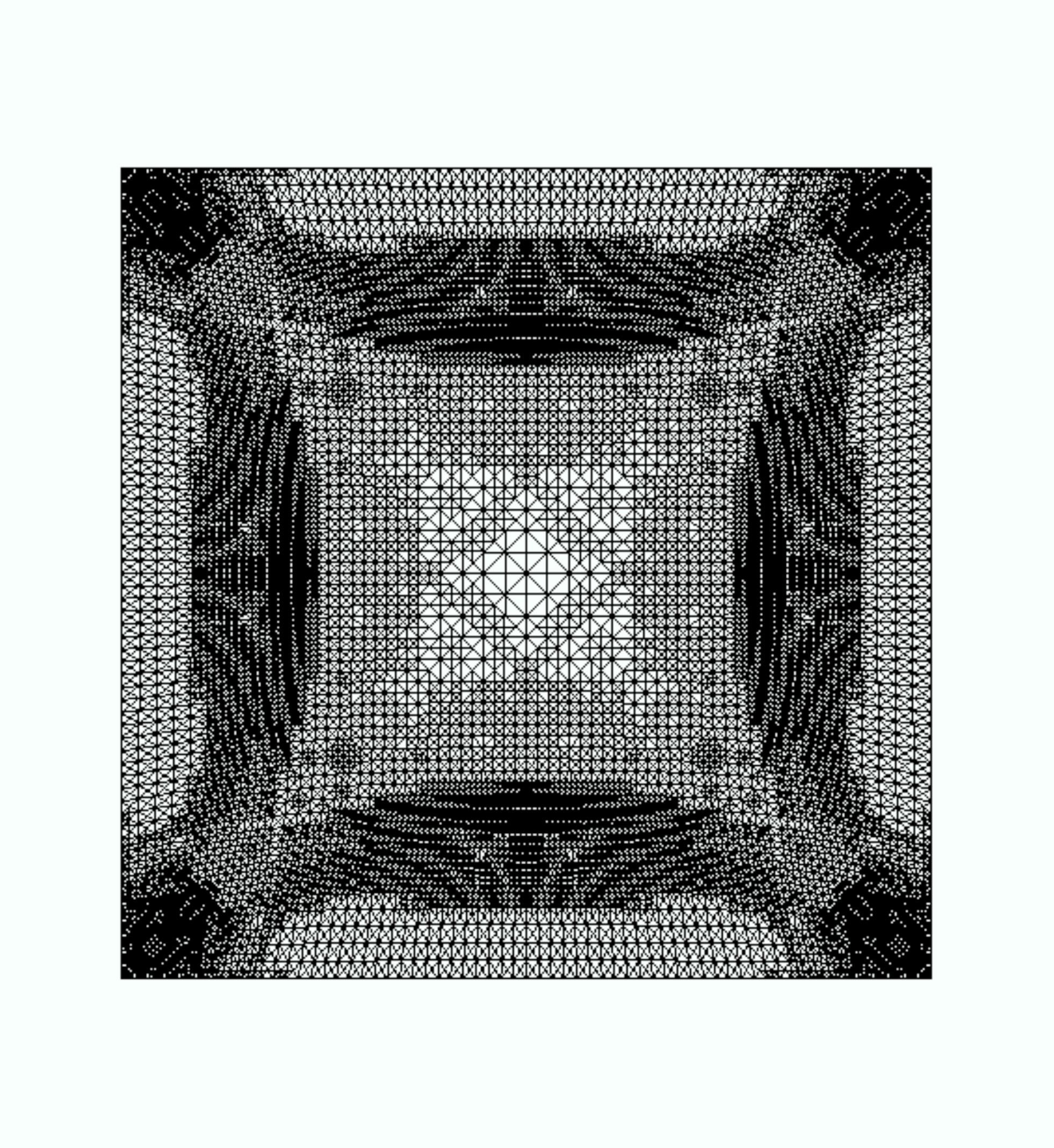}
    }
    \hfill 
    \subfigure[{
        \label{Fig:incompat-04}Mesh at time $t_n=0.393893$ with $\dim(\V{n}) = 3,525$ 
    }]{
      \includegraphics[scale=\figscale]{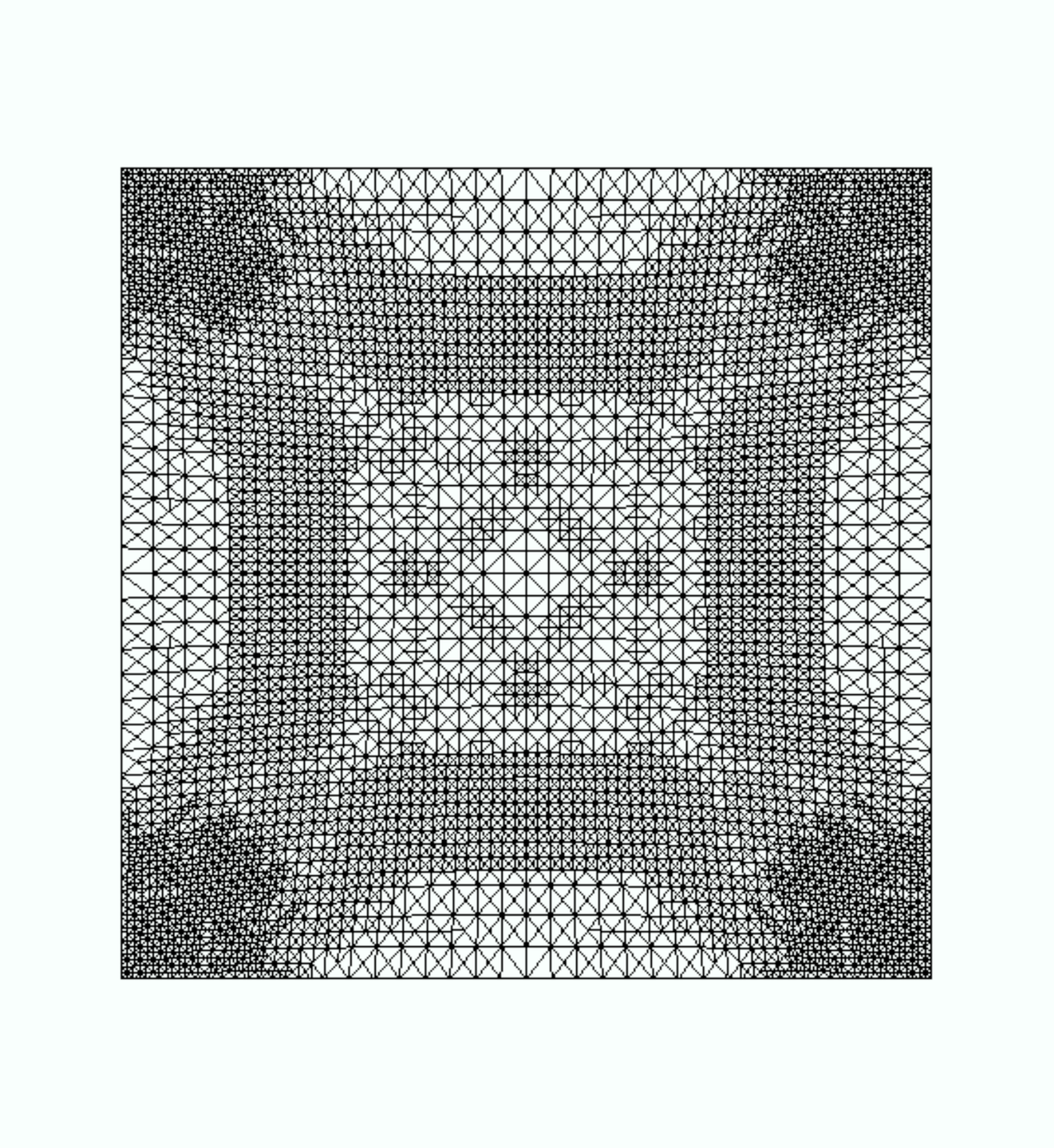}
    }
  \end{center}
\end{figure}
\section*{Acknowledgement}
O.L.'s research was partially supported by a Nuffield Young Researcher's Grant.
\\
T.P.'s research is fully supported by his EPSRC D.Phil. scholarship grant.
\\
Both authors wish to thank Alan Demlow for an interesting exchange of ideas.
  \appendix
  \section{Coarsening error preindicator implementation}
  \label{sec:building-coarsening-estimator}
  \newcommand{\coarse}{\mathbb{X}}
  \newcommand{\fine}[1]{\mathbb{Y}^{#1}}
  \newcommand{\getrowcol}[3]{\ensuremath{\smash{\qb{\smash{#1}}^{#2}_{#3}}}}
  \newcommand{\getcol}[2]{\ensuremath{\smash{\qb{\smash{#1}}_{#2}}}}
  \newcommand{\getrow}[2]{\ensuremath{\smash{\qb{\smash{#1}}^{#2}}}}
  We describe next a practical implementation of the \emph{coarsening
  error preindicator} (we use this term to emphasise the fact that
  this indicator can be computed \apriori, as opposed to the other
  indicators involved in the adaptive strategy).  Since we used
  \alberta for our computations, this section relies substantially on
  the principles described in the manual~\cite{alberta}. We briefly
  describe these principles in the next paragraph, in order to expose
  the main idea behind the coarsening preindicator.
  \subsection{Refinement, coarsening and interpolation in \alberta}
  Mathematically, a simplicial mesh (or partition, or triangulation)
  is a set of disjoint open simplexes, the union of the closure of
  which is $\closure\W$.  A mesh into a new mesh is refined by
  \emph{bisecting} a subset of its simplexes, following a special
  procedure which ensures mesh conformity (e.g., no hanging nodes) and
  does not deteriorate shape-regularity (on fully fitted polygonal
  domains).  A mesh is thus represented as a binary tree, where each
  node represents a simplex.  The children of each simplex are thus
  the $2$ subsimplexes obtained by bisection.  Hence, from a coding
  view-point, refinement means growing the binary tree.
  
  The inverse of refinement is coarsening.  Thus coarsening a mesh in
  \alberta consists in removing pairs of sibling simplexes (both
  marked for coarsening) and produces the new---coarsened---mesh where
  the pairs of siblings are replaced by their parent.
  
  The coarsening preindicator is a real number defined on each simplex, of
  the triangulation to be coarsened.  This estimator can in fact be
  \emph{precomputed} with respect to coarsening.  This is in contrast
  with usual \emph{\aposteriori} error estimators which can be
  postcomputed only (i.e., after the discrete solution has been
  computed).  To clarify this point, let us focus on the particular
  situation of interest.  Let $\uno$ be the solution from the previous
  timestep; $\uno\in\V{n-1}$, the finite element space with respect to
  mesh $\T{n-1}$.  The error due to coarsening appears in the term
  \begin{equation}
    \uno-\laginterpol n\uno.
  \end{equation}
  This term is nonzero only when simplexes are coarsened.  
  
  Furthermore, we assume that the new mesh $\T n$ is a refinement of
  $\coarseT n$, which is a coarsening of the old mesh $\T{n-1}$:
  \begin{equation}
    \label{eqn:diag:coarsen-refine}
    \T{n-1} 
    \xrightarrow{\text{coarsen}}
    \coarseT n 
    \xrightarrow{\text{refine}}
    \dotsb 
    \xrightarrow{\text{refine}}
    \T n
  \end{equation}
  
  If $\coarseinterpol n$ is the Lagrange interpolator onto the finite
  element space $\V n_0$, relative to the new coarse mesh $\coarseT n$,
  it is not very difficult to predict $\coarseinterpol n\uno$ without
  actually computing it.  Therefore this term can be predicted from (a)
  the simplexes of $\T{n-1}$ marked for coarsening which leads to
  $\coarseT n$ and (b) the values of $\uno$.
  
  Note that since $\coarseT n$ is subsequently \emph{refined but not
    coarsened} to produce $\T n$, as depicted in
  \eqref{eqn:diag:coarsen-refine}, then the additional coarsening error
  will be zero.  Namely, if $\laginterpol n$ denotes the Lagrange
  interpolant onto $\V n$, the finite element space over $\T n$, which
  is a refinement of $\coarseT n$, then $\laginterpol
  n\uno=\coarseinterpol n\uno$, and thus
  \begin{equation}
    \uno-\laginterpol n\uno=\uno-\coarseinterpol n\uno.
  \end{equation}
  The coarsening strategy therefore consists in choosing a subset of
  simplexes of $\T{n-1}$ which minimises term
  $\Norm{\uno-\coarseinterpol n\uno}$ \emph{before} producing the new
  coarse mesh $\coarseT n$.
  
  The rest of this section describes how $\uno-\coarseinterpol n\uno$
  can be precomputed.
  \subsection{Notation}
  \label{prelim}
  Let $K$ be an element of the new coarse mesh $\coarseT n$ resulting
  from the coarsening of its two children which we denote by $K^\pm$.
  (Note that $K^+$ and $K^-$ correspond to $\codename{child[0]}$ and
  $\codename{child[1]}$ of $K$ in the \alberta manual~\cite{alberta}.)
  Define the \emph{fine space}
  \begin{equation}
    \fine{}:=\ensemble{\restriction\Fi K}{\Fi\in\V{n-1}}.
  \end{equation}
  Likewise define
  the \emph{coarse space} $\coarse$ to be the local finite
  element space, i.e.,
  \begin{equation}
    \coarse:=\ensemble{\restriction\Fi K}{\Fi\in\V n_0};
  \end{equation}
  simply put we just have $\coarse=\poly p$.  We introduce also the
  \emph{fine spaces} $\fine{\pm}$, defined like $\fine{}$, but
  restricting functions over $K^\pm$, respectively (so functions in
  $\fine\pm$ are in fact the same as $\coarse=\poly p$, albeit with
  different domains).
  
  Denote by $\setof{\geovec x_0,\dotsc,\geovec x_L}$ and $\setof{\geovec
    x_0^{\pm},\dotsc,\geovec x_L^\pm}$ the set of Lagrange degrees of
  freedom on the simplex $K$ and its children $K^\pm$, respectively. We
  indicate with $\setof{\pi^0,\dotsc,\pi^L}$ and
  $\setof{\pi_\pm^0,\dotsc,\pi_\pm^L}$ the corresponding Lagrange
  polynomial bases of $\coarse{}$ and $\fine{\pm}$, respectively,
  whereby
  \begin{equation}
    \label{eqn:lagrange-basis:property}
    \pi^i(\geovec x_j)=\pi_{\pm}^i(\geovec x_j^\pm)=\delta^i_j.
  \end{equation}
  For short we will write these bases as column vectors
  $\vec{\pi}=\colvec{\pi^0,\dotsc,\pi^L}$, etc.  We also define the
  (local) \emph{coarse-on-fine matrixes} by
  \begin{equation}
    \label{eqn:coarse-on-fine-matrixes}
    \mat A^\pm:=\vecof{\vec\pi(\geovec x_0^\pm)\ \dotsc\ \vec\pi(\geovec x_L^\pm)}
    =\matof{\pi^i(\geovec x_j^\pm)}_{i,j=0,\dotsc,L}.
  \end{equation}
  These matrixes are closely related to \alberta's \emph{refine-interpolation}
  matrix ~\cite[matrix $A$ (1.5) in \S1.4.4 ]{alberta}.
  \begin{Pro}[coarse-on-fine matrix properties]
    \label{pro:coarse-on-fine-matrixes}
    The matrixes $\mat A^+$ and $\mat A^-$ are independent of
    $K,K^+,K^-$ and
    \begin{equation}
      \restriction{\vec\pi}{K^\pm} 
      = \mat A^\pm\vec\pi_\pm.
    \end{equation}
  \end{Pro}
  \begin{Proof}
    Fix $\rangefromto i0L$. Because $\pi^i$ is a polynomial and
    $\setof{\pi_+^0,\dotsc,\pi_+^L}$ is a polynomial basis, it follows that
    \begin{equation}
      \pi^i=\sumifromto j0L a^i_j\pi_+^j,
    \end{equation}
    for some vector $\vecof{a^i_0,\dotsc,a^i_L}$.  Applying $\pi^i$ to
    $\geovec x_j^+$, and recalling (\ref{eqn:lagrange-basis:property}), we
    obtain
    \begin{equation}
      a^i_j=\pi^i(\geovec x_j^+),
    \end{equation}
    and hence
    \begin{equation}
      \pi^i=\qb{\mat A^+\vec\pi_+}^i.
    \end{equation}
  \end{Proof}
  \begin{Example}[quadratic elements in $2$ dimensions]
    \label{Ex:quad2d}
    To make the discussion more accessible, we will illustrate it as we go
    with the concrete situation where $p=2$ (quadratic elements) and
    $d=2$.  Following the \alberta conventions the relation between the
    coarse and fine triangles is given by the following diagram.
    
    \begin{tikzpicture}      
      
      \path[fill=green!60] 
      (0,0) node (N0) {}
      -- (2,0) node (N5) {}
      -- (2,2) node (N2) {}
      -- cycle;
      
      \path[coordinate]
      (2,.67) node (K) {};
      
      \path[fill=green!60,xshift=6cm] 
      (0,0) node (NL1) {}
      -- (2,0) node (NL2) {}
      -- (2,2) node (NL0) {}
      -- cycle;
      
      \foreach \s in {0,1}{
	\path[coordinate,xshift=6cm*\s]
	(1.33,.67) node (L\s) {};
      }
      
      \path[fill=green!20]
      (4,0) node (N1) {}
      -- (2,2)
      -- (2,0)
      -- cycle
      ;
      
      \path[fill=green!20,xshift=7cm]
      (4,0) node (NR0) {}
      -- (2,2) node (NR1) {}
      -- (2,0) node (NR2) {}
      -- cycle
      ;
      
      \foreach \s in {0,1}{
	\path[coordinate,xshift=7cm*\s]
	(2.67,.67) node (R\s) {};
      }
      
      \path[coordinate]
      (3,1) node (N3) {}
      (1,1) node (N4) {}
      ;
      
      \path[coordinate,xshift=6cm]
      (1,0) node (NL3) {}
      (2,1) node (NL4) {}
      (1,1) node (NL5) {}
      ;

      \path[coordinate,xshift=7cm]
      (2,1) node (NR3) {}
      (3,0) node (NR4) {}
      (3,1) node (NR5) {}
      ;      
      
      \foreach \s in {0,1,...,5}{
	\node (K\s) at (N\s) [circle,inner sep=1pt,fill=blue!30,draw] {$\s$};
      }
      
      \foreach \s in {0,1,2,5}{
	\node (KL\s) at (NL\s) [circle,inner sep=1pt,fill=blue!30,draw] {$\s$};
	\node (KR\s) at (NR\s) [circle,inner sep=1pt,fill=blue!30,draw] {$\s$};
      }
      
      \foreach \s in {3,4}{
	\node (KL\s) at (NL\s) [circle,inner sep=1pt,fill=yellow!60,draw] {$\s$};
	\node (KR\s) at (NR\s) [circle,inner sep=1pt,fill=yellow!60,draw] {$\s$};
      }
      
      \node at (K) [rectangle,inner sep=1pt] {$K$};
      
      \foreach \s in {0,1}{
	\node at (L\s) [rectangle,inner sep=1pt] {$K^+$};
	\node at (R\s) [rectangle,inner sep=1pt] {$K^-$};
      }
      
      
      \draw[->] (4,0.9) -- (6,0.9)
      node[pos=0.5,below] {refine$\,K$};
      
      \draw[<-] (3.6,1.3) -- (6.4,1.3)
      node[pos=0.5,above] {coarsen$(K^+,K^-)$};
      
    \end{tikzpicture}
    
    In this case, the coarse-on-fine matrixes are computed as follows:
    \begin{equation}
      \begin{split}
	\vec A^+
	=
	\begin{bmatrix}
	  {} &  1 & {} & \phantom+3/8 &         -1/8 & {}\\
	  {} & {} & {} &         -1/8 &         -1/8 & {}\\
	  1 & {} & {} &           {} &           {} & {}\\
	  {} & {} & {} &           {} & \phantom+1/2 & {}\\
	  {} & {} & {} &           {} & \phantom+1/2 &  1\\
	  {} & {} &  1 & \phantom+3/4 & \phantom+1/4 & {}
	\end{bmatrix}
	,\:
	\vec A^-
	=
	\begin{bmatrix}
      	  {} & {} & {} &         -1/8 &         -1/8 & {}\\
	  1 & {} & {} &         -1/8 & \phantom+3/8 & {}\\
	  {} &  1 & {} &           {} &           {} & {}\\
	  {} & {} & {} & \phantom+1/2 &           {} &  1\\
	  {} & {} & {} & \phantom+1/2 &           {} & {}\\
	  {} & {} &  1 & \phantom+1/4 & \phantom+3/4 & {}
	\end{bmatrix}
      \end{split}
    \end{equation}
  \end{Example}
  \subsection{Degrees of freedom and global--local relations}
  \label{sec:dof-glob-loc}
  Denote by $U$ the generic finite element function in the old space
  $\V{n-1}$ and let $V:=\coarseinterpol n U$.  Then we have
  \begin{equation}
    U
    =\transposevec u\vec\Ps
    \AND
    V
    =\transposevec v\vec\Fi,
  \end{equation}
  where $\vec\Ps=\colvec{\Ps^0,\dotsc,\Ps^N}$ and
  $\vec\Fi=\colvec{\Fi^0,\dotsc,\Fi^M}$, are the columns of nodal
  Lagrange piecewise polynomial bases of $\V{n-1}$ and $\V n_0$,
  respectively, and $\vec u$ and $\vec v$ are the corresponding
  vectors of DOF values.
  
  There are $L+1$ degrees of freedom (DOF) per simplex, e.g., $L=5$
  for $p=2=d$.  The simplex $K$ in $\T n_0$ comes
  with a \emph{local-to-global} index relation
  $\funk{g=g^{\T n_0}_K}{\fromto0L}{\fromto0M}$ whereby
  \begin{equation}
    \restriction{\Fi^{g(i)}}K=\pi^i\Foreach\rangefromto j0L.
  \end{equation}
  It follows that the finite element function $V$ is locally
  represented on $K$ by
  \begin{equation}
    Y:=\restriction V{K}=\sumifromto i0L v_{g(i)}\pi^i=:\transposevec y\vec\pi.
  \end{equation}
  Similarly we have $\funk{g^\pm=g^{\T{n-1}}_{K^\pm}}{\fromto0L}{\fromto0N}$
  such that
  \begin{equation}
    Y^\pm:=\restriction U{K^\pm}
    =\sumifromto j0L u_{g^\pm(j)}\pi^j_{\pm}
    =:\transposevec{y^\pm}\vec\pi_\pm.
  \end{equation}
  The relation between the DOF coefficients $\vec u$ and $\vec v$ will
  be described next.
  \subsection{Local fine--coarse DOF relations}
  \label{sec:local-fine-coarse}
  Some degrees of freedom---that is those depicted in yellow or
  bright---are removed during coarsening.  The others, which are kept,
  have their local index change.  This information is fully encoded in
  the \emph{fine-to-coarse} index maps $\funk{c^\pm}{D^\pm}{C^\pm}$
  where
  \begin{equation}
    D^\pm:=
    \ensemble{\rangefromto j0L}
	     {\geovec x_j^\pm\in\setof{\geovec x_0,\dotsc,\geovec x_L}}.
  \end{equation}
  and
  \begin{equation}
    C^\pm:=c^\pm(D^\pm)\subseteq\fromto0L.
  \end{equation}
  A basic property of the fine-to-coarse maps is that
  \begin{equation}
    C^+\join C^-=\fromto0L,
  \end{equation}
  but $C^+$ and $C^-$ need not be disjoint (in fact, for conforming
  methods these are never disjoint).  The fine-to-coarse maps $c^\pm$
  are injective and we denote their inverses, the \emph{coarse-to-fine} maps,
  by $\funk{d^\pm}{C^\pm}{D^\pm}$.
  
  In the example above, $p=2=d$, the fine-to-coarse maps
  $\funk{c^\pm}{D^\pm}{\fromto05}$, satisfy $D^+=D^-=\setof{0,1,2,5}$
  (though $D^+$ and $D^-$ do not generally coincide, as seen for
  $p=3,\,d=2$, e.g.) and evaluated by the schedule
  \begin{equation}
    \begin{array}{r c c c c c c}
      j=& 0 & 1 & 2 & 3 & 4 & 5,\\
      c^+(j)=& 2 & 0 & 5 & - & - & 4,\\
      c^-(j)=& 1 & 2 & 5 & - & - & 3.
    \end{array}
  \end{equation}
  It follows that $C^+=\setof{0,2,4,5}$ and $C^-=\setof{1,2,3,5}$ and
  \begin{equation}
    \begin{array}{r c c c c c c}
      i     =& 0 & 1 & 2 & 3 & 4 & 5,\\
      d^+(i)=& 1 & - & 0 & - & 5 & 2,\\
      d^-(i)=& - & 0 & 1 & 5 & - & 2.
    \end{array}
  \end{equation}
  \begin{Obs}[redundancy of the coarse-to-fine maps]
    The coarse-to-fine maps $c^\pm$ and their inverses $d^\pm$ are
    partially redundant with $\mat A^\pm$.  Namely, if $j\in D^\pm$,
    then $j=d^\pm(i)$ and $i=c^\pm(j)$, for some $\rangefromto i0L$.
    By definition of $c^\pm$ it follows that $\geovec x^\pm_j=\geovec
    x_i$.  Therefore
    \begin{equation}
      \getrowcol{\mat A^\pm}k j
      =\pi^k(\geovec x^\pm_j)
      =\pi^k(\geovec x_i)
      =\delta^k_i.
    \end{equation}
    We have thus proved the following result that will be used to
    compress $\mat A^\pm$ in the sequel.
  \end{Obs}
  \begin{Pro}[redundant coarse-on-fine columns]
    \label{pro:redundant-columns}
    If $j\in D^\pm$, then $\mat A^\pm$'s $j$-th column is described by
    \begin{equation}
      \getrowcol{\mat A^\pm}k j=\delta^k_{c^\pm(j)}.
    \end{equation}
  \end{Pro}
  \subsection{Precomputing the coarsening error}
  \label{sec:precompute-coarse-error}
  The coarsening error is the difference between $U$, to which we have
  access via $\vec u$, and its interpolation on the locally coarser
  mesh $V$, to which we have no direct access.  Working locally at the
  coarsening-marked element $K^+$ (and similarly for $K^-$), all we
  need is to compute $\restriction V{K^+}$ and subtract it from
  $\restriction U{K^+}$.
  
  Recalling that in \alberta $V=\laginterpol n_0U$ is built by simply
  ``dropping'' the coefficients of the DOF removed by coarsening we have
  \begin{equation}
    \transposevec y\vec\pi=Y=\restriction V{K}
    =\sumsu{u_{g_+(d^+(i))}\pi}i{C^+}
    +\sumsu{u_{g_-(d^-(i))}\pi}i{C^-\take C^+},
  \end{equation}
  that is, for $\rangefromto j0L$, we set
  \begin{equation}
    \label{eqn:coarse-dof-calculation}
    v_{g(i)}:= y_i :=
    \begin{cases}
      u_{g_+(d^+(i))}=y^+_{d_+(i)}&\text{ if }i\in C^+
      \\
      u_{g_-(d^-(i))}=y^-_{d_-(i)}&\text{ otherwise }.
    \end{cases}
  \end{equation}
  (Note that the vector $\vec y$ is the same for the two siblings
  $K^\pm$ and needs to be calculated only once.)  Following the
  example with $p=2=d$, we see that
  \begin{equation}
    \begin{aligned}
      \vec y
      &=\colvec{y_1^+,y_0^-,y_0^+,y_5^-,y_5^+,y_2^+}
      \\
      &=\colvec{y_1^+,y_0^-,y_1^-,y_5^-,y_5^+,y_2^-}.
    \end{aligned}
  \end{equation}
  
  To conclude we rewrite the coarse basis, $\vec\pi$, in
  terms of the fine one, $\vec\pi_+$, using Proposition
  \ref{pro:coarse-on-fine-matrixes} as follows:
  \begin{equation}
    \restriction V{K^+}=\restriction Y{K^+}
    =\transposevec y\restriction{\vec\pi}{K^+}
    =\transposevec y\mat A^+\vec\pi_+.
  \end{equation}
  Thus the coarsening error on $K^+$ is calculated as
  \begin{equation}
    \restriction{\qb{U-V}}{K^+}
    =\Transpose{\vec y^+}\vec\pi_+
    -\transposevec y\mat A^+\vec\pi_+
    =\Transpose{\vec\pi_+}
    \qp{\vec y^+-\Transpose{\mat A^+}\vec y}
    =\sumifromto j0L
    \qp{y^+_j-\transposevec y\getcol{\mat A^+}j}\pi^+_j.
  \end{equation}
  Recalling Proposition \ref{pro:redundant-columns}, if $j\in D^+$ we have
  \begin{equation}
    \transposevec y\getcol{\mat A^+}j
    =
    \sumifromto i0L y_k\delta^k_{c^+(j)}
    =y_{c^+(j)}=y^+_j,
  \end{equation}
  and thus the coefficient for $\pi^+_j$ is $0$, and it needs not be
  calculated.  Proceeding similarly on $K^-$ we may summarise the
  findings as follows.
  \begin{The}[coarsening error calculation]
    Let $U\in\V{n-1}$ with the notation of \S\ref{sec:dof-glob-loc},
    to calculate the coarsening error that would result from
    coarsening the elements $K^+,K^-\in\T{n-1}$ into $K\in\T{n}$
    \begin{enumerate}[(1)\ ]
    \item
      calculate $\vec y$ following (\ref{eqn:coarse-dof-calculation})
      using the coarse-to-fine map $d^+$ defined in \secref{sec:local-fine-coarse},
    \item
      obtain the error using
      \begin{equation}
	\begin{gathered}
	  \restriction{\qb{U-\laginterpol n_0U}}{K^+}
	  =\sum_{j\in\fromto0L\take D^+}\qp{y^+_j-\transposevec y\getcol{\mat A^+}j}\pi_+^j,
	  \\
	  \restriction{\qb{U-\laginterpol n_0U}}{K^-}
	  =\sum_{j\in\fromto0L\take D^-}\qp{y^-_j-\transposevec y\getcol{\mat A^-}j}\pi_-^j.
	\end{gathered}
      \end{equation}
    \end{enumerate}
  \end{The}
  \begin{Obs}
    Note that the $j$-th coefficient of the coarsening error's local
    DOF vector is zero when $j\in D^\pm$, respectively.  So the
    calculation needs to be carried out only for those $j\not\in D^\pm$.
    
    Also, the coefficients for the DOF that are common to $K^+$ and
    $K^-$ must be equal, so they can be in fact computed once.
  \end{Obs}
  For example in the case of quadratic elements in $d=2$ we have
  \begin{equation}
    \begin{split}
      Y^+ - \restriction Y{K^+} 
      =&
      \pi^3_+ 
      \qp{y^+_3 -\frac38y^+_1 +\frac18y^-_0 -\frac34y_2^+ }
      \\
      +&\pi_+^4
      \qp{y^+_4 +\frac18y^+_1 +\frac18y^-_0 -\frac14 y^+_2  
	-\frac12y^+_5  -\frac12y^-_5},
      \\
      Y^- -\restriction Y{K^-}
      =&
      \pi_-^3
      \qp{y^-_3 + \frac{1}{8} y^+_1 + \frac{1}{8} y^-_0
	- \frac{1}{2} y^-_5 -\frac{1}{2} y^+_5 - \frac{1}{4} y^+_2}
      \\
      +& \pi^4_-
      \qp{y^-_4 + \frac{1}{8} y^+_1 - \frac{3}{8} y^-_0 
	- \frac{3}{4} y^+_2}  
    \end{split}
  \end{equation}
  \subsection{Coarsening error algorithm}
  As seen in \secref{sec:precompute-coarse-error}, the information
  needed for the coarsening error computation for Lagrange finite
  elements of degree $p$ in dimension $d$, is contained in the
  coarse-on-fine matrixes $\vec A^\pm$ defined by
  (\ref{eqn:coarse-on-fine-matrixes}) and the fine-to-coarse maps,
  $d^\pm$, and their domains $C^\pm$ defined in
  \ref{sec:local-fine-coarse}. This information is independent of the
  particular pair of simplex siblings $K^\pm$ and their parent $K$ and
  can be included in the code via given index permutations and efficient
  matrix-vector multiplication.
  
  With this information at hand and the notation previously introduced
  in this section, we formulate an \alberta-implementable algorithm to
  precompute the coarsening error on all seimplexes.
  \subsection*{\Algoname{Coarsening Preindicator}}
  \label{alg:coarsening-error}
  \begin{algorithmic}
    \Require $(\U{}=\transnumvec u\numvec\Phi, \fespace, \T{})$ 
    
    \Ensure $\vec\gamma=(\gamma_K:K\in\T{})$
    
    \ForAll{$K\in\T{}$}
    
    \If{$\operatorname{child order}(K)=0$\footnote{The element information in
	\alberta is quite local and to determine whether an element is
	left or right child is not trivial. In \alberta 1.2 this can
	be done utilising \codename{EL->index} which provides a global
	indexing of elements. Testing the
	\codename{EL\_INFO->parent->child[0]->index} against
	\codename{EL->index} gives the correct child order of $K$. In
	\alberta 2.0 \codename{EL->index} is unavailable so we check
	the global index of DOF for both parent and children.} }
    
    \State $D\assignvalue D^+$, $D'\assignvalue D^-$, $c\assignvalue c^+$, $c'\assignvalue c^-$, $\mat A\assignvalue\mat A^+$
    
    \Else
    
    \State $D\assignvalue D^-$, $D'\assignvalue D^+$, $c\assignvalue c^-$, $c'\assignvalue c^+$, $\mat A\assignvalue\mat A^-$
    
    \EndIf
 
    \State $K'\assignvalue\operatorname{sibling} K$
       
    \State initialise two local DOF vectors $\numvec y$ and $\numvec r$

    \ForAll{$j\in D$}
    
    \State $y_{c(j)} = u_{g_K(j)}$
    
    \EndFor
    
     \ForAll{$j\in D'$}
    
    \State $y_{c'(j)} = u_{g_{K'}(j)}$
    
    \EndFor
    
    \ForAll{$j\notin D \join D'$}
    
    \State $r_j=u_{g_K(j)}-\transnumvec y\getcol{\nummat A}j$
    
    \EndFor
    
    \State $\gamma_K = 0$
    
    \ForAll{$i \notin D \cup D'$}
    
    \ForAll{$j \notin D \cup D'$}
    
    \State $\gamma_K =\gamma_K+ r_ir_j\ltwop{\Fi_i}{\Fi_j}_K$
    
    \EndFor
    
    \EndFor
    
    \EndFor
  \end{algorithmic}
  \subsection{Coarsening preindicator matrixes}
  To close, we provide here the information needed to implement Algorithm
  \ref{alg:coarsening-error} for Lagrange piecewise $\poly p$ finite
  elements in dimension $d=2$.  (For dimension $3$ the situation is
  complicated by the ``types'' of tetrahedrons, whereby the matrixes
  $A^\pm$ and the maps $c^\pm$ may depend on the type and is not covered
  in this appendix.)
  \subsection{$\poly1$ elements}
  The coarse-on-fine matrixes (omitting $0$ entries for clarity) are given by
  \begin{equation}
    \vec{A^+} =
    \begin{bmatrix}
      {} &  1  & 1/2\\
      {} & {}  & 1/2\\
      1 & {}  & {} 
    \end{bmatrix}
    ,
    \vec{A^-} =
    \begin{bmatrix}
      {} & {} & 1/2\\
      1 & {} & 1/2\\
      {} &  1 & {}
    \end{bmatrix}
    ,
  \end{equation}
  the fine-to-coarse maps and the coarse-to-fine maps are respectively given by
  \begin{equation}
    \begin{array}{rccc}
      i=&0&1&2,
      \\
      c^+(i)=&2&0&-,
      \\
      c^-(i)=&1&2&-,
    \end{array}
    \AND
    \begin{array}{rccc}
      i=&0&1&2,
      \\
      d^+(i)=&1&-&0,
      \\
      d^-(i)=&-&0&1.
    \end{array}
  \end{equation}
  \subsection{$\poly2$ elements}
  See the worked example in \secref{sec:building-coarsening-estimator}.
  \subsection{$\poly3$ elements}
  The coarse-on-fine matrixes are given by
  \begin{align}
    \mat A^+ 
    &=
    \begin{bmatrix}
      {} &  1 &          -1/16 &  \phantom+5/16 & {} & {} & \phantom+1/16         & {} & {} &         -1/16         \\
      {} & {} &          -1/16 &  \phantom+1/16 & {} & {} & \phantom+1/16         & {} & {} & \phantom+1/16         \\
      1 & {} &             {} &             {} & {} & {} &            {}         & {} & {} &           {}          \\
      {} & {} &             {} &             {} & {} & {} &         -1/4\phantom0 & {} & {} &         -1/8\phantom0 \\
      {} & {} &             {} &             {} & {} & {} & \phantom+1/2\phantom0 & {} & {} &           {}          \\
      {} & {} &             {} &             {} & {} & {} & \phantom+1/2\phantom0 &  1 & {} &           {}          \\
      {} & {} &             {} &             {} & {} & {} &         -1/4\phantom0 & {} &  1 & \phantom+3/8\phantom0 \\
      {} & {} &  \phantom+9/16 &          15/16 &  1 & {} &         -1/16         & {} & {} & \phantom+3/16         \\
      {} & {} &  \phantom+9/16 &          -5/16 & {} & {} &         -1/16         & {} & {} &         -3/16         \\
      {} & {} &             {} &             {} & {} &  1 & \phantom+1/2\phantom0 & {} & {} & \phantom+3/4\phantom0 \\
    \end{bmatrix}
    \intertext{and}
    \mat A^-
    &=
    \begin{bmatrix}
      {} & {} &          -1/16 & \phantom+1/16         & {} &  \phantom+1/16 & {} & {} & {} &         -1/16         \\
      1 & {} &          -1/16 & \phantom+1/16         & {} &  \phantom+5/16 & {} & {} & {} &           {}          \\
      {} &  1 &           {}   &            {}         & {} &           {}   & {} & {} & {} & \phantom+1/16         \\
      {} & {} &           {}   &         -1/4\phantom0 & {} &           {}   & {} &  1 & {} &           {}          \\
      {} & {} &           {}   & \phantom+1/2\phantom0 & {} &           {}   & {} & {} &  1 &         -1/8\phantom0 \\
      {} & {} &           {}   & \phantom+1/2\phantom0 & {} &           {}   & {} & {} & {} &         -3/16         \\
      {} & {} &           {}   &         -1/4\phantom0 & {} &           {}   & {} & {} & {} & \phantom+3/16         \\
      {} & {} & \phantom+9/16 &          -1/16         & {} &          -5/16 & {} & {} & {} & \phantom+3/8\phantom0 \\
      {} & {} & \phantom+9/16 &          -1/16         & {} &          15/16 &  1 & {} & {} &           {}          \\
      {} & {} &           {}   & \phantom+1/2\phantom0 &  1 &           {}   & {} & {} & {} & \phantom+3/4\phantom0
    \end{bmatrix}
  \end{align}
  the fine-to-coarse maps
  \begin{equation}
    \begin{array}{rcccccccccc}
      i=&0&1&2&3&4&5&6&7&8&9,\\
      c^+(i)=&2&0&-&-&7&9&-&5&6&-,\\
      c^-(i)=&1&2&-&-&9&8&-&3&4&-.
    \end{array}
  \end{equation}
  and the coarse-to-fine maps
  \begin{equation}
    \begin{array}{rcccccccccc}
      i=&0&1&2&3&4&5&6&7&8&9,\\
      d^+(i)=&1&-&0&-&-&7&8&4&-&5,\\
      d^-(i)=&-&0&1&7&8&-&-&-&5&4.
    \end{array}
  \end{equation}
  \newpage
  \subsection{$\poly4$ elements}
  The coarse-on-fine matrixes are given by
  \\
  \begin{sideways}{
      \small
      \begin{minipage}{1.4\linewidth}
	\begin{align*}
	  \vec{A^+} &=
	  \left[
	    \begin{array}{c c c c c c c c c c c c c c c}
	      {}    & 1  & {} & 35/128 & {} &         -5/128 & \phantom+3/128 & {} & -5/128 & {} & {} & {} & \phantom+3/128 & -5/128 & {}
	      \\ {} & {} & {} & -5/128 & {} & \phantom+3/128 & \phantom+3/128 & {} & -5/128 & {} & {} & {} & -5/128 & -5/128 & {}
	      \\  1 & {} & {} & {}     & {} & {} & {} & {} & {} & {} & {} & {} & {} & {} & {}
	      \\ {} & {} & {} & {}     & {} & {} & -1/16 & {} & 3/16 & {} & {} & {} & \phantom+1/8 & \phantom+1/16 & {}
	      \\ {} & {} & {} & {}     & {} & {} & {} & {} & -3/8 & {} & {} & {} &-1/8 & {} & {}
	      \\ {} & {} & {} & {}     & {} & {} & {} & {} & \phantom+1/2 & {} & {} & {} & {} & {} & {}
	      \\ {} & {} & {} & {}     & {} & {} & {} & {} & \phantom+1/2 & 1 & {} & {} & {} & {} & {}
	      \\ {} & {} & {} & {}     & {} & {} & {} & {} & -3/8 & {} & 1 & {} & \phantom+3/8 & {} & {}
	      \\ {} & {} & {} & {}     & {} & {} & -1/16 & {} & \phantom+3/16 & {} & {} & 1 & -1/8 & \phantom+5/16 & {}
	      \\ {} & {} & {} & \phantom+35/32 & 1 & \phantom+15/32 & -3/32 & {} & \phantom+1/32 & {} & {} & {} & -1/32 & \phantom+5/32 & {}
	      \\ {} & {} & 1 & -35/64  & {} & \phantom+45/64 & \phantom+9/64  & {} & \phantom+1/64 & {} & {} & {} & \phantom+3/64 & \phantom+15/64
	      \\ {} & {} & {} & \phantom+7/32 & {} & -5/32 & -3/32 & {} & \phantom+1/32 & {} & {} & {} & \phantom+3/32 & \phantom+5/32 & {} 
	      \\ {} & {} & {} & {}     & {} & {} & \phantom+9/16 & {} & -3/16 & {} & {} & {} & \phantom+3/8 & \phantom+15/16 & 1
	      \\ {} & {} & {} & {}     & {} & {} & \phantom+9/16 & {} & -3/16 & {} & {} & {} & -3/8 & -5/16 & {}
	      \\ {} & {} & {} & {}     & {} & {} & {} & 1 & \phantom+3/4 & {} & {} & {} & 3/4 & {} & {}
	    \end{array}
	    \right
	  ]
	  \\
	  \vec{A^-} &=
	  \left[
	    \begin{array}{c c c c c c c c c c c c c c c}
	      {} & {} & {} & -5/128 & {} & \phantom+3/128 & \phantom+3/128 & {} & -5/128 & {} & {} & {} & -5/128 & -5/128 & {}
	      \\ 1 & {} & {} & -5/128 & {} & \phantom+3/128 & -5/128 & {} & \phantom+35/128 & {} & {} & {} & -5/128 & \phantom+3/128 & {}
	      \\ {} & 1 & {} & {} & {} & {} & {} & {} & {} & {} & {} & {} & {} & {} & {}
	      \\ {} & {} & {} & \phantom+3/16 & {} & -1/16 & {} & {} & {} & 1 & {} & {} & \phantom+5/16 & -1/8 & {}
	      \\ {} & {} & {} & -3/8 & {} & {} & {} & {} & {} & {} & 1 & {} & {} & 3/8 & {}
	      \\ {} & {} & {} & \phantom+1/2 & {} & {} & {} & {} & {} & {} & {} & 1 & {} & {} & {}
	      \\ {} & {} & {} & \phantom+1/2 & {} & {} & {} & {} & {} & {} & {} & {} & {} & {} & {}
	      \\ {} & {} & {} & -3/8 & {} & {} & {} & {} & {} & {} & {} & {} & {} & -1/8 & {}
	      \\ {} & {} & {} & \phantom+3/16 & {} & -1/16 & {} & {} & {} & {} & {} & {} & \phantom+1/16 & \phantom+1/8 & {}
	      \\ {} & {} & {} & \phantom+1/32 & {} & -3/32 & -5/32 & {} & \phantom+7/32 & {} & {} & {} & \phantom+5/32 & \phantom+3/32 & {}
	      \\ {} & {} & 1 & \phantom+1/64 & {} & \phantom+9/64 & \phantom+45/64 & {} & -35/64 & {} & {} & {} & -15/64 & -3/64 & {}
	      \\ {} & {} & {} & \phantom+1/32 & {} & -3/32 & \phantom+15/32 & 1 & \phantom+35/32 & {} & {} & {} & \phantom+5/32 & -1/32 & {}
	      \\ {} & {} & {} & -3/16 & {} & \phantom+9/16 & {} & {} & {} & {} & {} & {} & -5/16 & -3/8 & {}
	      \\ {} & {} & {} & -3/16 & {} & 9/16 & {} & {} & {} & {} & {} & {} & \phantom+15/16 & \phantom+3/8 & 1
	      \\ {} & {} & {} & \phantom+3/4 & 1 & {} & {} & {} & {} & {} & {} & {} & {} & \phantom+3/4 & {}
	    \end{array}
	    \right
	  ]. 
	\end{align*}
      \end{minipage}
    }
  \end{sideways}
  
  The fine-to-coarse maps are given by
  \begin{equation*}
    \begin{array}{rccccccccccccccc}
      i=&0&1&2&3&4&5&6&7&8&9&10&11&12&13&14,\\
      c^+(i)=&2&0&10&-&9&-&-&14&-&6&7&8&-&-&12,\\
      c^-(i)=&1&2&10&-&14&-&-&11&-&3&4&5&-&-&13.
    \end{array}
  \end{equation*}
  and the coarse-to-fine maps by
  \begin{equation*}
    \begin{array}{rccccccccccccccc}
      i=&0&1&2&3&4&5&6&7&8&9&10&11&12&13&14,\\
      d^+(i)=&1&-&0&-&-&-&9&10&11&4&2&-&14&-&7,\\
      d^-(i)=&-&0&1&9&10&11&-&-&-&-&2&7&-&14&4.
    \end{array}
  \end{equation*}
\providecommand{\bysame}{\leavevmode\hbox to3em{\hrulefill}\thinspace}
\providecommand{\MR}{\relax\ifhmode\unskip\space\fi MR }
\providecommand{\MRhref}[2]{%
  \href{http://www.ams.org/mathscinet-getitem?mr=#1}{#2}
}
\providecommand{\href}[2]{#2}

\end{document}